\documentclass[preprint,11pt]{elsarticle}

\usepackage{amssymb}
\usepackage{amsmath}
\usepackage{multirow}
\usepackage{booktabs}
\begin{document}

\begin{frontmatter}

%% Title, authors and addresses

%% use the tnoteref command within \title for footnotes;
%% use the tnotetext command for theassociated footnote;
%% use the fnref command within \author or \address for footnotes;
%% use the fntext command for theassociated footnote;
%% use the corref command within \author for corresponding author footnotes;
%% use the cortext command for theassociated footnote;
%% use the ead command for the email address,
%% and the form \ead[url] for the home page:
%% \title{Title\tnoteref{label1}}
%% \tnotetext[label1]{}
%% \author{Name\corref{cor1}\fnref{label2}}
%% \ead{email address}
%% \ead[url]{home page}
%% \fntext[label2]{}
%% \cortext[cor1]{}
%% \affiliation{organization={},
%%             addressline={},
%%             city={},
%%             postcode={},
%%             state={},
%%             country={}}
%% \fntext[label3]{}

%\title{A COMPARISON OF SIX NUMERICAL METHODS FOR INTEGRATING A COMPARTMENTAL HODGKIN-HUXLEY TYPE MODEL}
\title{A comparison of six numerical methods for integrating a compartmental Hodgkin-Huxley type model}

%% use optional labels to link authors explicitly to addresses:
%% \author[label1,label2]{}
%% \affiliation[label1]{organization={},
%%             addressline={},
%%             city={},
%%             postcode={},
%%             state={},
%%             country={}}
%%
%% \affiliation[label2]{organization={},
%%             addressline={},
%%             city={},
%%             postcode={},
%%             state={},
%%             country={}}

\author[label1]{R. Park\corref{cor1}}
\ead{rpark@asu.edu}
\cortext[cor1]{Corresponding author}

\affiliation[label1]{organization={Home},%Department and Organization
            addressline={}, 
            city={Chandler},
            postcode={85226}, 
            state={AZ},
            country={U.S.A.}}

\begin{abstract}
We compare six numerical integrators' performance when simulating a regular spiking cortical neuron model whose 74-compartments
are equipped with eleven membrane ion channels and Calcium dynamics. Four methods are explicit and two are implicit; three are finite 
difference PDE methods, two are Runge-Kutta methods, and one an exponential time differencing method. Three methods are first-, two 
commonly considered second-, and one commonly considered fourth-order. Derivations show, and simulation data confirms, that 
Hodgkin-Huxley type cable equations render multiple order explicit RK methods as first-order methods. Illustrations compare accuracy, 
stability, variations of action potential phase and waveform statistics. Explicit methods were found unsuited for our model given their inability 
to control spiking waveform consistency up to 10 microseconds less than the step size for onset of instability.  While the backward-time central 
space method performed satisfactorily as a first order method for step sizes up to 80 microseconds, performance of the Hines-Crank-Nicolson 
method, our only true second order method, was unmatched for step sizes of 1-100 microseconds. 
\end{abstract}

%%Graphical abstract
%\begin{graphicalabstract}
%\includegraphics{grabs}
%\end{graphicalabstract}

%%Research highlights
%\begin{highlights}
%\item Seven integrator's accuracy and stability visualized, $1$ to $100\mu$sec at $1\mu$sec resolution   
%\item Explicit, multi-order RK methods rendered $1^{st}$-O time by Hodgkin-Huxley cable equations
%\item Hines-Crank-Nicolson decaying oscillations visualized as a function of step size
%\end{highlights}

\begin{keyword}
%% keywords here, in the form: keyword \sep keyword
Numerical Integration \sep Cable Equation \sep Compartmental Neural Model \sep Stability \sep Hines-Crank-Nicolson \sep Runge-Kutta Methods
%% PACS codes here, in the form: \PACS code \sep code

%% MSC codes here, in the form: \MSC code \sep code
%% or \MSC[2008] code \sep code (2000 is the default)

\end{keyword}

\end{frontmatter}

%% \linenumbers

%% main text
\section{Introduction}
\label{sec:theIntro}
While finite difference (FD), Runge-Kutta (RK) and exponential-time-differencing (ETD) methods have numerically integrated Hodgkin-Huxley (HH) 
type neural models for over fifty years, it could be argued a measure of uncertainty remains regarding the differences one can expect by choosing 
one numerical method or integration step-size over another. Such uncertainties include: How does simulation inaccuracy present? How do action 
potential (AP) waveforms change as integration step size increases? Do ETD methods obviate the need for implicit integration methods \cite{Borgers}? 
Does oscillation onset appear only on the rising edge of spiking membrane voltage waveforms \cite{Borgers}? If the Hines-Crank-Nicolson (HCN) 
method allows decaying oscillations as step sizes grow large, when and where do these oscillations present and to what extent 
\cite[ch.4]{CarnevaleHines}? We address these questions as they relate to integrating a branched, $74$ compartment, single cell, neuronal model 
over integration step sizes ranging from $1$--$100\mu$s. \\ \\
Rallpacks formalized the measurement of similar data to answer the question, ``How does one simulator's performance compare to others in 
terms of accuracy, model size, simulation speed, and core memory needed?''\cite{BBB_Rall}. In contrast this study graphically illustrates how 
performance varies for different numerical methods as integration step size increases. Rallpacks' first two benchmarks were passive models 
and the third added two active ion membrane channels. Each compartment in our study's model has eleven active membrane ion channels 
and includes calcium dynamics.\\ \\
Section \ref{sec:TestProbAndNumMethods} begins by defining our test problem and listing the numerical methods we test. Section \ref{sec:AccStab} 
reviews the predicted accuracy and stability of the integration methods evaluated in this study. Here we show that the HH partial differential cable 
equation renders multiple order explicit RK methods mere first-order methods. Section \ref{sec:AnalysisMethods} describes data analysis techniques. \\ \\
The literature suggests significant AP waveform structures and measurements include action-potential (AP) spikes and bursts \cite{Feld12, Gol06, 
Hage15, Kepecs00, Kuznet10, Miller81, Pinsky94, Shad98, T2005}, spiking multiplicities, amplitudes and spike rise times \cite{Varela2012}, 
after-depolarization-pulses (ADP) \cite{Chu06, Kepecs00, Metz07, Varela2012}, and power spectrums \cite{Bahr13}. Results presented in Section 
\ref{sec:Results} show numerical integration method and step size combinations capable of inducing inaccuracies able to either masquerade as or 
distort these same membrane voltage waveform measures. We show that while ETD method stability makes it an excellent choice when integrating 
HH equations for single compartment models, the first-order exponential Euler method is not suited for multi-compartment models. We then identify 
integration step sizes and locations where oscillations appear in AP cycle waveforms integrated by the HCN method as well as examine their magnitude 
and duration. Finally, Section \ref{sec:SumAndDiscuss} discusses our results and how they relate to other studies.
\section{Test Problem and Numerical Methods} 
\label{sec:TestProbAndNumMethods}
\subsection{Governing Equations} 
The membrane voltage of compartmental excitable nerve cell models with varying cross-section diameters depends on neighboring compartments' 
membrane voltages as well as their own membrane ion channel states as described by the partial differential cable equation \cite[ch.6]{DayanAbbott} 
\begin{equation}
\label{CableEquationRoot}
c_m\frac{\partial V}{\partial t} + \frac{(V-E_L)}{r_m} + \sum_i g_i(V-E_i) - \frac{1}{2ar_L}\frac{\partial}{\partial x}\bigg(a^2\frac{\partial V}{\partial x}\bigg)=0.
\end{equation}
Parameters of equation (\ref{CableEquationRoot}) are defined in Table \ref{tab:cableEqnParms}.
\begin{table}[h!]\small
\caption{Cable Equation Parameters}
\ \\
\centering
\setlength{\tabcolsep}{1.0mm}
\begin{tabular}{l l l}
\toprule
Symbol & Parameter & Units   \\
\midrule \\[-2.0ex]
$a$ & compartment radius & $m$ \\
$c_m$  & specific membrane capacitance  &  $F/m^2$ \\   
$C_m$  & compartment membrane capacitance  &  $F$ \\     
$E_i$ & membrane ion channel $i$ Nernst Potential & $volts$ \\
$E_L$ & membrane leakage current Reversal Potential & $volts$ \\
$g_i$ & membrane ion channel $i$ specific conductance & $S/m^2~=~(\Omega\cdot m^2)^{-1}$ \\
$r_L$  & specific axial resistance  &  $\Omega \cdot m$ \\  
$R_a$  & compartment axial resistance  &  $\Omega$ \\  
$r_m$  & specific membrane resistance  &  $\Omega \cdot m^2$ \\     
$V$ & membrane voltage & $volts$ \\
\bottomrule
\end{tabular}
\label{tab:cableEqnParms}
\end{table}
\\ %\\
Membrane ion channel $i$ specific conductance, 
\begin{equation}
\label{ionChanSpecCond}
g_i=\overline{g}_if_i(m_i,h_i),~\text{for}~i=1,2,...,11,
\end{equation}
is the product of channel specific concentration, $\overline{g}_i$ and the channel specific gate states function, $f_i(m_i,h_i)$, for each membrane 
ion channel. Activation, $m_i$, and inactivation, $h_i$, states of channel constituent gates are governed by   
\begin{equation}
\label{eqn:gateStateODE}
\frac{dy}{dt}=\frac{y_{\infty}(V)-y}{\tau_y(V)}~\text{for }y=m~\text{or}~h,
\end{equation}
where gate steady state, $y_{\infty}(V)$, and gate time constant, $\tau_y(V)$, are functions of membrane voltage. 
\subsection{Neural Model} 
\label{sec:TraubL23rsModel}
Our model is Traub's cortical L2/3 regular spiking pyramidal neuron whose parameters were defined in \cite{T2005}.  It is deterministic and 
without hyper/de-polarizing or synaptic currents.
\subsection{Numerical Methods} 
Numerical integration methods are listed in Table \ref{tab:schemeTruncError}. Forward-time central-space, backward-time central-space and 
Hines-Crank-Nicolson, are numerical PDE methods. Traditional ODE methods, exponential Euler and RK numerical integrators, were included 
to demonstrate how not accounting for a PDE's changing spatial derivative affects an integrator's truncation error and our ability to predict step 
size upper limits guaranteeing stability. HCN is the only method to integrate the membrane voltage and gate states in a staggered fashion, 
temporally centering the gate state integration between membrane voltage time-steps as described in \cite{CarnevaleHines,Hines}. 
\subsection{Simulation Specification}
Simulation parameters are listed in Table \ref{tab:simSpec}.
\begin{table}[h!]\small
\caption{Simulation Specification}
\ \\
\centering
\begin{tabular}{l | l}
\toprule
Duration & $3.0$ seconds ($\sim40$ AP cycles) \\
\hline \\[-2.5ex]
Step Sizes & $1,~2,~3,~...~,~99~~\mu$s (given stability)\\
\hline \\[-2.0ex]
Initial Values & Compartment membrane voltage, $V=-0.07$ volts\\
& Calcium ion concentration, [Ca$^{2+}]=0.0$ \\
& Gate states: steady-state, $y_{\infty}\big(V=-0.07\big)$ \\
\bottomrule
\end{tabular}
\label{tab:simSpec}
\end{table}
\section{Accuracy and Stability}
\label{sec:AccStab}
\subsection{Accuracy}
\label{subsec:Acc} 
Integration methods listed in Table \ref{tab:schemeTruncError} are intended for different contexts. At the same time there exist myriad examples 
of these methods being used to integrate generic HH models \cite{GEN, Borgers, CarnevaleHines, CooleyDodge66, MooreRamon74, T2005}. 
To be clear, we integrated ion channel states represented by equation (\ref{eqn:gateStateODE}) with ODE methods forward Euler, backward 
Euler and trapezoidal rule in place of PDE methods forward-time central-space (FTCS), backward-time central-space (BTCS) and HCN, respectively. 
\begin{table}[h!]\small
\caption{Membrane Voltage Integration Truncation Error, $1^{st}$ and $2^{nd}$ Order Terms}
\ \\
\centering
\setlength{\tabcolsep}{1.0mm}
\begin{tabular}{l | l}
\toprule
Method & $P_{k,h}(\cdot)-P(\cdot)$  \\
\midrule \\[-2.0ex]
\parbox{2.5cm}{Forward-Time \\ Central-Space} 
& $\displaystyle ~~\frac{V_{tt}}{2}k + \frac{V_{ttt}}{6}k^2 -\gamma\bigg(\overline{a}^2V_{xxxx}+2\big(\overline{a}^2\big)_xV_{xxx}\bigg)\frac{h^2}{12}$ 
\\[2.0ex] \hline \\[-2.0ex]
\parbox{2.5cm}{Backward-Time Central-Space} 
& $\displaystyle -\frac{V_{tt}}{2}k-\frac{V_{ttt}}{3}k^2-\gamma\bigg(\overline{a}^2V_{xxxx}+2\big(\overline{a}^2\big)_xV_{xxx}\bigg)\frac{h^2}{12}$
\\ [1.0ex] \hline \\[-2.0ex]
\parbox{2.5cm}{Exponential \\ Euler} 
& $\displaystyle ~~\frac{V_{tt}}{2}k+\frac{V_{ttt}}{6}k^2-\gamma\bigg(\overline{a}^2V_{xxxx}+2\big(\overline{a}^2\big)_xV_{xxx}\bigg)\frac{h^2}{12}
+\frac{BV_t}{2}k - \frac{BV_t}{6}k^2$ \\[2.0ex] \hline \\[-2.0ex]
\parbox{2.5cm}{Hines-Crank-Nicolson} 
& $\displaystyle \qquad\qquad \frac{V_{ttt}}{24}k^2-\gamma\bigg(\overline{a}^2V_{xxxx}+2\big(\overline{a}^2\big)_xV_{xxx}\bigg)\frac{h^2}{12}$ 
\\ [1.0ex] \hline \\[-2.0ex]
\parbox{2.5cm}{RK21} 
& $\displaystyle ~~\frac{V_{tt}}{2}k + \frac{V_{ttt}}{6}k^2 -\gamma\bigg(\overline{a}^2V_{xxxx}+2\big(\overline{a}^2\big)_xV_{xxx}\bigg)\frac{h^2}{12} 
+\frac{BV_t}{2}k$\\ [1.0ex] \hline \\[-2.0ex]
\parbox{2.5cm}{RK41} 
& $\displaystyle ~~\frac{V_{tt}}{2}k + \frac{V_{ttt}}{6}k^2 -\gamma\bigg(\overline{a}^2V_{xxxx}+2\big(\overline{a}^2\big)_xV_{xxx}\bigg)\frac{h^2}{12}
+\frac{BV_t}{2}k$ \\ [1.0ex] \hline \\[-2.0ex]
Definitions & $\displaystyle \alpha=\frac{1}{r_mc_m}, \quad \beta_i=\frac{g_i}{c_m} \quad \text{and} \quad \gamma=\frac{1}{2ar_Lc_m},$ \\ \\[-1.0ex]
& $\displaystyle B=\Big(\alpha+\sum_i\beta_i\Big)+\bigg(\frac{1}{R_aC_m}+\frac{1}{R'_aC_m}\bigg)$ \\ \\[-1.0ex]  
& $\overline{a}:$ Average compartment radius, $~~ \big(\overline{a}^2\big)_x\equiv\big(a(x)^2\big)_x|_{x=\overline{a}}$ \\[1.0ex] 
\bottomrule
\end{tabular}
\label{tab:schemeTruncError}
\end{table}
\\ \\
For integrating equation (\ref{CableEquationRoot}) we normalized all our methods by recasting traditional ODE methods as quasi-FD methods to 
assess their accuracy and stability on an equal footing. Truncation errors are summarized in Table \ref{tab:schemeTruncError} not just to point out 
the order of a scheme's time and space error, but also to show differences between similarly ordered methods. Table \ref{tab:schemeTruncError} 
includes the claim that HH partial differential cable equations render multi-order explicit RK integrators mere first-order methods. The following 
derivation of RK21 truncation error describes the basis for this claim.
\subsubsection{RK21 Truncation Error}
\label{subsubsec:RK21truncErr} 
Using process and notation from \cite{Str} we define truncation error as the difference between the cable equation's differential operator and, 
what we define below, RK21's quasi-FD operator.  Modifying derivative notation in equation (\ref{CableEquationRoot}) lets us write the cable 
equation's differential operator as
\begin{multline}
\label{txtCableEquationOperator}
P(V) = V_t + \Big(\alpha+\sum_i\beta_i\Big)V -\Big(\alpha E_L+\sum_i\beta_iE_i\Big) -\gamma\big(a^2V_x\big)_x  = 0,
\end{multline}
where \[ \alpha=\frac{1}{r_mc_m}, \quad \beta_i=\frac{g_i}{c_m} \quad \text{and} \quad \gamma=\frac{1}{2ar_Lc_m}. \] \\
Relationships between specific properties $r_L$, $c_m$ and a compartment's actual physical properties are
\begin{equation}
\label{txtSpecificAxialResDefs}
\frac{1}{R_aC_m}=\frac{a^2/a}{2r_Lc_mh^2}~~\text{and}~~\frac{1}{R'_aC_m}=\frac{a'^2/a}{2r_Lc_mh^2},
\end{equation}
\\
where $a$ is the radius of a given compartment with subscript $j$, $a'$ the radius of neighboring compartment, subscripted $j+1$, and $h$ 
is the compartmentally defined spatial resolution. Relationships between $a$, $a'$ and the diffusion term's spatial derivatives in equation  (\ref{txtCableEquationOperator}) needed below include
\begin{multline}
\label{txtIdentitiesForTE}
\big(a^2V_x\big)_x\equiv\overline{a}^2V_{xx}+\big(\overline{a}^2\big)_xV_x, \\ \\[-1.5ex]
\quad \text{where} \quad\overline{a}^2\equiv\frac{a'^2+a^2}{2}, \quad \text{and}
\quad\big(\overline{a}^2\big)_x\equiv\displaystyle\lim_{h\rightarrow 0}\frac{a'^2-a^2}{h}. %\hfill
\end{multline}
\\
Next, we spatially discretize equation (\ref{txtCableEquationOperator}) and rewrite its diffusion term as a function of the compartment's 
actual physical properties giving us the ODE
\begin{multline}
\label{cableODEtxt}
\big(V_t\big)_j + \Big(\alpha+\sum_i\beta_i\Big)V_j -\Big(\alpha E_L+\sum_i\beta_iE_i\Big) ~ ... \hfill \\ \\[-2.0ex]
-\Bigg(\frac{V_{j-1}-V_j}{R_aC_m}+\frac{V_{j+1}-V_j}{R'_aC_m}\Bigg)=0,% ~ \cite[p.332,B.19]{Ster}
\end{multline}
\begin{multline}
\label{simpleODEtxt}
\text{which we may abreviate as } \quad \big(V_t\big)_j=A-BV_j, \hfill
\end{multline}
\begin{multline}
\label{simpleODEtxtDefs}
\text{where } A=\alpha E_L +\sum_i \beta_i E_i +\bigg(\frac{V_{j-1}}{R_aC_m} +\frac{V_{j+1}}{R'_aC_m}\bigg), \\
\text{ and } \quad B=\alpha +\sum_i \beta_i +\bigg(\frac{1}{R_aC_m} +\frac{1}{R'_aC_m}\bigg).
\end{multline}
For the context of a generic HH multi-compartment model, we define the format of an RK method as 
\begin{multline}
\label{rk21TruncErr0text}
V_j^{n+1}=V_j^n+kF(t_n,V_j^n,k;f) ~~\text{where}~~F(t_n,V_j^n,k;f)=\sum_{i=1}^s b_iK_i, \\ \\[-2.0ex]
K_i=f\big(t_n+c_ik,V_j^n+k\sum_{j=1}^sa_{ij}K_j\big), ~~ f(t_n,V_j^n) = \big(V_t\big)_j^n = A-BV_j^n, 
\end{multline}
$k$ is the integration step size and terms $A,B$ are defined in equation (\ref{simpleODEtxtDefs}). 
\begin{table}[h!]\small
\caption{RK21 Tableau}
\ \\
\centering
\begin{tabular}{l l}
\toprule \\[-2.0ex]
\begin{tabular}{l | l l}
$\displaystyle c_1$ & $a_{11}$  & \\
\hline \\[-2.0ex]
 & $b_1$ & $b_2$
\end{tabular}
&
$~ \Rightarrow ~$
\begin{tabular}{l | l l}
$1$ & $1$ & \\
\hline \\[-2.0ex]
 & $\frac{1}{2}$ &  $\frac{1}{2}$
\end{tabular}
\\ \\[-2.0ex]
\bottomrule
\end{tabular}
\label{tabInText:RK21}
\end{table}
%\\ %\\[-2.0ex]
RK21 rendered a quasi-FD method then becomes
\begin{multline}
\label{rk21TruncErr1text}
V_j^{n+1}=V_j^n+\frac{k}{2}\Big( K_1+K_2 \Big), ~~ \text{where} \hfill \\ \\[-2.0ex]
K_1=f\big(t_n,V_j^n\big)= A-BV_j^n, ~~ \text{and} \hfill \\ \\[-2.0ex]
K_2= f\big(t_n,V_j^n+kK_1\big)  = A-B\Big(V_j^n+k\big(A-BV_j^n\big)\Big). \hfill
\end{multline}
\\[-1.0ex]
After combining terms we may express the RK21 quasi-FD operator as
\begin{multline}
\label{rk21TruncErr2text}
P_{h,k}(V) = V_j^{n+1}-V_j^n-k\Big(A-BV_j^n \Big) + \frac{k^2}{2}\Big(B\big(A-BV_j^n\big)\Big)=0. \hfill 
\end{multline}
\\
Taylor expanding the operator, simplifying, dividing by $k$ and substituting the first occurrence of $A$ and $B$ with definitions in 
equation (\ref{simpleODEtxtDefs}) leads to
\begin{multline}
\label{rk21TruncErr3text}
V_t+\frac{k}{2}V_{tt}+\frac{k^2}{3!}V_{ttt}+O(k^3) %~ ... \\ \\[-2.0ex]
-\Bigg(\bigg[\alpha E_L+\sum_i\beta_iE_i +\bigg(\frac{V^n_{j-1}}{R_aC_m} +\frac{V^n_{j+1}}{R'_aC_m}\bigg)\bigg]\Bigg. ~...\\
\Bigg. - \bigg[\alpha +\sum_i \beta_i +\bigg(\frac{1}{R_aC_m} +\frac{1}{R'_aC_m}\bigg) \bigg] V_j^n \Bigg) +\frac{Bk}{2}\Big( A-BV_j^n \Big)=0. 
\end{multline}
\\
Taylor expanding neighboring compartment voltage terms, simplifying and applying definitions in equation (\ref{txtIdentitiesForTE}), brings us to
\begin{multline}
\label{rk21TruncErr4text}
V_t+\frac{k}{2}V_{tt}+\frac{k^2}{6}V_{ttt}
-\Bigg(\Big(\alpha E_L+\sum_i\beta_iE_i\Big)-\Big(\alpha+\sum_i\beta_i\Big)V^n_j\Bigg. ~... \hfill  \\
+\frac{a^2/a}{2r_Lc_mh^2}\bigg(-hV_x +\frac{h^2}{2}V_{xx}-\frac{h^3}{6}V_{xxx}+\frac{h^4}{24}V_{xxxx} \bigg)~ ... \hfill \\
\Bigg. +\frac{a'^2/a}{2r_Lc_mh^2}\bigg(hV_x+\frac{h^2}{2}V_{xx}+\frac{h^3}{3!}V_{xxx}+\frac{h^4}{4!}V_{xxxx}\bigg) \Bigg)%~...\\
+\frac{Bk}{2}\Big( A-BV_j^n \Big)=O(k^3,h^3).   
\end{multline}
Rearranging, again applying definitions in equation (\ref{txtIdentitiesForTE}) and eliminating third-order terms gives us
\begin{multline}
\label{rk21TruncErr5text}
V_t+\frac{V_{tt}}{2}k+\frac{V_{ttt}}{6}k^2 -\Bigg(\Big(\alpha E_L+\sum_i\beta_iE_i\Big)- \Big(\alpha+\sum_i\beta_i\Big)V^n_j  \Bigg.~... \hfill \\
\Bigg. +\gamma\bigg[\Big(\overline{a}^2V_{xx}+\big(\overline{a}^2\big)_xV_x\Big)+\Big(\overline{a}^2V_{xxxx} 
+2\big(\overline{a}^2\big)_xV_{xxx}\Big)\frac{h^2}{12} \bigg] \Bigg) %~ ... \\
+\frac{Bk}{2}\Big( A-BV_j^n \Big) =0.
\end{multline}
After replacing the remaining occurrence of $A-BV^n_j$ with what the first occurrence in equation (\ref{rk21TruncErr2text}) eventually 
became in equation \ref{rk21TruncErr5text} we have
\begin{multline}
\label{rk21TruncErr6text}
V_t+\frac{V_{tt}}{2}k+\frac{V_{ttt}}{6}k^2 
-\Bigg(\Big(\alpha E_L+\sum_i\beta_iE_i\Big) - \Big(\alpha+\sum_i\beta_i\Big)V^n_j \Bigg.~... \hfill \\
\hfill \Bigg. +\gamma\bigg[\big(\overline{a}^2V_x\big)_x+\Big(\overline{a}^2V_{xxxx} +2\big(\overline{a}^2\big)_xV_{xxx}\Big)\frac{h^2}{12} \bigg] \Bigg) ~ ... \\
+\frac{Bk}{2}\Bigg(\Big(\alpha E_L+\sum_i\beta_iE_i\Big) - \Big(\alpha+\sum_i\beta_i\Big)V^n_j \Bigg.~... \hfill \\
\Bigg. +\gamma\bigg[\big(\overline{a}^2V_x\big)_x+\Big(\overline{a}^2V_{xxxx} +2\big(\overline{a}^2\big)_xV_{xxx}\Big)\frac{h^2}{12} \bigg] \Bigg) =0.
\end{multline}
\begin{multline}
\label{rk21TruncErr7text}
\Rightarrow V_t+\Big(\alpha+\sum_i\beta_i\Big)V^n_j - \Big(\alpha E_L+\sum_i\beta_iE_i\Big)-\gamma\big(\overline{a}^2V_x\big)_x~... \hfill \\
+\frac{V_{tt}}{2}k+\frac{V_{ttt}}{6}k^2-\gamma\Big(\overline{a}^2V_{xxxx} +2\big(\overline{a}^2\big)_xV_{xxx}\Big)\frac{h^2}{12} ~ ... \\
+\frac{Bk}{2}\bigg(V_t+\gamma\Big(\overline{a}^2V_{xxxx} +2\big(\overline{a}^2\big)_xV_{xxx}\Big)\frac{h^2}{12}\bigg) =0,
\end{multline}
\\ \\[-2.0ex]
where the term $V_t$ in the last line of equation (\ref{rk21TruncErr7text}) came by using equation (\ref{txtCableEquationOperator}) to eliminate 
terms in the last two lines of equation (\ref{rk21TruncErr6text}). Then after discarding the third-order $kh^2$ term we are left with
\begin{multline}
\label{rk21TruncErr8text}
V_t+\Big(\alpha+\sum_i\beta_i\Big)V^n_j - \Big(\alpha E_L+\sum_i\beta_iE_i\Big)-\gamma\big(\overline{a}^2V_x\big)_x~... \hfill \\
+\Big(BV_t+V_{tt}\Big)\frac{k}{2}+\frac{V_{ttt}}{6}k^2-\gamma\Big(\overline{a}^2V_{xxxx} +2\big(\overline{a}^2\big)_xV_{xxx}\Big)\frac{h^2}{12}=0.
\end{multline}
\\
Terms $A$ and $B$ in equation (\ref{simpleODEtxt}) are both functions of the compartment's membrane voltage, $V^n_j$, through ion channel
conductance term $\sum\beta_i$, and, of course, $V^n_j$ changes over time. Hence, $V_{tt}=A_t-(BV_t+B_tV)$. Clearly the first-order term is 
not eliminated and truncation error for RK21 is
\begin{multline}
\label{rk21TruncErr9text}
P_{k,h}(V)-P(V) = ~ ... \hfill \\
\Big(BV_t+V_{tt}\Big)\frac{k}{2}+\frac{V_{ttt}}{6}k^2-\gamma\Big(\overline{a}^2V_{xxxx} +2\big(\overline{a}^2\big)_xV_{xxx}\Big)\frac{h^2}{12},
\end{multline}
which is only first-order in time, second in space. 
\subsection{Stability} 
\label{ssec:SchemeStab}
Results of Von Neumann stability analyses \cite[ch.2]{Str} for this study's six membrane voltage integration methods on an unbounded 
domain are listed in Table \ref{tab:schemeStability}. The bounded stability analysis presented in \cite[ch.4]{CarnevaleHines}, confirmed the 
same results shown in Table \ref{tab:schemeStability} for FTCS, BTCS and HCN.
\begin{table}[h!]\small
\caption{Integration Method Stability}
\ \\
\centering
\setlength{\tabcolsep}{1.0mm}
\begin{tabular}{l | l}
\toprule
Method & Condition for Growth Factor Magnitude, $|g(\theta,k)|^2<1$ \\
\midrule \\[-3.0ex] 
\parbox{2.5cm}{Forward-Time Central-Space} & $\displaystyle k<\min_{\text{AP cycle}}\bigg(\frac{2}{K+2L}\bigg)=7\mu$s. 
\\[2.0ex] \hline \\[-2.0ex]
\parbox{2.5cm}{Exponential \\ Euler} & Unconditionally Stable. 
\\[2.0ex] \hline \\[-2.0ex]
\parbox{2.5cm}{Backward-Time Central-Space} & Unconditionally Stable. 
\\ [2.0ex] \hline \\[-2.0ex] 
\parbox{2.5cm}{Hines-Crank-Nicolson} & Unconditionally Stable, but $\displaystyle k >\min_{\text{AP cycle}}\bigg(\frac{2}{K+2L}\bigg)$  \\ \\[-2.0ex]
 & $\Rightarrow-1<g(k,\theta)<0\Rightarrow$ sol'n oscillates with decreasing amplitude. \\ [1.0ex] \hline \\[-2.0ex]
\parbox{2.5cm}{RK21} & $\displaystyle k<\min_{\text{AP cycle}}\bigg(\frac{2}{B}\bigg)=14\mu$s.$~~$ \\ [2.0ex] \hline \\[-2.0ex] 
\parbox{2.5cm}{RK41} & $\displaystyle k<\min_{\text{AP cycle}}\bigg(\frac{\sim 2.7853}{B}\bigg)=20\mu$s. \\ [2.0ex] \hline \\[-2.0ex]
Definitions & $\displaystyle K=\Big(\alpha+\sum_i\beta_i\Big),~L=\bigg(\frac{1}{R_aC_m}+\frac{1}{R'_aC_m}\bigg)\sin^2{\frac{\theta}{2}},
~B$: see Table \ref{tab:schemeTruncError} \\
\bottomrule
\end{tabular}
\label{tab:schemeStability}
\end{table}
\section{Analysis Methods}
\label{sec:AnalysisMethods}
While our model neuron is composed of $74$ compartments, we assess only somatic compartment data as in \cite{T2005}. 
\subsection{Waveform Definitions}
\label{sec:waveformDefs}
We categorize AP spiking waveforms according to two distinct features. {\it Spikes} refer to pulse peaks that rise monotonically from 
less than $-0.04$ volts to over $-0.01$ volts.  {\it ADP} refers to an after-depolarization-pulse, a modest local maximum rarely rising above 
$-0.055$ volts following the last {\it Spike}.  AP waveforms are classified, as in Figures \ref{fig:p23rsAPclassFE_2OT_RK21RK41} and 
\ref{fig:p23rsAPclassBE_HCN}, using $Spikes$ and $ADP$ numbers as follows.
\begin{equation}
\label{eqn:classificationFormat}
\text{(number of $Spikes$ - number of $ADP$)}
\end{equation}
\subsection{Accuracy}
\label{sec:AccMeasure}
The correct solution for the active Rallpacks $3$ model could not be found analytically. Instead, Rallpacks provided simulator code ``used 
to generate the results in NEURON \cite{CarnevaleHines} and GENESIS \cite{GEN}. These models were'' integrated with the HCN method 
using a time-step of $1\mu$s for which correspondence between the two was better than $1\%$ \cite{BBB_Rall}. To assess this correspondence 
``... spike peaks [were] aligned for the voltage calculations, and interval differences [were] separately added to the root mean squared total 
\cite{BBB_Rall}.'' \\ \\
Likewise, our model's exact waveform solution cannot be determined analytically. Furthermore, Rallpacks' choice of an AP cycle waveform 
integrated by HCN using a $1\mu$s time-step as the standard of accuracy works well only if different step size waveforms compared 
to it are also integrated by HCN. This because each integration method generates different AP spiking waveform phases as shown in Figure 
\ref{fig:SpikePhaseComp}, even for $1\mu$s step sizes. For this reason we have chosen the standard of accuracy for each method as the 
AP cycle waveform integrated by the same method when using $1\mu$s step sizes. \\ \\
Note that because AP cycle absolute minimums and periods do not stabilize until the twentieth cycle, as shown in Figures 
\ref{fig:minsVapidx}-\ref{fig:periodsVapidx}, model statistics in Figure \ref{fig:somaStatsWithStd} were only based on cycles following the nineteenth AP. 
Likewise, accuracy of simulation waveforms, shown in Figure \ref{fig:methodErrVsts}, was determined by comparing the synchronized twentieth 
AP cycle from each step size's simulation to the reference waveform, that being the twentieth AP cycle generated by the same integration 
method using a step size of $1\mu$s.
\subsection{Spectral Analysis}
\label{subsection:SpecAnalysis}
Power spectral density (PSD) graphs in Figures \ref{fig:HCNandBEPSD}-\ref{fig:HCNandExpoEPSD} are periodograms. Integration step sizes, 
$1$ to $100\mu$sec, correspond to, respectively, sampling frequencies of $1$ MHz and $10$ KHz. But the latter is higher than necessary since 
neuroscience power spectrum interests coincide with e.e.g. bands, all less than $64$ Hz.  Simulation voltages were therefore downsampled to 
achieve a sampling frequency of $250$ Hz allowing us to see the distribution of signal power as high as $125$ Hz. Power spectral densities were 
then computed with the Welch periodogram in MATLAB; sampling frequency was $250$ Hz, Hamming window and DFT widths were $250$ 
samples with $150$ sample overlaps. \\ \\
We define the spectral centroid of an FFT computed from compartmentally sampled membrane voltages, used to compute step size upper limits 
shown in Figures \ref{fig:fEstepSizeLimit}-\ref{fig:RKstepsizeLimits}, by writing the magnitude of a plane wave's spatial frequency spectrum as 
$|H(\omega)|$, a function of frequency in radians where 
\begin{equation}
\label{eqn:omegaSpan}
\omega \equiv \omega_i=\frac{i\pi}{16}\quad\text{for } i = 0, 1, 2, ... , 16.
\end{equation}
The spectral centroid is then
\begin{equation}
\label{eqn:spectralCentroid}
\omega_c=\frac{\sum_i \omega_i|H(\omega_i)|}{\sum_i |H(\omega_i)|}.
\end{equation}
\subsection{Oscillation Magnitude and Duration}
\label{sec:OscDet}
\begin{figure}[h!]
\centering
\includegraphics[height=2.0in,width=4.0in]{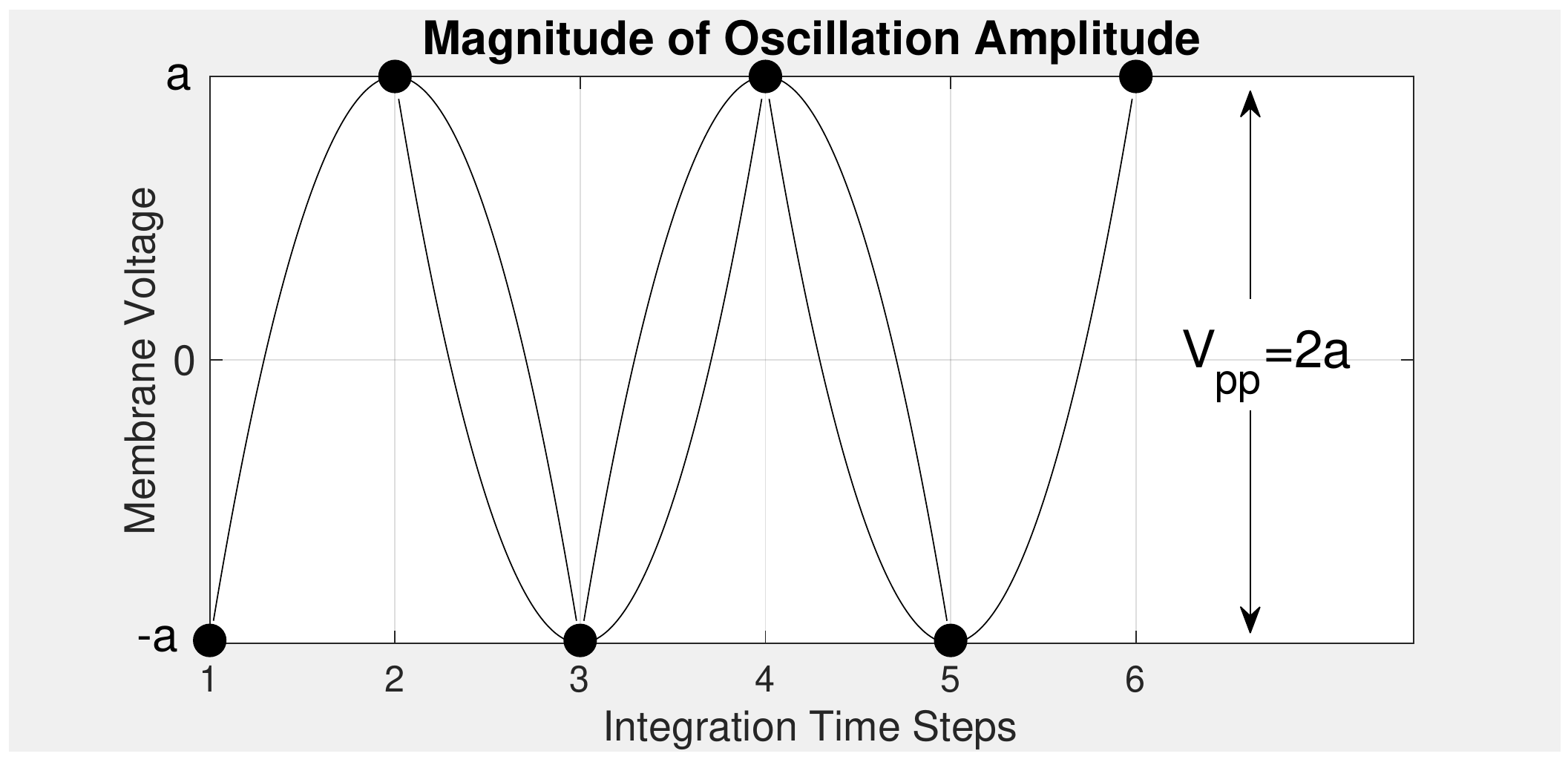}
\caption{Oscillation amplitude is one quarter the second undivided difference.}
\label{fig:oscMagMeasure}
\end{figure}

\begin{figure}[h!]
\centering
\includegraphics[height=4.5in,width=4.0in]{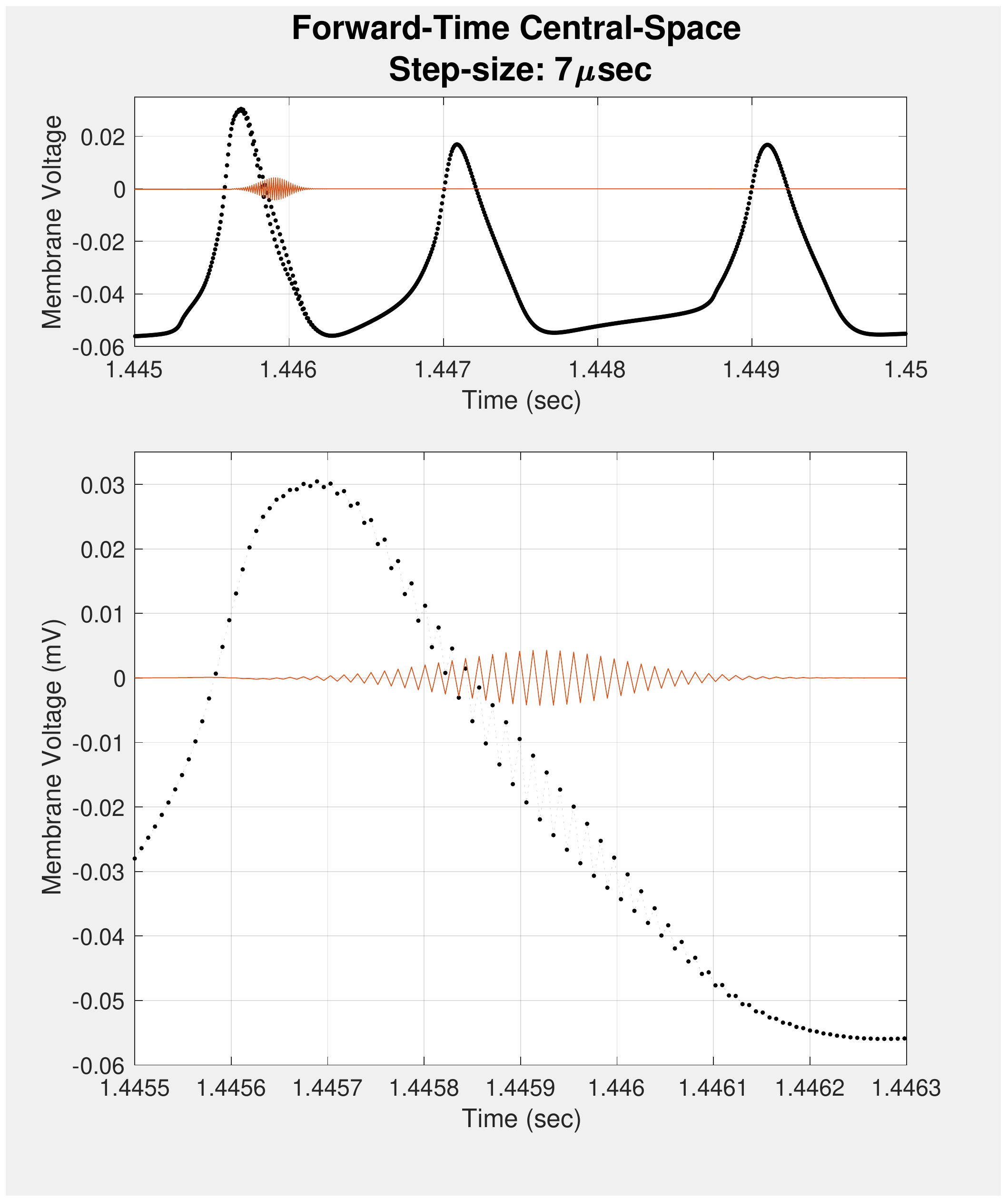}
\caption{The measure of oscillation amplitude is plotted in red. The lower pane is a closeup of the first AP spike in the upper frame. 
Note that, contrary to observations described in \cite{Borgers}, instability can present on the descending slopes of an AP spike.}
\label{fig:exOsc}
\end{figure}

Our oscillation search method looks for concavities with alternating polarities. Two consecutive nonzero concavities with alternate polarities 
are however common in AP cycle waveforms near membrane voltage spiking swings, especially as integration step size grows larger than 
$50\mu$s. So, we require greater than two consecutive nonzero concavities with alternate polarities for a membrane voltage waveform 
to be considered an oscillation. \\ \\
Next we consider how to best approximate concavity. Divided differences, the discrete approximation of derivatives, represent the estimate of 
an instantaneous rate of change. In contrast, undivided differences are an approximation of a differential, the estimate of an amount, not a rate.  
Since we are most interested in the magnitude an oscillating membrane voltage swings to reach the opposite polar extreme in one integration 
step, regardless of step size, we will estimate concavity with the second undivided difference. \\ \\
The shortest oscillation is a sequence of membrane voltages whose second undivided differences change sign three or more consecutive 
integration steps. Figure \ref{fig:oscMagMeasure} illustrates a sequence of four concavities with alternating polarity to demonstrate that amplitude 
magnitude is one half the peak-to-peak voltage, or identically one quarter the second undivided difference. 
Figure \ref{fig:exOsc} shows an example of membrane voltage oscillation induced by FTCS using a step size of $7\mu$s. The magnitude of 
second undivided differences, in red, is the basis for computing the RMS of oscillations shown in Figures 
\ref{fig:fE2OToscVsTimestep}-\ref{fig:hCNoscVsTimestep}.  
%\clearpage
\section{Results}
\label{sec:Results}
Simulation software was developed by the author in MATLAB. Table \ref{tab:SomaWaveSummary} lists graphs of results and their 
corresponding figure numbers. 
\begin{table}[h!]\footnotesize
\caption{P23RS Somatic Waveform Performance Figures}
\ \\
\centering
\renewcommand{\arraystretch}{0.5}% Tighter  
\begin{tabular}{l|c}
\toprule 
Figure Summary & Figure \\
\midrule
\parbox{10.5cm}{Stable step size intervals versus method} &  \ref{fig:successfulStepSizes} \\ \\
\parbox{10.5cm}{Ideal AP cycle, integrated by HCN, step size $1\mu$s} & \ref{fig:idealAPcyclesHCNsts1} \\ \\
\parbox{12.5cm}{Absolute minimums between AP cycles reach steady state by $20$ AP cycles} & \ref{fig:minsVapidx} \\ \\
\parbox{12.5cm}{AP periods reach steady state by $20$ AP cycles}  & \ref{fig:periodsVapidx} \\ \\
\parbox{10.5cm}{Integration method accuracy comparison} & \ref{fig:methodErrVsts} \\ \\
\parbox{10.5cm}{Integration method statistics with standard deviations} &  \ref{fig:somaStatsWithStd} \\ \\
\parbox{10.5cm}{Spiking classifications for FTCS, RK21 and RK41 methods} & \ref{fig:p23rsAPclassFE_2OT_RK21RK41} \\ \\
\parbox{10.5cm}{Spiking classifications for BTCS and HCN versus step size} & \ref{fig:p23rsAPclassBE_HCN} \\ \\ 
\parbox{10.5cm}{Spiking classifications for exponential Euler versus step size} & \ref{fig:p23rsAPclassExpoE} \\ \\ 
\parbox{10.5cm}{Examples of exponential Euler AP spiking distortions}  &  \ref{fig:corruptExpoEulerAPs} \\ \\
\parbox{10.5cm}{Integration method AP spike phases at step sizes $1,9,50,80$ and $99\mu$s} & \ref{fig:SpikePhaseComp} \\ \\
\parbox{12.5cm}{Comparison of power spectral density (PSD) for HCN and BTCS methods} & \ref{fig:HCNandBEPSD} \\ \\
\parbox{12.5cm}{Comparison of PSD for HCN and exponential Euler methods} & \ref{fig:HCNandExpoEPSD} \\ \\
\parbox{12.5cm}{Maximum step size for FTCS versus AP cycle time} &  \ref{fig:fEstepSizeLimit} \\ \\
\parbox{12.5cm}{RK step size limits with/without consideration of plane wave phase angle, $\theta$} &  \ref{fig:RKstepsizeLimits} \\ \\
\parbox{10.5cm}{RMS of FTCS, RK21 and RK41 oscillation magnitudes} & \ref{fig:fE2OToscVsTimestep} \\ \\
\parbox{10.5cm}{RMS of exponential Euler and BTCS oscillation magnitudes vs step size} & \ref{fig:eE_BTCSoscVsTimestep} \\ \\
\parbox{10.5cm}{RMS of HCN oscillation versus step size} & \ref{fig:hCNoscVsTimestep} \\ \\
\parbox{10.5cm}{HCN Growth factor span versus step size} & \ref{fig:hCNGrwthFctrVsTimestep} \\ \\
\parbox{12.5cm}{HCN oscillations during Spiking phase} & \ref{fig:Spiking} \\ \\
\parbox{10.5cm}{HCN oscillations during ADP phase} & \ref{fig:ADP} \\ \\
\parbox{10.5cm}{HCN oscillations during Max Polarization phase} & \ref{fig:maxPolarization} \\ \\
\parbox{10.5cm}{HCN oscillation magnitude increases with step size and decays before the next AP cycle begins} & \ref{fig:hCNOscDecay} \\ 
\bottomrule
\end{tabular}
\label{tab:SomaWaveSummary}
\end{table}
\subsection{Accuracy}
\label{sec:Accuracy}
\begin{figure}[h!]
\centering
\begin{tabular}{l}
\includegraphics[height=3.5in,width=3.5in]{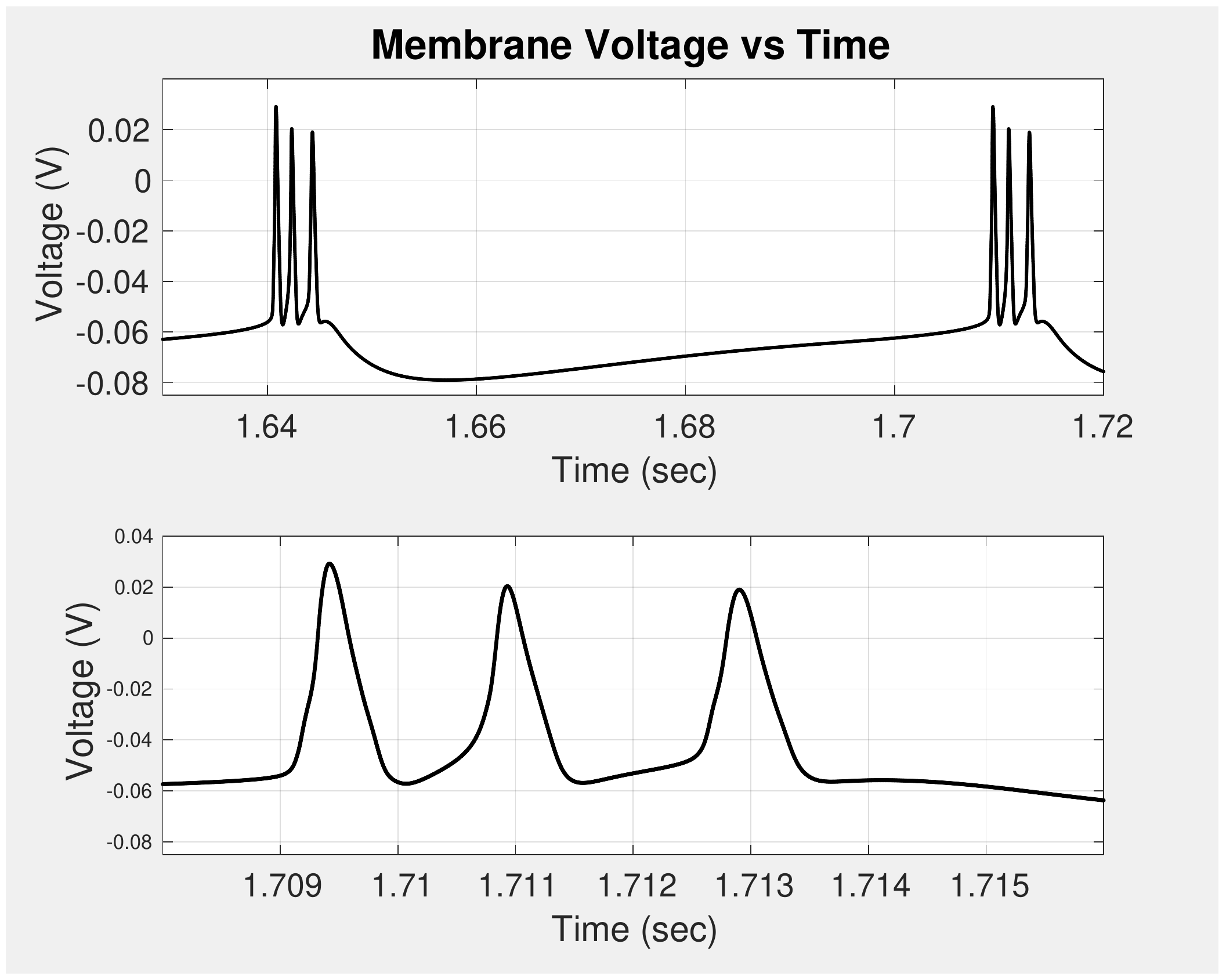}
\end{tabular}
\caption{The ideal waveform, a triplette, whose spike heights descend in order, and one modest after-depolarization-pulse (ADP).}
\label{fig:idealAPcyclesHCNsts1}
\end{figure}

\begin{figure}[h!]
\centering
\setlength{\tabcolsep}{0.7mm}
\begin{tabular}{l}
\includegraphics[height=6.0in,width=4.5in]{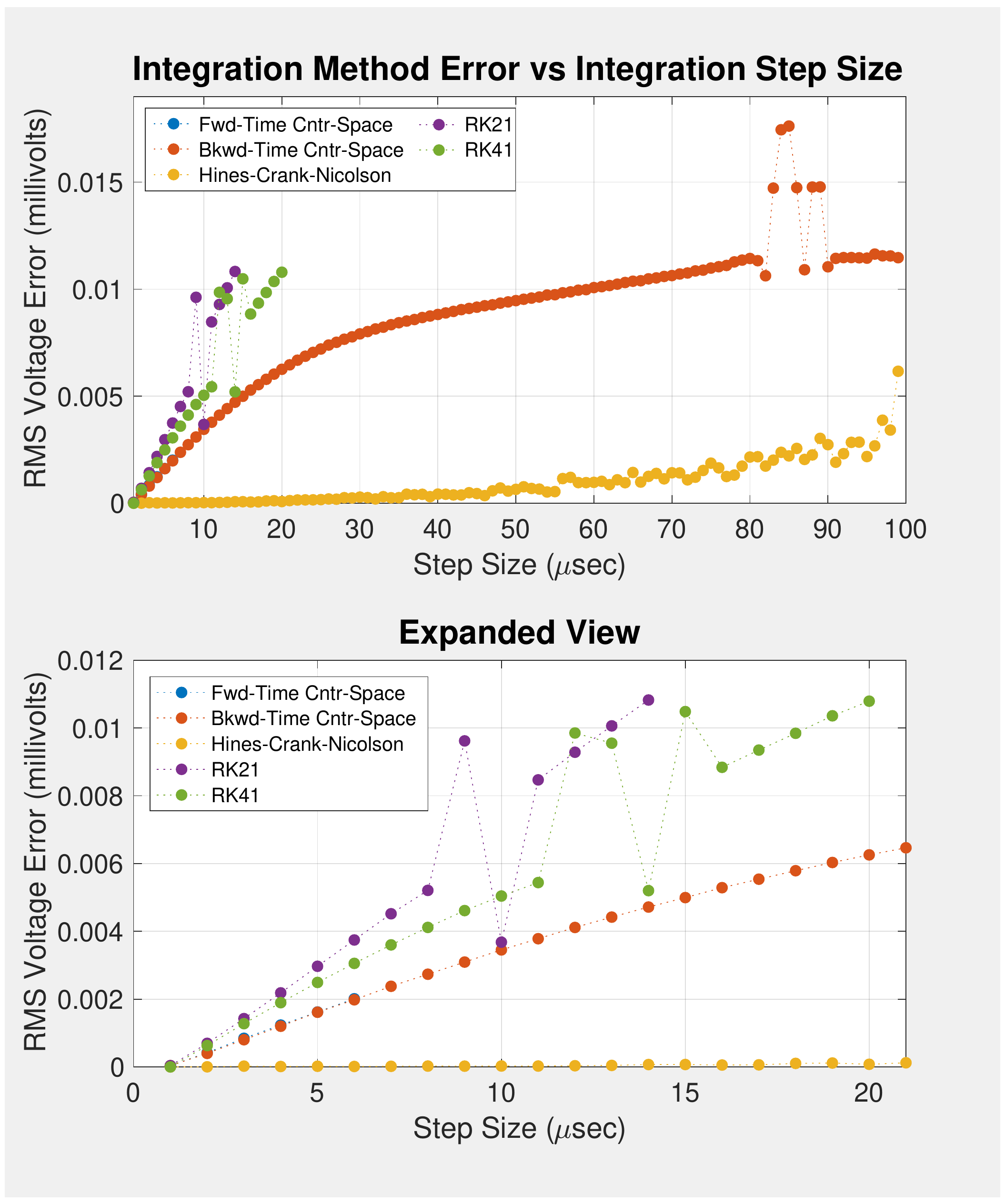} \\
\end{tabular}
\caption{Accuracy of membrane voltage integration found by comparing the twentieth AP cycle at each step size to the same AP cycle, 
integrated by the same method with a $1\mu$s step size. FTCS data presents typical first-order, straight line, growth of error, and 
is hidden behind BTCS data. BTCS data presents square-root-like, $\frac{1}{2}$-order, error growth. HCN presents typical second-order, 
concave-up, growth of error. As predicted, both RK methods present straight line, first-order, error growth. Of particular interest are shocks 
in BTCS, RK21 and RK41 data that indicate step sizes beyond which loss of accuracy for each first order method leads to spiking morphology 
changes. The expanded view on the bottom supports approximately $\frac{1}{2}$-Order growth of error for BTCS as well as RK41.}
\label{fig:methodErrVsts}
\end{figure}

\begin{figure}[h!]
\centering
\setlength{\tabcolsep}{0.7mm}
\begin{tabular}{l}
\includegraphics[height=2.25in,width=4.4in]{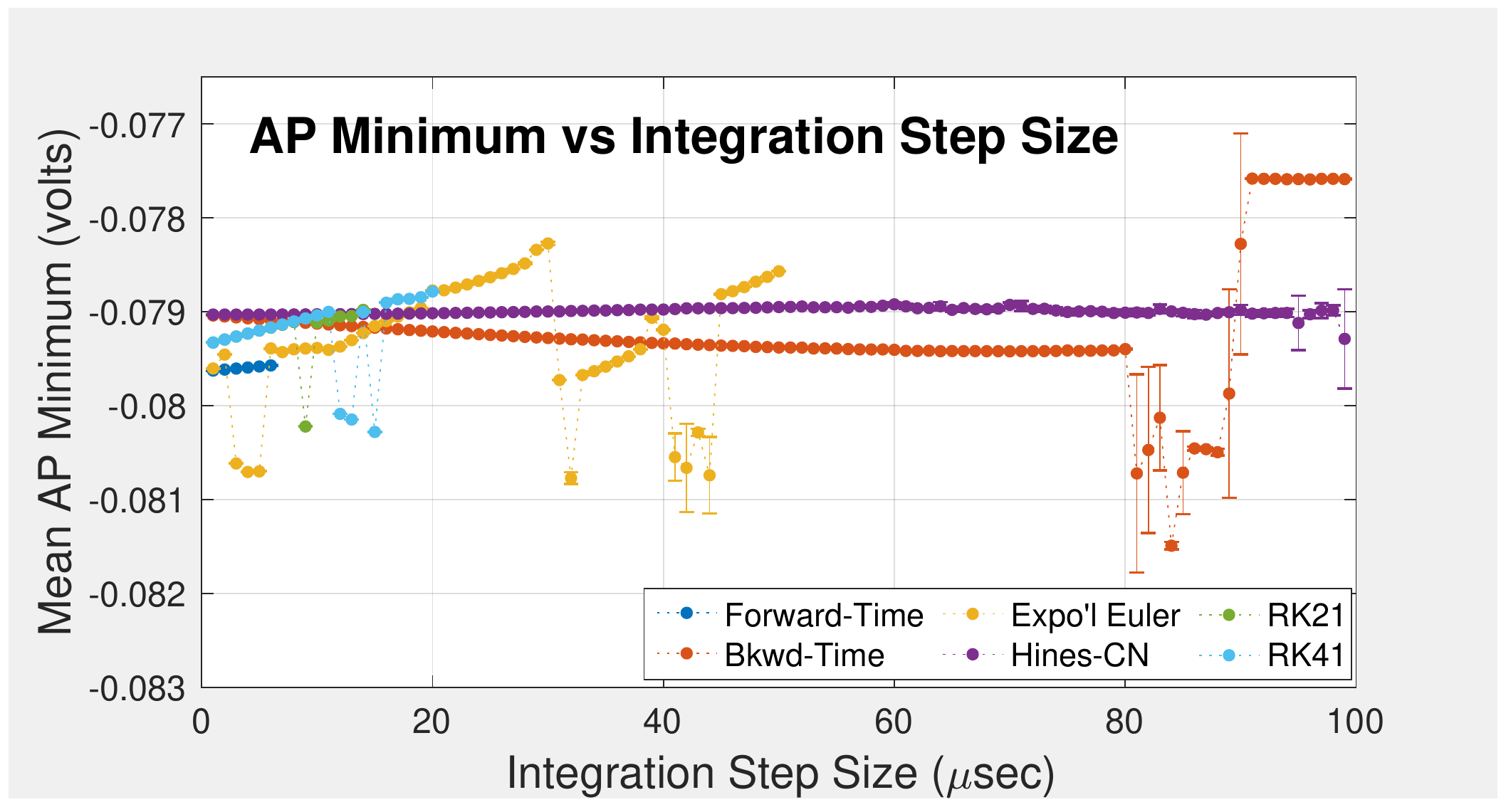} \\
\includegraphics[height=2.25in,width=4.4in]{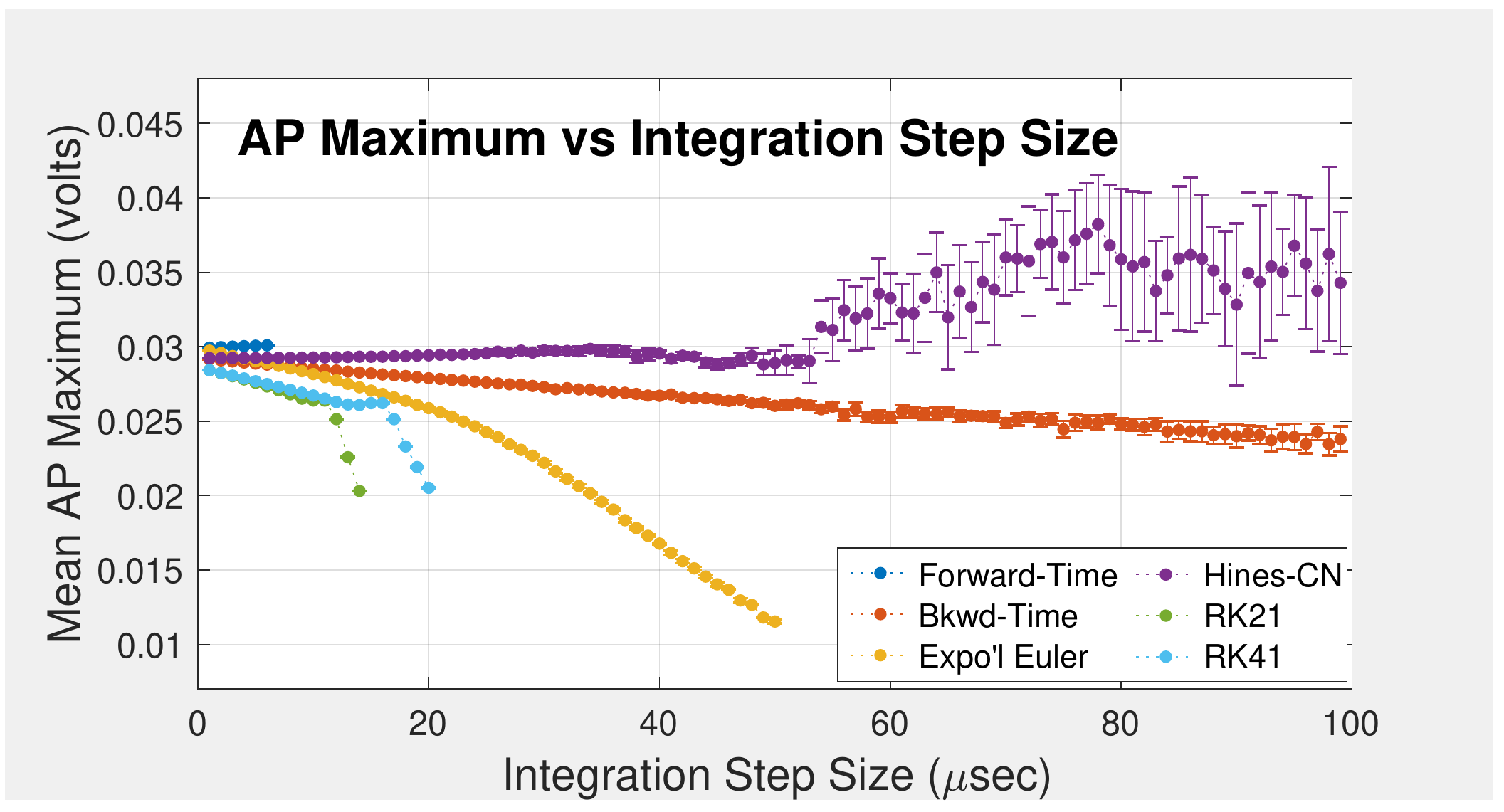} \\
\includegraphics[height=2.25in,width=4.4in]{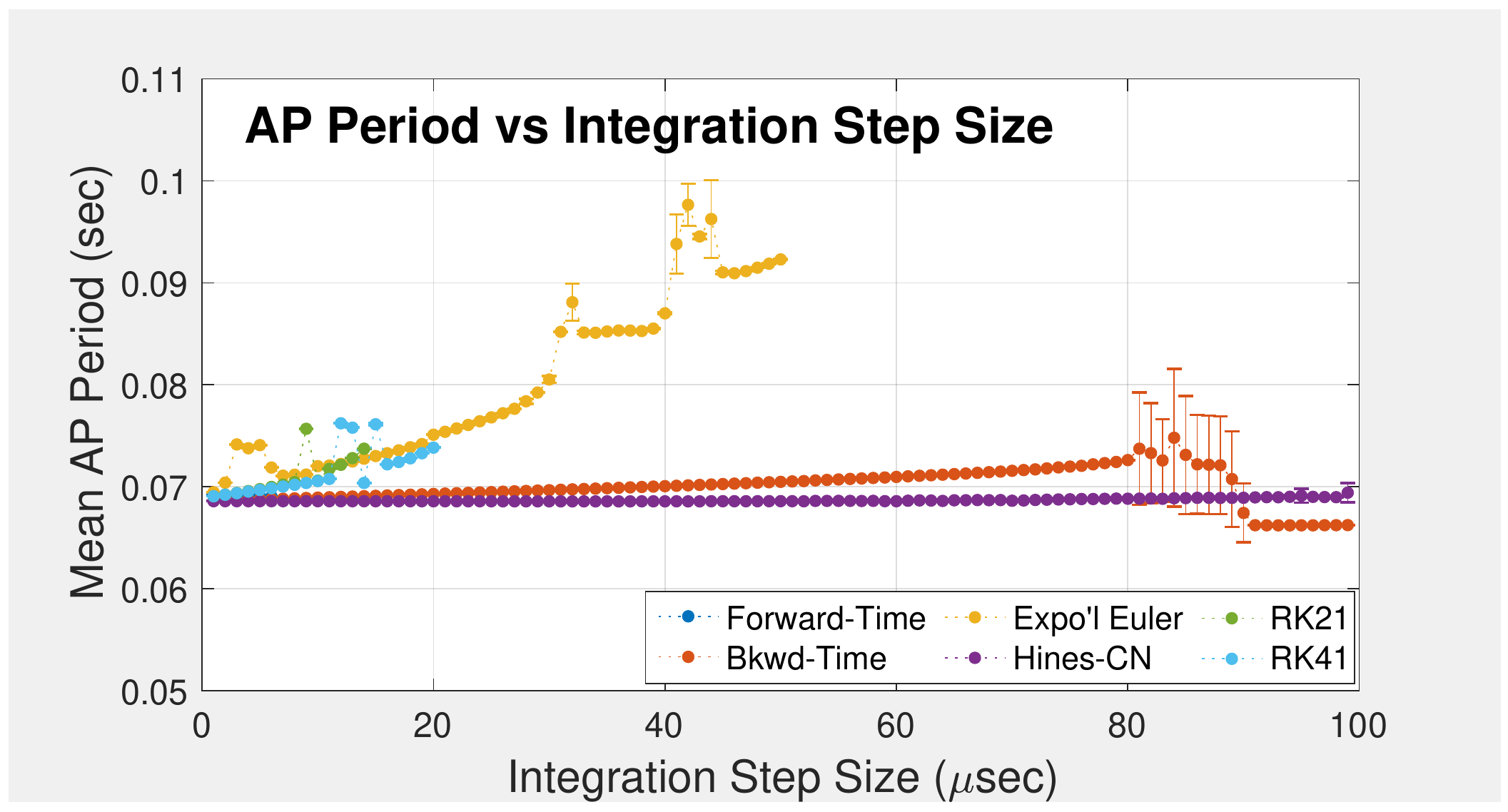}
\end{tabular}
\caption{Comparison of somatic compartment statistical means and standard deviations based on \mbox{AP cycles $20-$end}.}
\label{fig:somaStatsWithStd}
\end{figure}

\begin{figure}[h!]
\centering
\begin{tabular}{l}
\includegraphics[height=2.7in,width=3.5in]{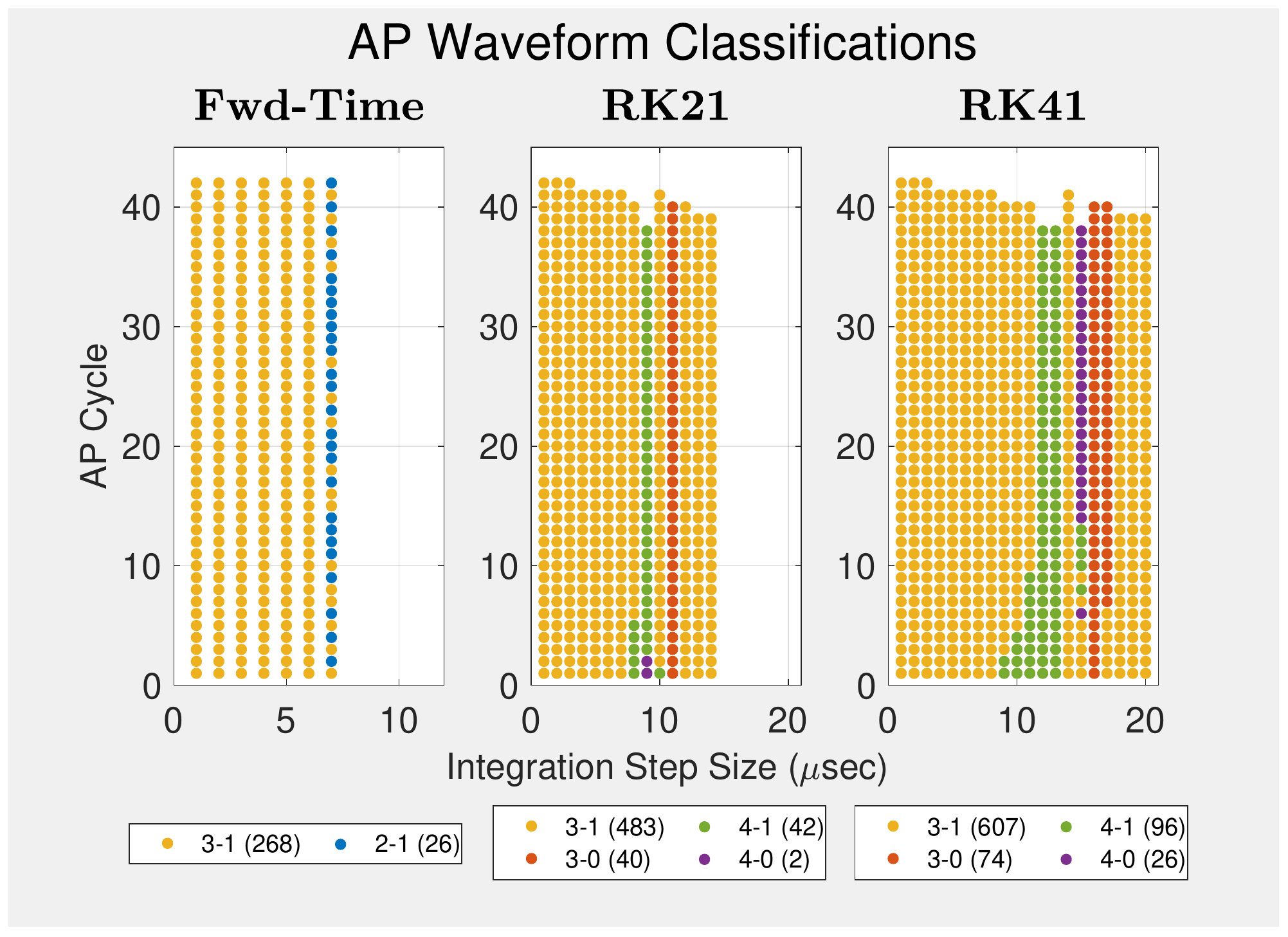} 
\end{tabular}
\caption{FTCS spiking morphology deteriorates just $1\mu$s prior to the onset of instability. RK21 and RK41, both rendered first-order 
methods, present deteriorating spiking morphologies for step sizes larger than $8\mu$s and $9-10\mu$s, respectively.}
\label{fig:p23rsAPclassFE_2OT_RK21RK41}
\end{figure}

\begin{figure}[h!]
\centering
\begin{tabular}{l}
\includegraphics[height=2.8in,width=2.4in]{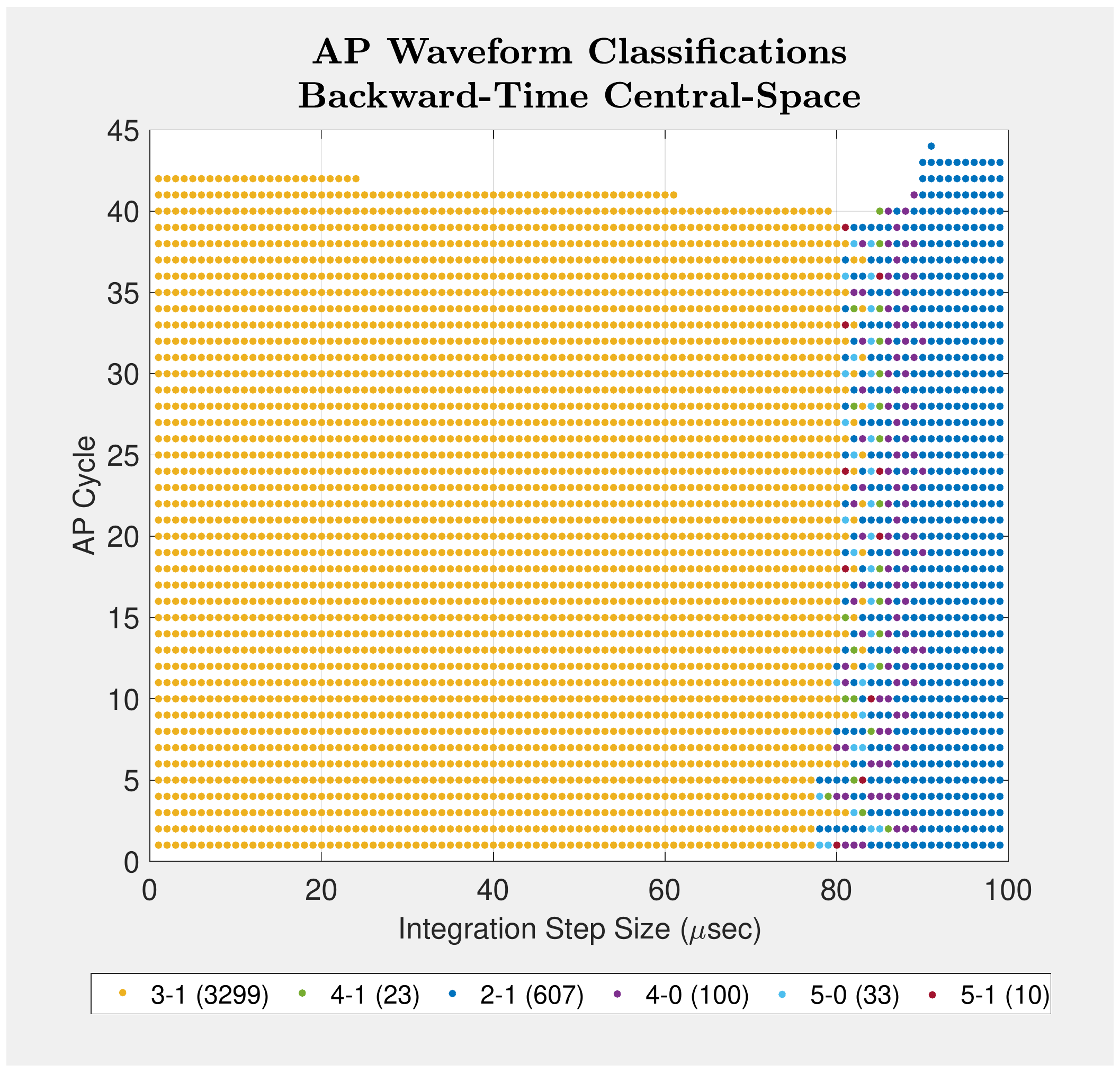}
\includegraphics[height=2.8in,width=2.4in]{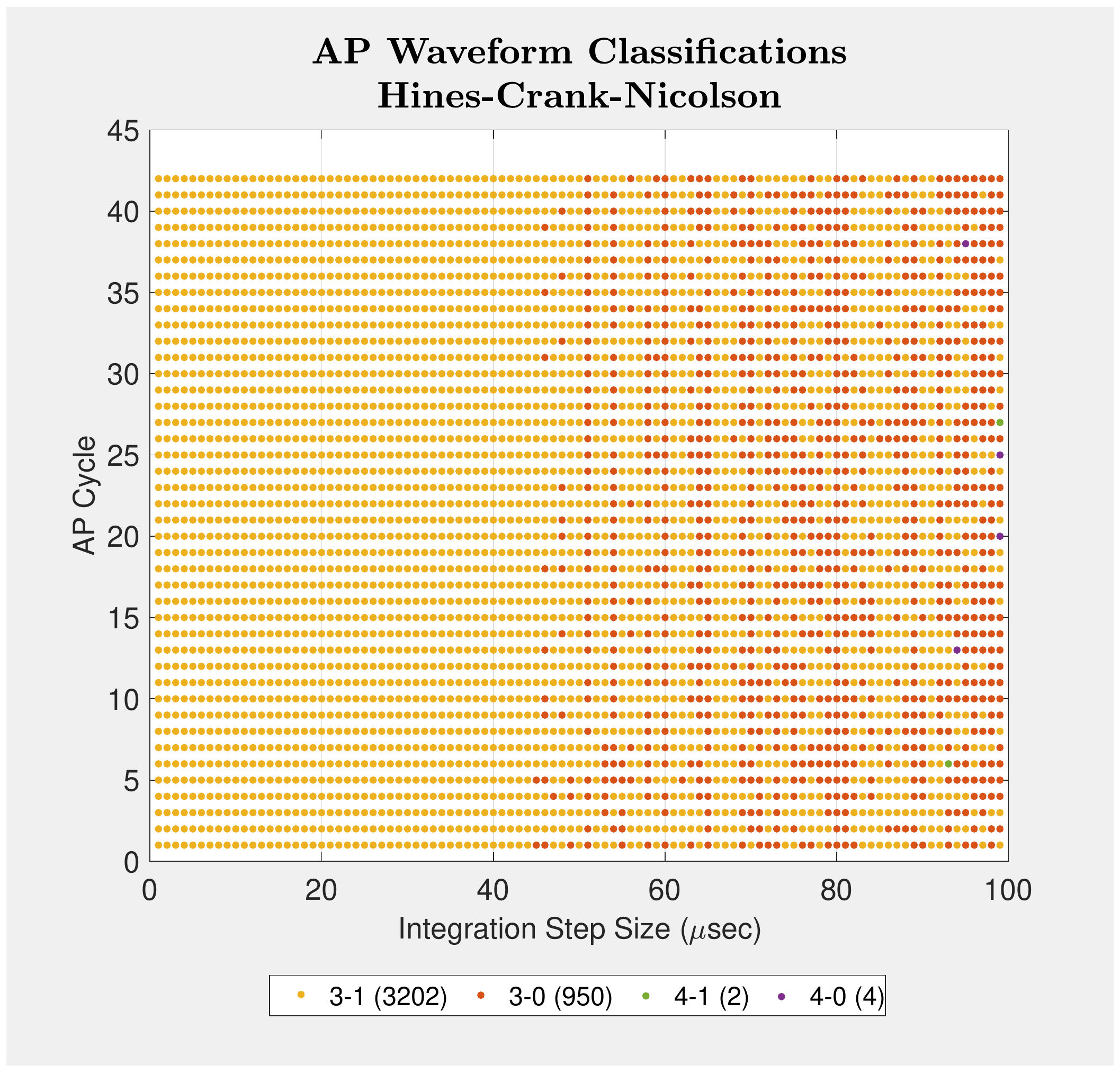}
\end{tabular}
\caption{BTCS presents deteriorating spiking morphologies for step sizes $>78\mu$s. HCN's class $3$-$0$ waveforms, i.e 
apparent loss of ADPs, are an inaccurate artifact. They are in fact class $3$-$1$ whose oscillations prevented software from detecting 
the ADP local maximum. This was intentionally left to visually corroborate oscillations in \mbox{Figure \ref{fig:hCNoscVsTimestep}.}}
\label{fig:p23rsAPclassBE_HCN}
\end{figure}

\begin{figure}[h!]
\centering
\begin{tabular}{l}
\includegraphics[height=3.5in,width=3.8in]{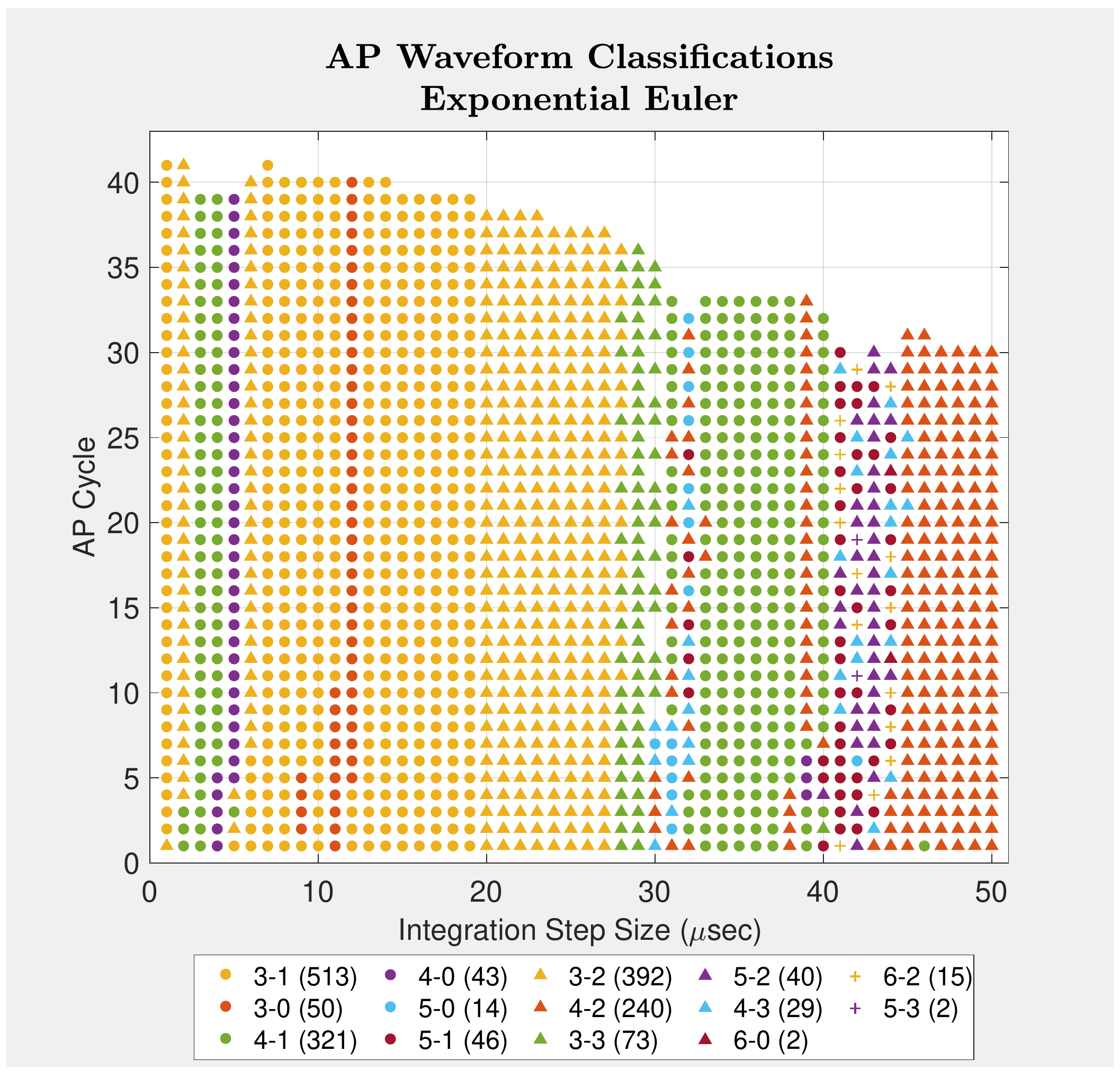}
\end{tabular}
\caption{AP spiking classification map for exponential Euler integration method .}
\label{fig:p23rsAPclassExpoE}
\end{figure}

\begin{figure}[h!]
\centering
\begin{tabular}{l}
\includegraphics[height=3.0in,width=4.5in]{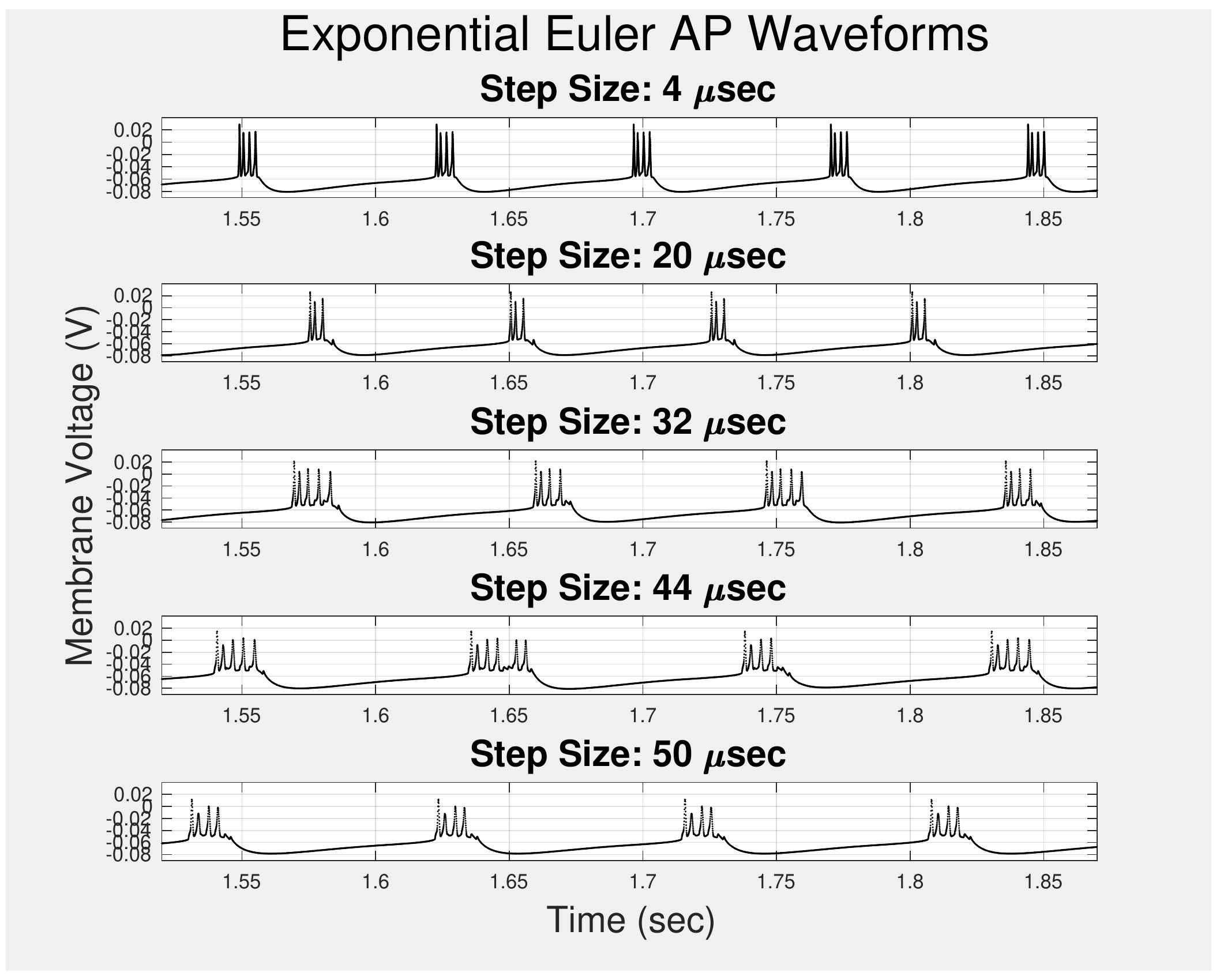}
\end{tabular}
\caption{Corrupted AP cycle waveforms integrated by the exponential Euler method.}
\label{fig:corruptExpoEulerAPs}
\end{figure}

\begin{figure}[h!]
\centering
\setlength{\tabcolsep}{0.7mm}
\begin{tabular}{l}
\includegraphics[height=7.0in,width=5.0in]{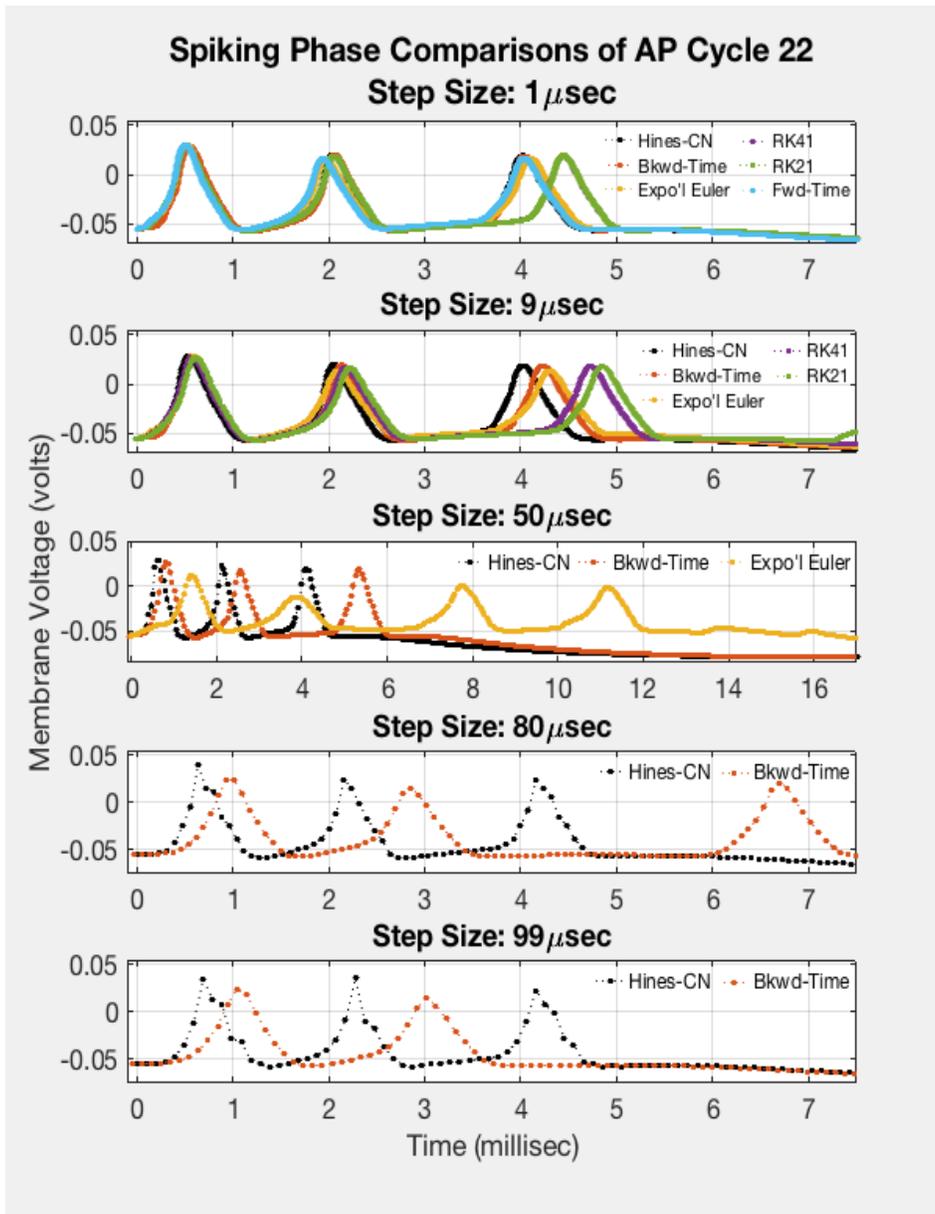}
\end{tabular}
\caption{AP cycle 22 spike phase comparisons at step sizes of $1,9,50,80$ and $99\mu$s.}
\label{fig:SpikePhaseComp}
\end{figure}

\begin{figure}[h!]
\centering
\setlength{\tabcolsep}{0.7mm}
\begin{tabular}{l}
\includegraphics[height=5.9in,width=4.4in]{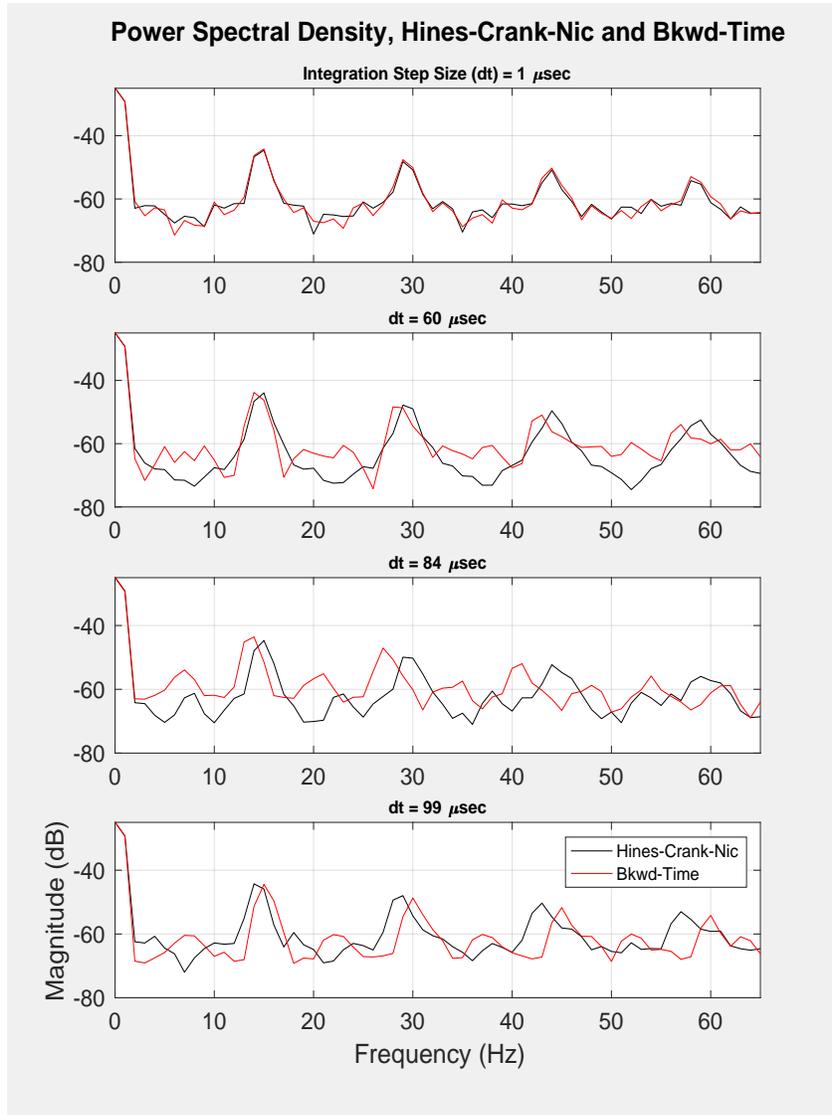}
\end{tabular}
\caption{Spectral density comparison showing the power distribution for HCN remaining largely fixed regardless of integration step size. 
BTCS is just as reliable up until $80\mu$sec step sizes where a likely bifurcation incursion changes spiking multiplicities for all AP cycles as 
shown in Figure \ref{fig:p23rsAPclassBE_HCN}. For both methods, the lowest frequency spectral peaks, those near $15$ Hz correspond 
to each numerical method's respective fundamental AP frequency. All other peaks are mere harmonics. }
\label{fig:HCNandBEPSD}
\end{figure}

\begin{figure}[h!]
\centering
\setlength{\tabcolsep}{0.7mm}
\begin{tabular}{l}
\includegraphics[height=5.9in,width=4.4in]{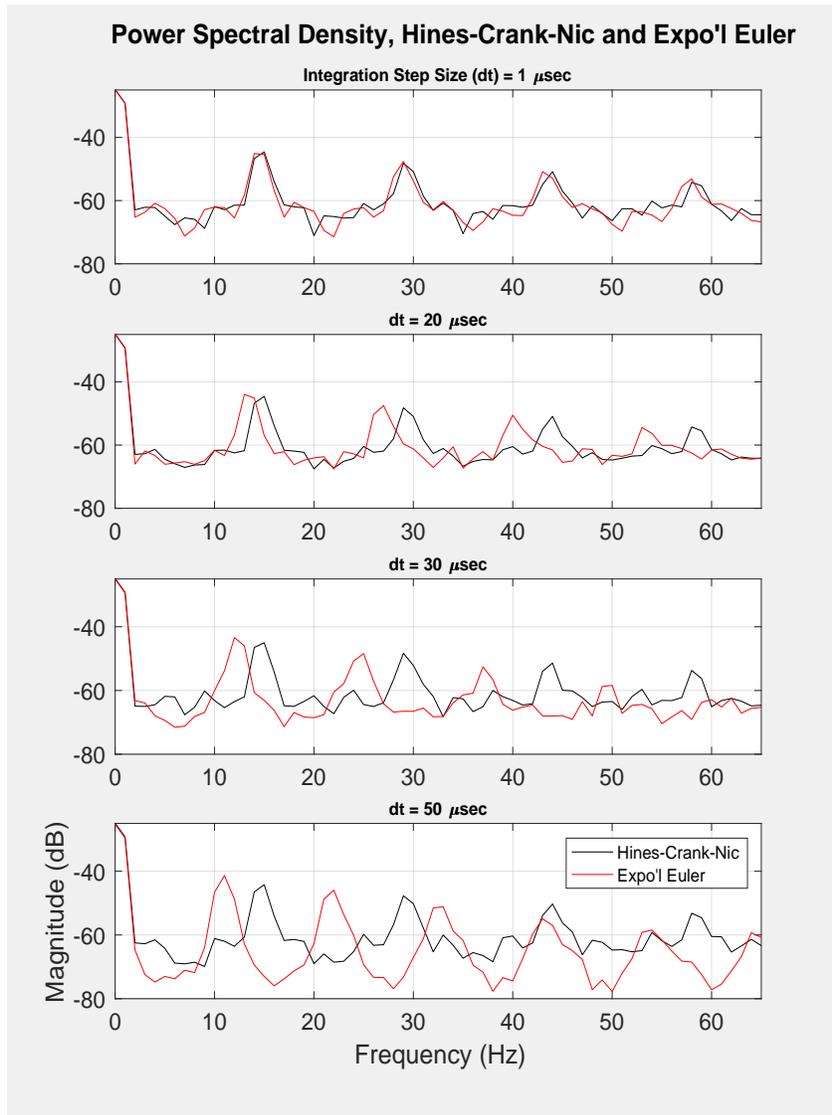}
\end{tabular}
\caption{Spectral density comparison showing the frequency of exponential Euler's power peaks decrease proportionally as integration 
step size increases.  This frequency shift is also corroborated by the decreasing AP cycle frequency, or increasing AP cycle period shown 
in Figures \ref{fig:somaStatsWithStd} and \ref{fig:p23rsAPclassExpoE}.}
\label{fig:HCNandExpoEPSD}
\end{figure}

Figure \ref{fig:idealAPcyclesHCNsts1} shows the ideal AP cycle and a closer view of its three spikes as generated by the HCN method with 
a step size of $1\mu$s. AP cycle waveform minimum's and periods for all methods required almost twenty AP cycles to reach macro view 
equilibriums, regardless of step size, exemplified by three methods shown in Figures \ref{fig:minsVapidx}-\ref{fig:periodsVapidx}. Membrane ion 
gate phase portraits included in the supplemental resources also support this observation.  AP waveform statistics shown in \mbox{Figure 
\ref{fig:somaStatsWithStd}} were therefore based on APs following the nineteenth cycle. \\ \\
Accuracy comparisons, as described in Section \ref{sec:AccMeasure}, appear in Figure \ref{fig:methodErrVsts}. For every integration method, 
spiking variations increased as step size increased.  For example, the center pane in Figure \ref{fig:somaStatsWithStd} shows HCN's maximum 
spike height stable until step sizes reached $53\mu$s, after which both mean and standard variation increased. In contrast, maximum spike 
heights for every other integration method decreased steadily as step size grew larger than $1\mu$s. Increasing step sizes led to a varying 
number of spikes in each AP cycle as shown in Figures \ref{fig:p23rsAPclassFE_2OT_RK21RK41}--\ref{fig:p23rsAPclassExpoE}. Prominent among 
them, waveforms generated by the exponential Euler method, as shown in Figure \ref{fig:corruptExpoEulerAPs}, were the most variable and 
inconsistent. Artifacts of accuracy loss also include changing AP spike phase presented by each method, except HCN, as shown in Figure 
\ref{fig:SpikePhaseComp}. PSDs of the exponential Euler and BTCS methods were each compared to the PSD of HCN in Figures 
\ref{fig:HCNandBEPSD}-\ref{fig:HCNandExpoEPSD}.
\subsection{Stability}
\label{sec:Stability}
\begin{figure}[h!]
\centering
\begin{tabular}{l}
\includegraphics[height=1.75in,width=4.0in]{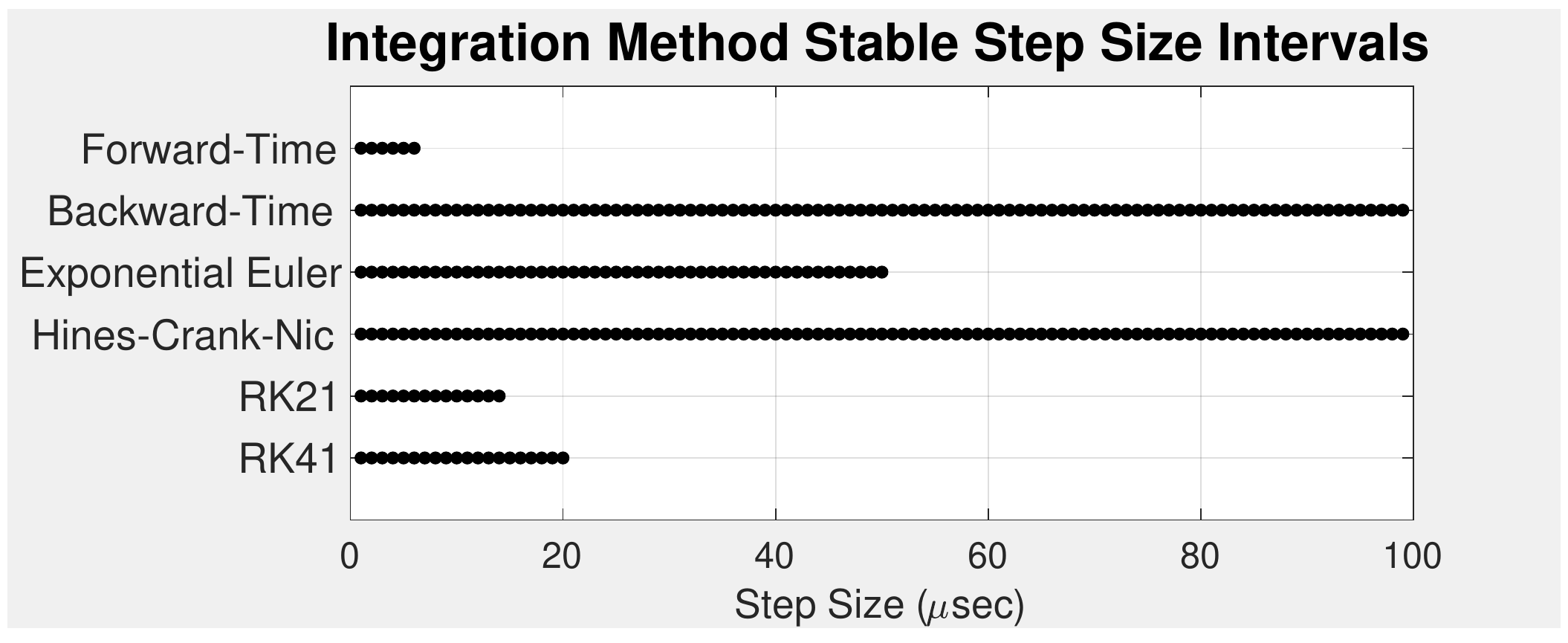} %\\
\end{tabular}
\caption{Stable step size intervals integrated by each numerical method.}
\label{fig:successfulStepSizes}
\end{figure}

\begin{figure}[h!]
\centering
\setlength{\tabcolsep}{0.7mm}
\begin{tabular}{l}
\includegraphics[height=4.5in,width=3.9in]{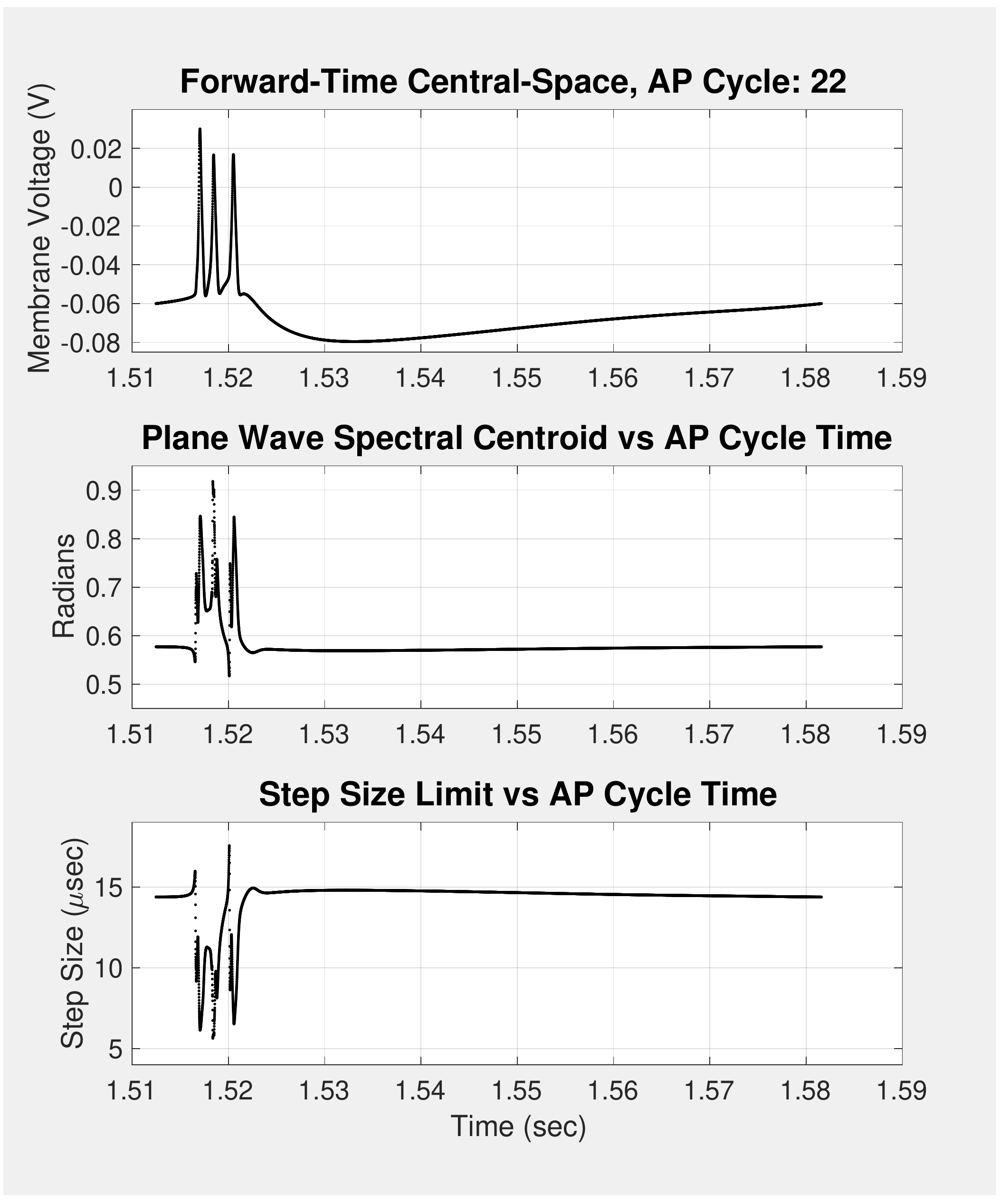} 
\end{tabular}
\caption{Graphs above illustrate correspondence between membrane voltage, phase angle of the model's plane wave from tip of 
distal dendrites to axon tip, and the integration step size limit for FTCS. The step size limit shown here for FTCS of $6\mu$sec matches 
simulation data shown in Figure \ref{fig:successfulStepSizes}.}
\label{fig:fEstepSizeLimit}
\end{figure}

\begin{figure}[h!]
\centering
\setlength{\tabcolsep}{0.7mm}
\begin{tabular}{l}
\includegraphics[height=1.6in,width=3.7in]{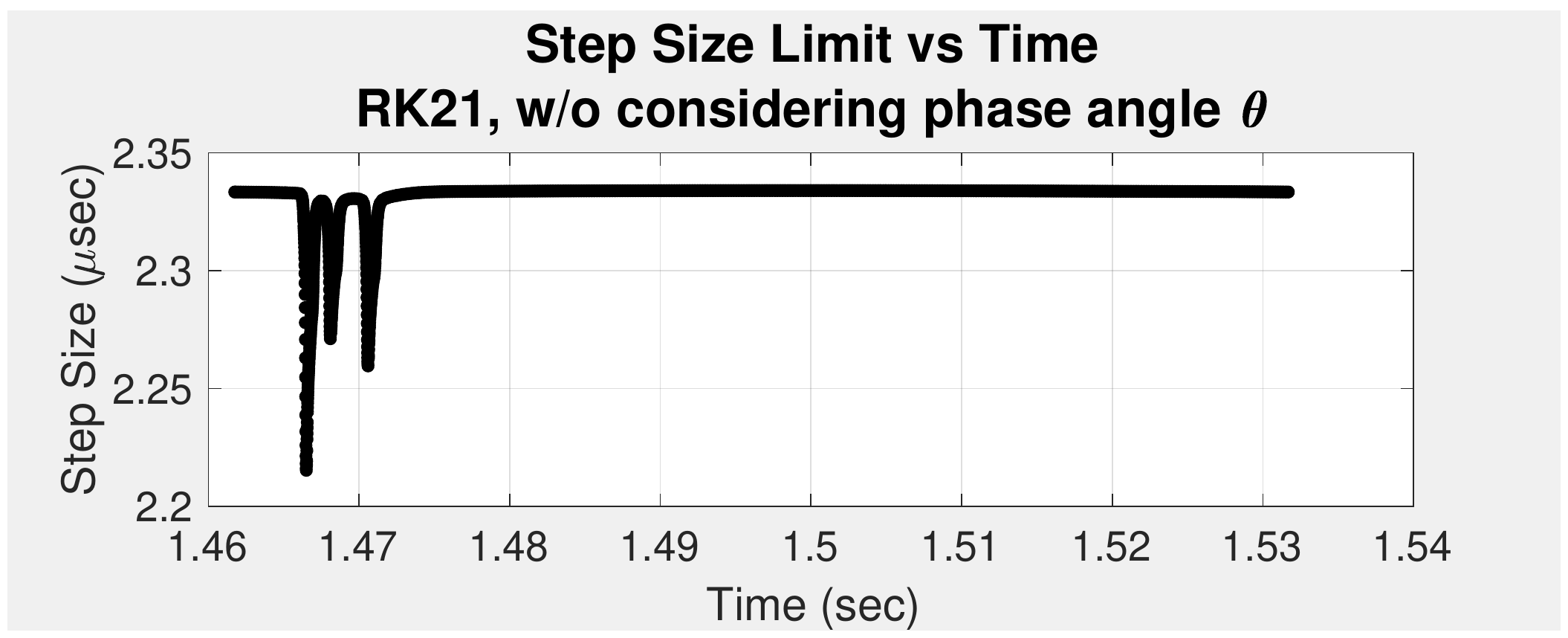} \\
\includegraphics[height=1.6in,width=3.7in]{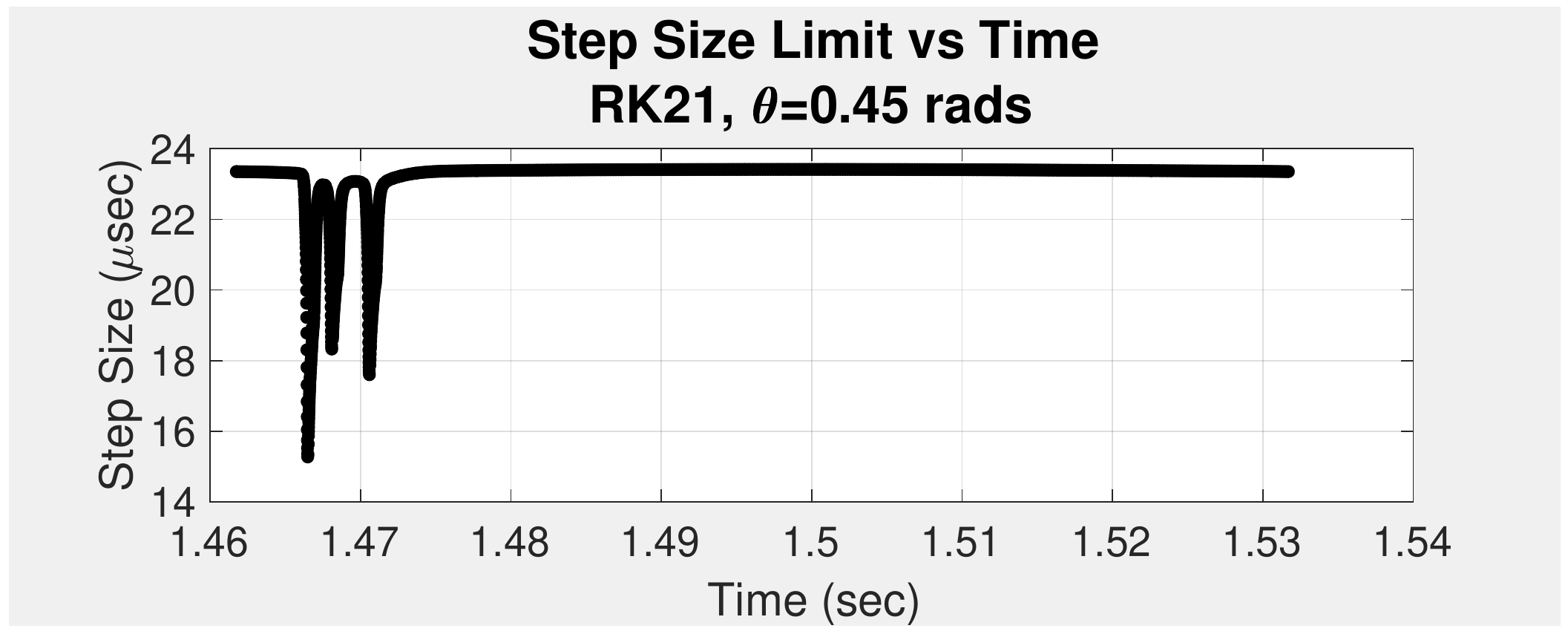} \\
\includegraphics[height=1.6in,width=3.7in]{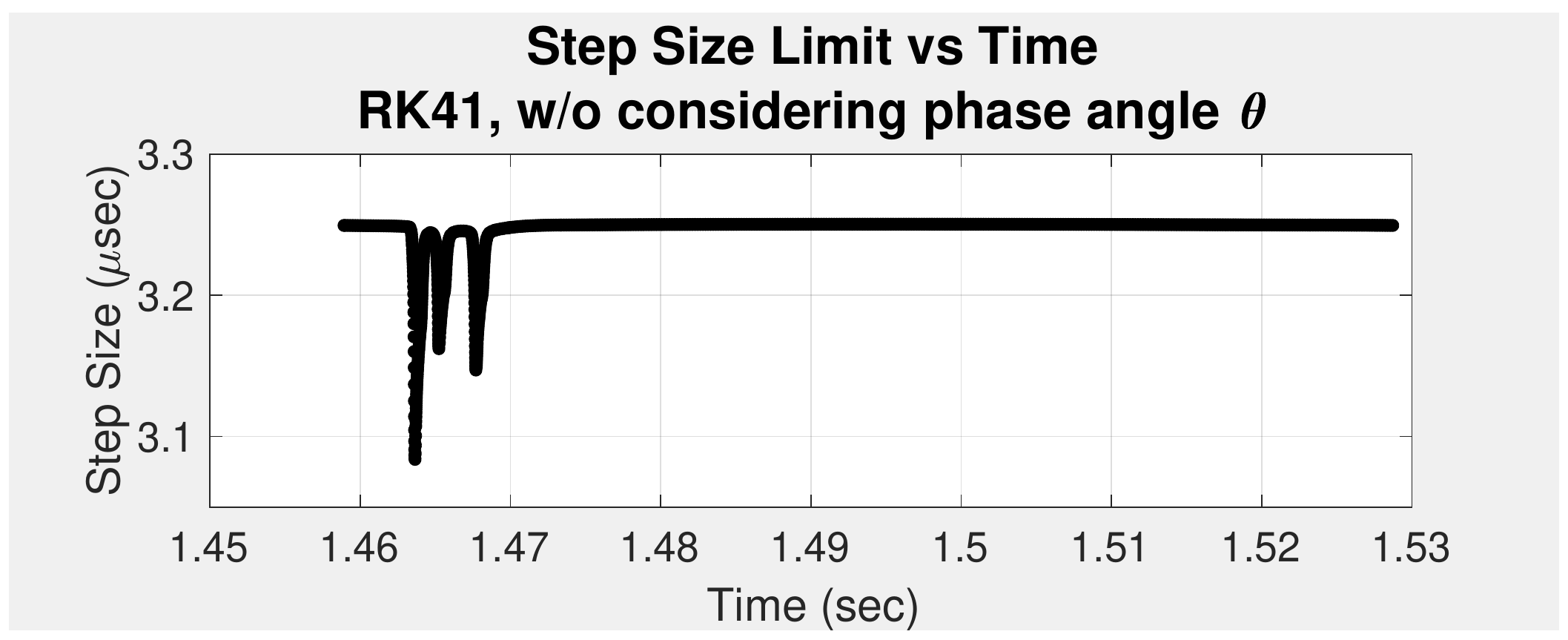} \\
\includegraphics[height=1.6in,width=3.7in]{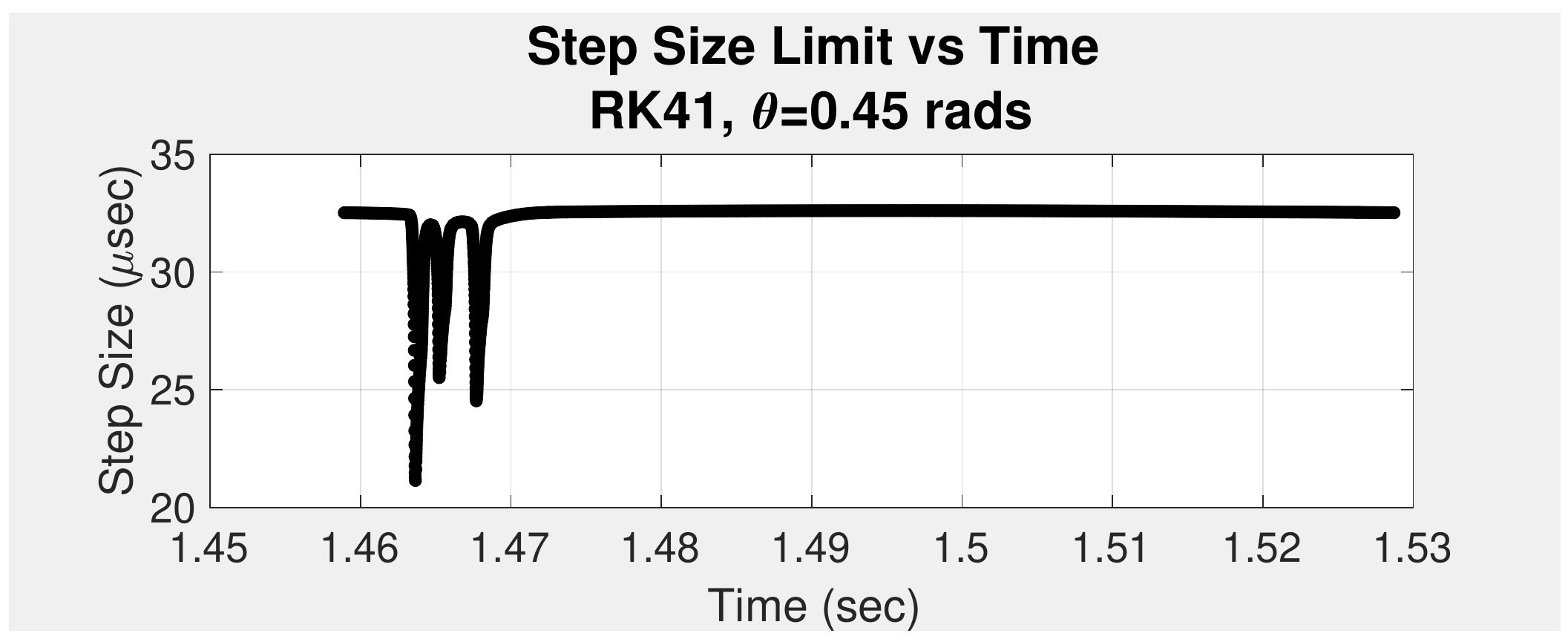}
\end{tabular}
\caption{Integration step size limits shown above for RK21 and RK41 predicted two ways. The first using RK stability analysis from \cite{Butcher} 
for a true ODE, not a spatially discretized PDE, does not consider the plane wave's phase angle. The second, Von Neumann stability analysis 
applied to RK methods cast as a quasi-FD methods, is in terms of the phase angle. As shown above, and in Figure \ref{fig:successfulStepSizes}, 
step size limits from Von Neumann analysis match those observed for RK21 and RK41, $14\mu$sec and 
$20\mu$sec respectively.}
\label{fig:RKstepsizeLimits}
\end{figure}

\begin{figure}[h!]
\centering
\setlength{\tabcolsep}{0.7mm}
\begin{tabular}{l}
\includegraphics[height=3.1in,width=4.0in]{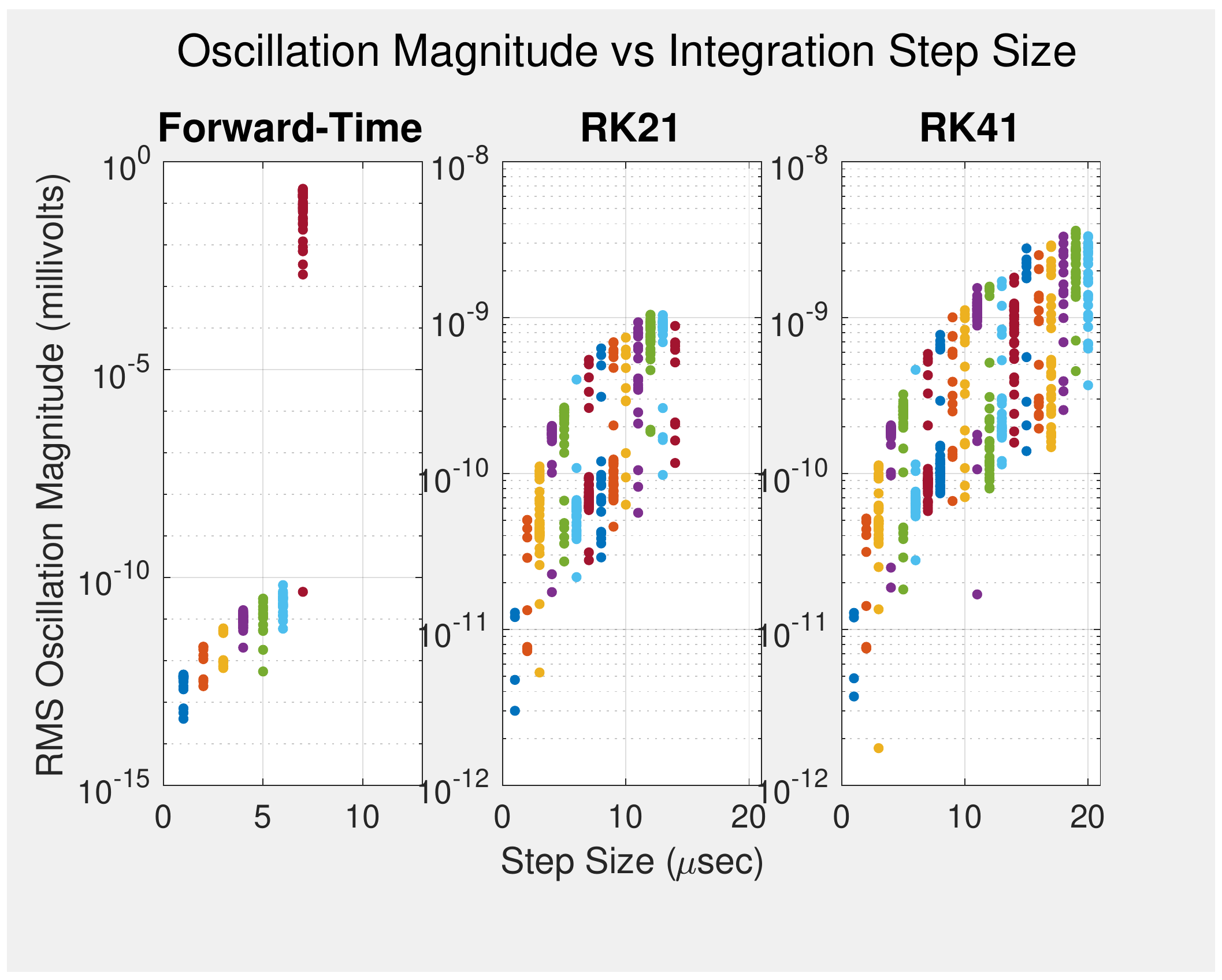}
\end{tabular}
\caption{Oscillation magnitudes versus step size for FTCS, RK21 and RK41. Each data point represents one of, as many as, $42$ AP cycles.}
\label{fig:fE2OToscVsTimestep}
\end{figure}

\begin{figure}[h!]
\centering
\setlength{\tabcolsep}{0.7mm}
\begin{tabular}{l}
\includegraphics[height=3.2in,width=2.5in]{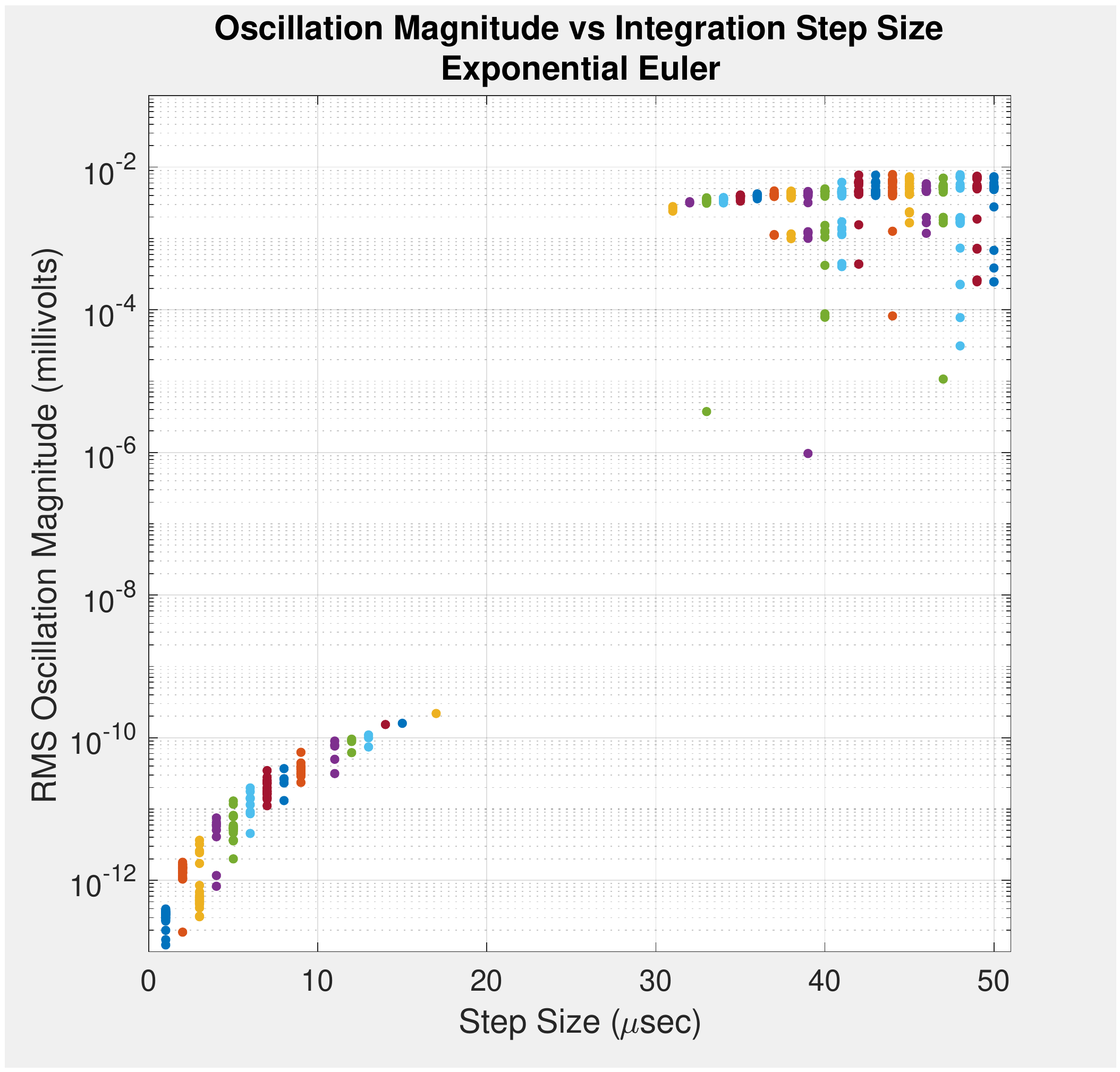}
\includegraphics[height=3.2in,width=2.5in]{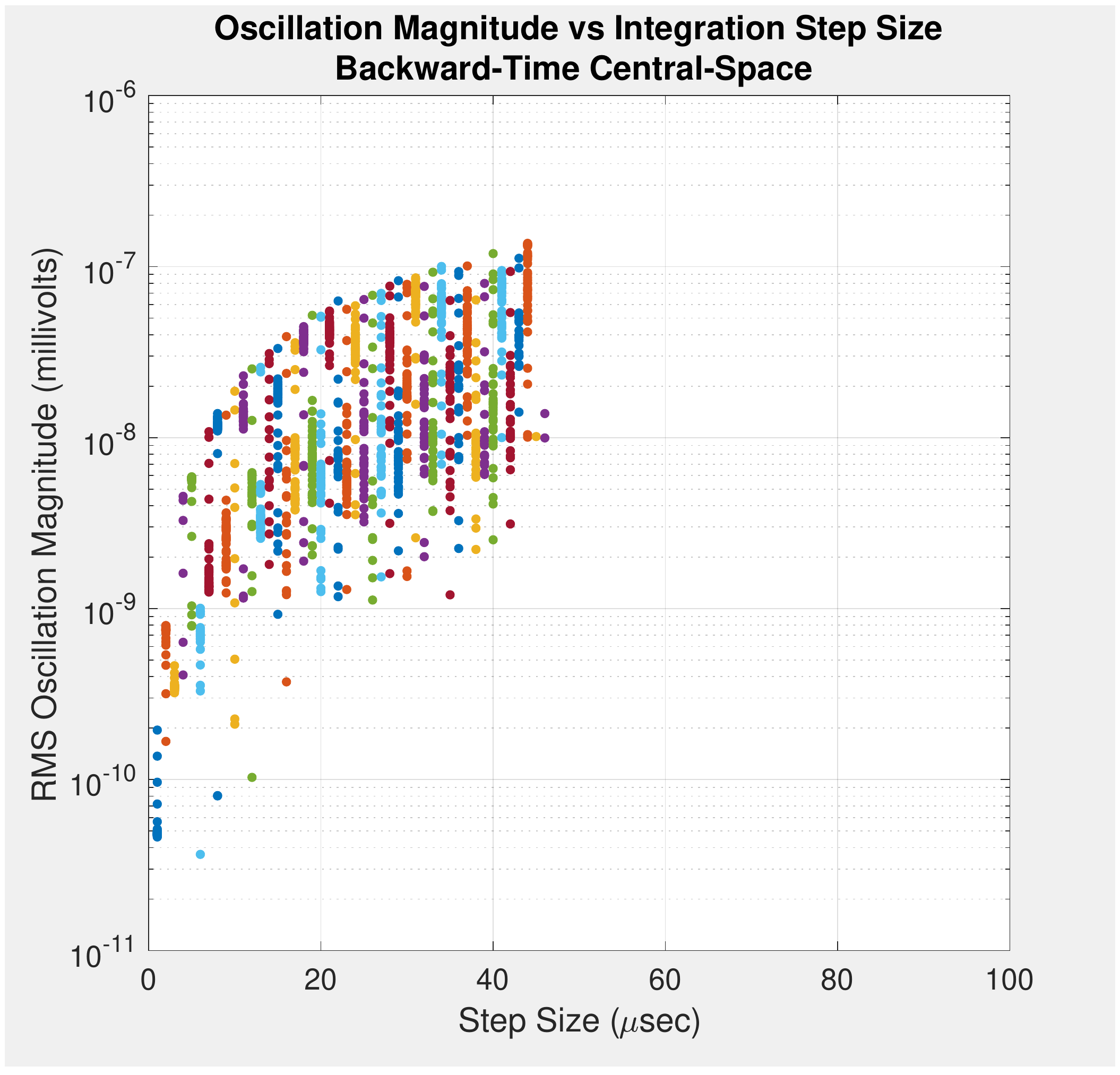}
\end{tabular}
\caption{Oscillation magnitudes versus step size for exponential Euler and BTCS. Each data point represents one of, as many as, $42$ AP cycles.}
\label{fig:eE_BTCSoscVsTimestep}
\end{figure}

\begin{figure}[h!]
\centering
\setlength{\tabcolsep}{0.7mm}
\begin{tabular}{l}
\includegraphics[height=3.2in,width=3.3in]{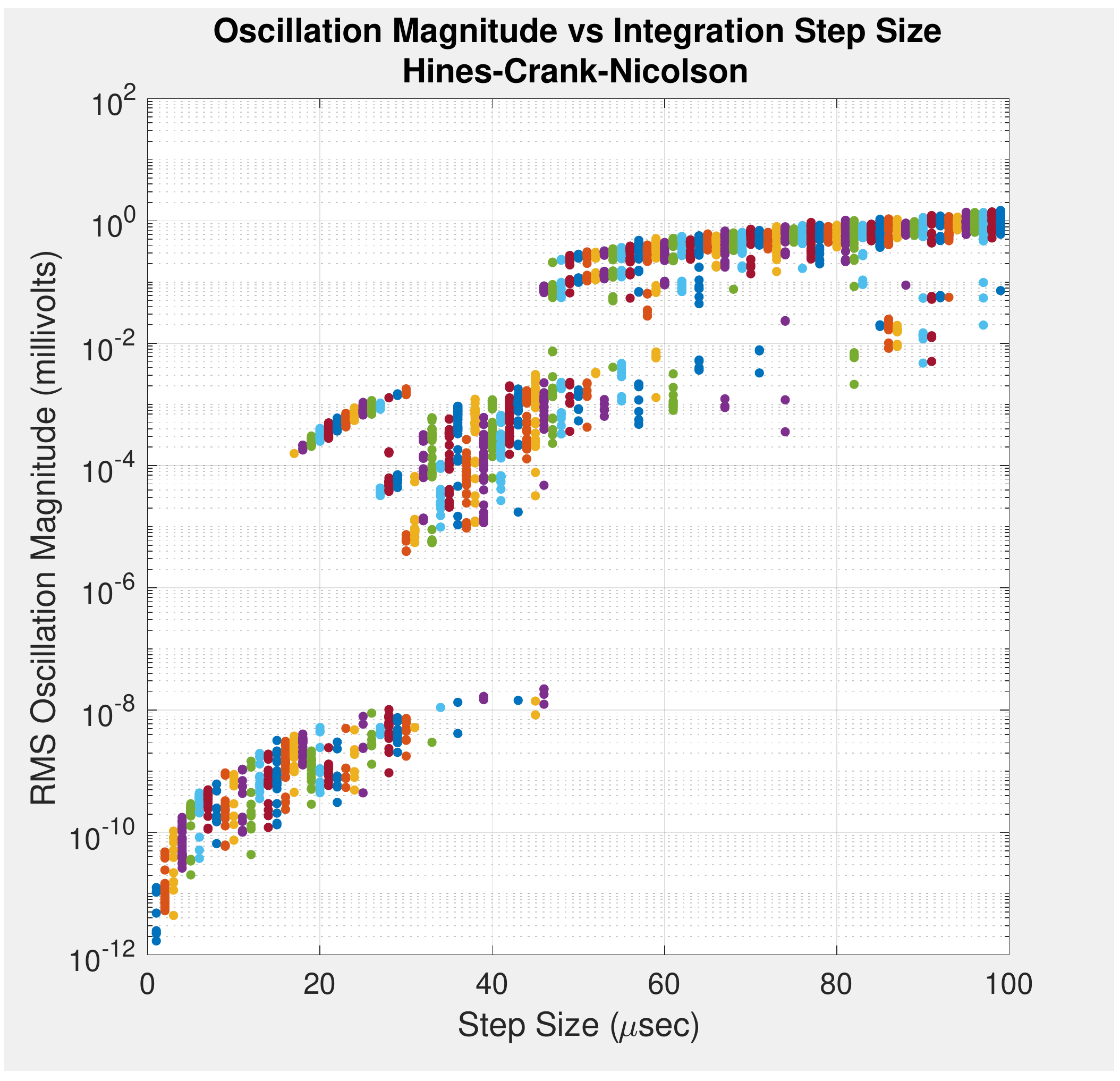} %\\
\end{tabular}
\caption{Oscillation magnitudes versus step size for HCN. Each data point represents one of $42$ AP cycles.}
\label{fig:hCNoscVsTimestep}
\end{figure}

\begin{figure}[h!]
\centering
\setlength{\tabcolsep}{0.7mm}
\begin{tabular}{l}
\includegraphics[height=3.2in,width=3.3in]{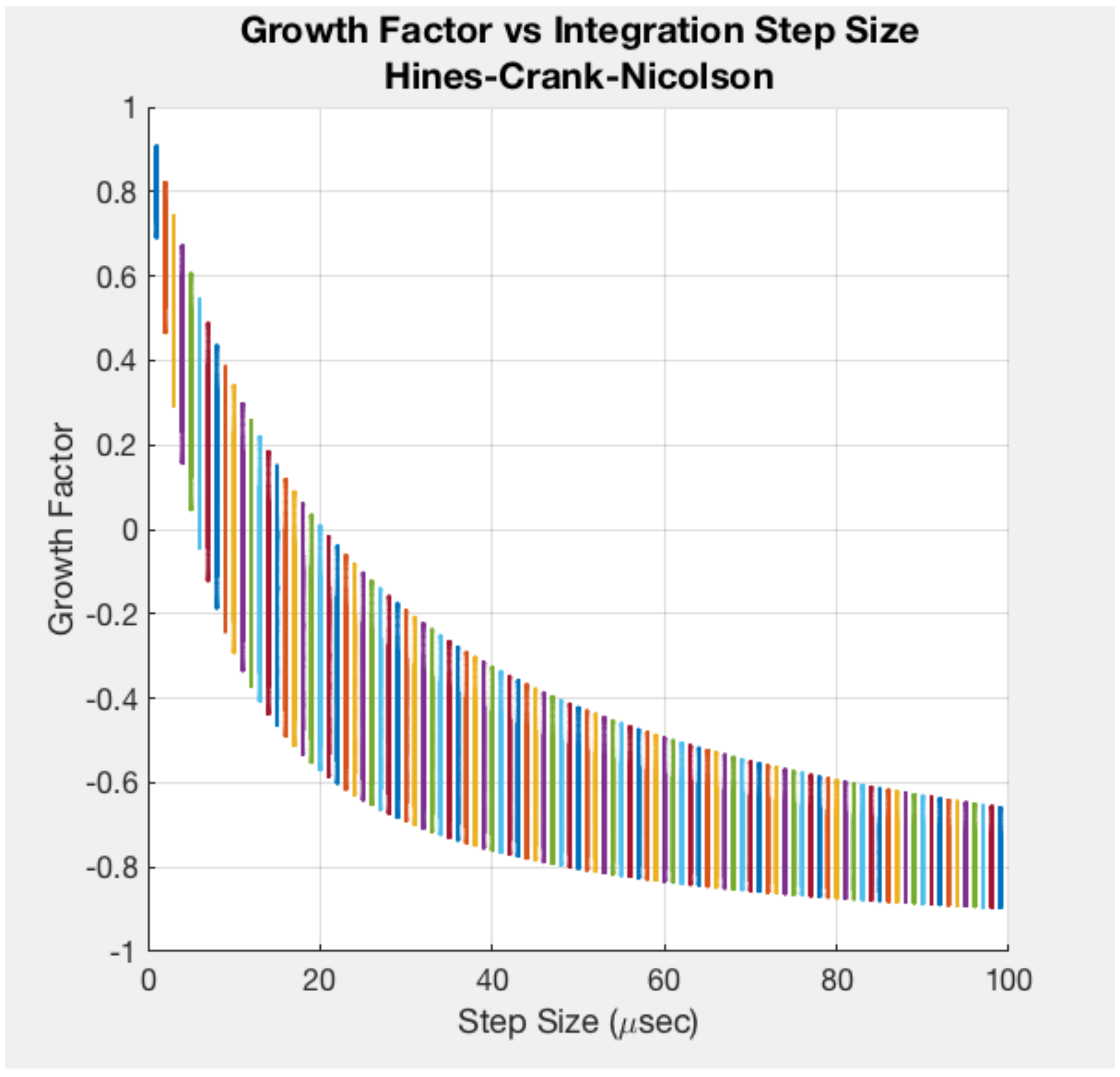} \\
\end{tabular}
\caption{Span of Von Neumann growth factor vs step size for HCN.}
\label{fig:hCNGrwthFctrVsTimestep}
\end{figure}

\begin{figure}[h!]
\centering
\setlength{\tabcolsep}{0.7mm}
\begin{tabular}{l}
\includegraphics[height=6.0in,width=4.5in]{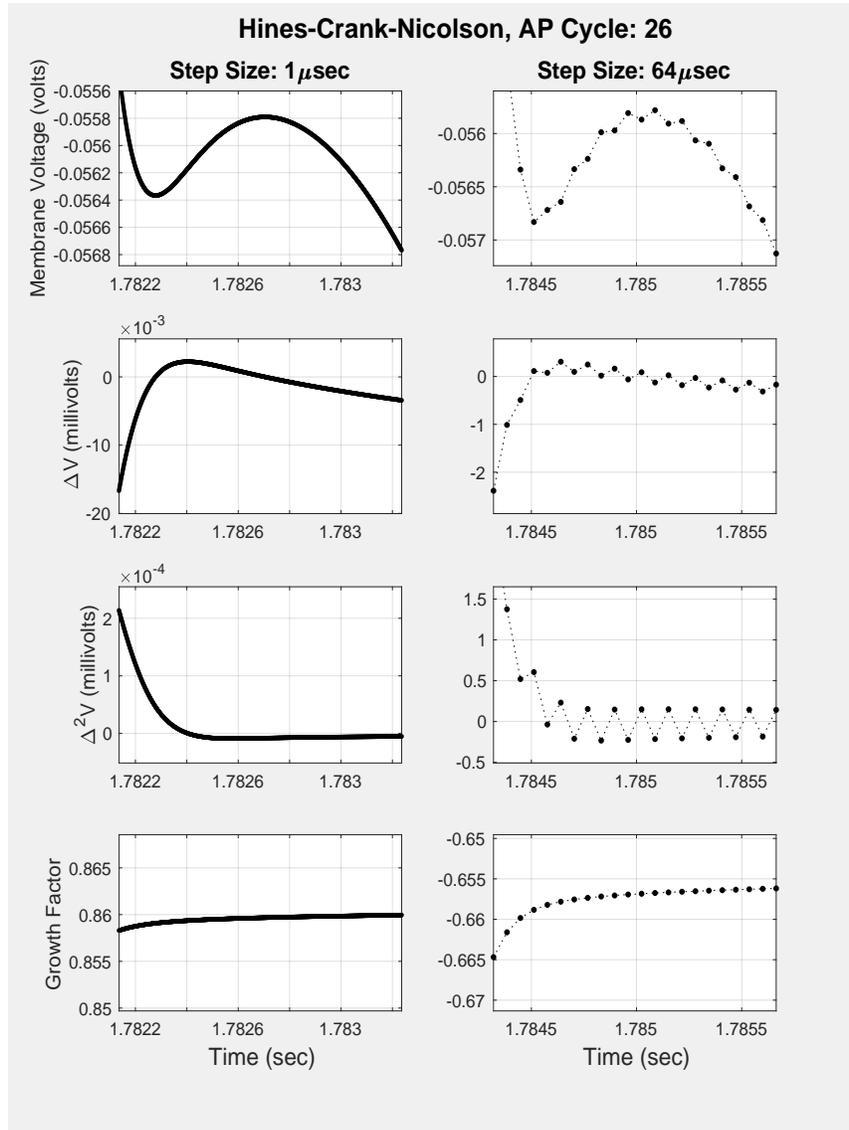}
\end{tabular}
\caption{Comparison of the ADP phase membrane voltage, $V$, the first and second undivided differences of membrane voltage, and 
the Von Neumann growth factor for integration step sizes $1$ and $64\mu$s. The onset of oscillation by the time step size has 
increased to $64\mu$s is also demonstrated by evidence of ADP oscillation in the lower right graph of Figure \ref{fig:p23rsAPclassBE_HCN}, 
by RMS oscillation magnitudes in Figure \ref{fig:hCNoscVsTimestep} and corroborated by the trend of growth factor illustrated in Figure 
\ref{fig:hCNGrwthFctrVsTimestep}.}
\label{fig:ADP}
\end{figure}

\begin{figure}[h!]
\centering
\setlength{\tabcolsep}{0.7mm}
\begin{tabular}{l}
\includegraphics[height=6.0in,width=4.5in]{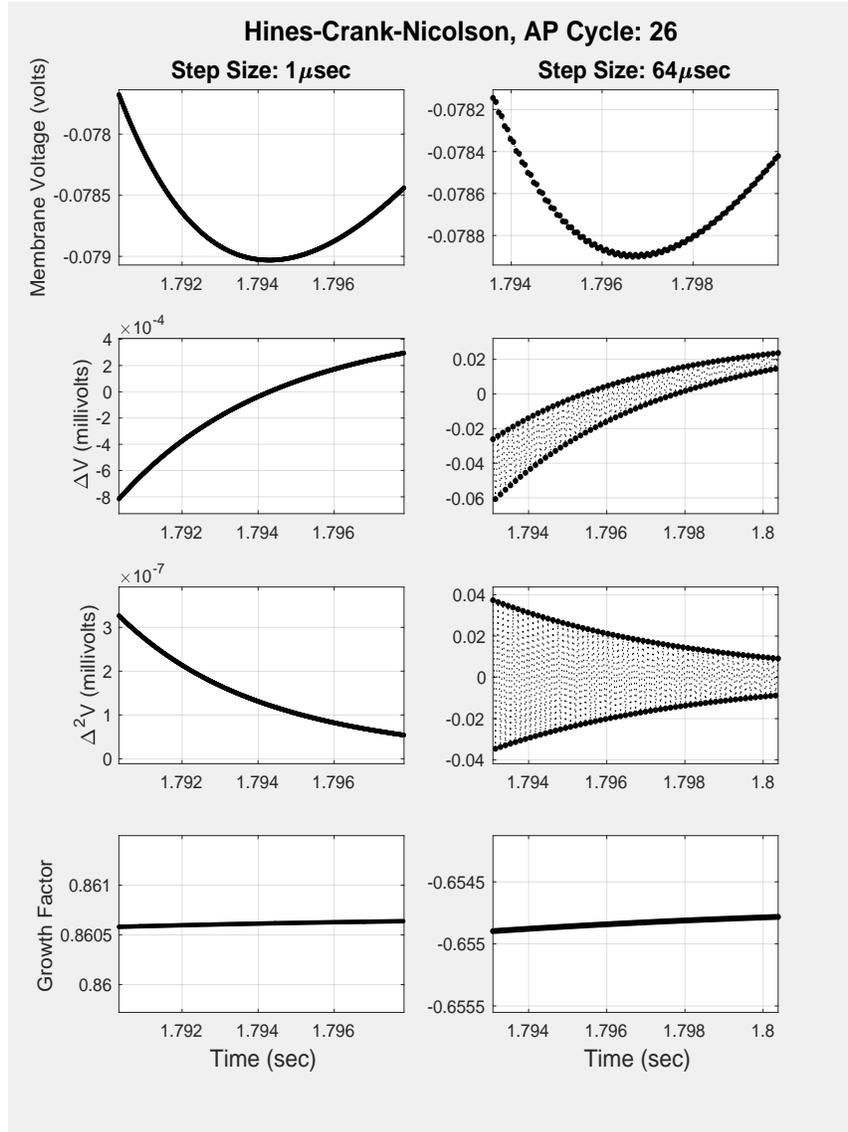}
\end{tabular}
\caption{Comparison at maximum polarization of membrane voltage, $V$, the first and second undivided differences of membrane voltage, 
and the Von Neumann growth factor for integration step sizes $1$ and $64\mu$s. The onset of oscillation by the time step size has increased 
to $64\mu$s is demonstrated by RMS oscillation magnitudes in Figure \ref{fig:hCNoscVsTimestep} and corroborated by the trend of growth factor 
illustrated in Figure \ref{fig:hCNGrwthFctrVsTimestep}.}
\label{fig:maxPolarization}
\end{figure}

\begin{figure}[h!]
\centering
\setlength{\tabcolsep}{0.7mm}
\begin{tabular}{l}
\includegraphics[height=6.0in,width=4.5in]{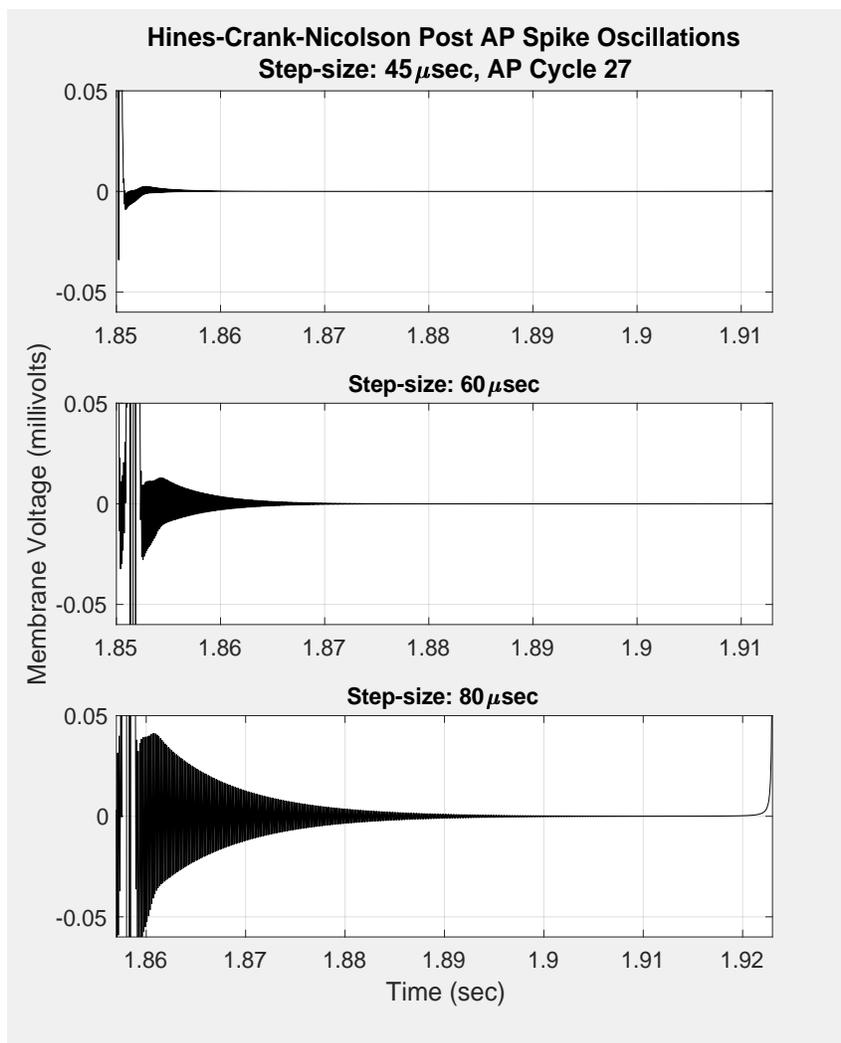}
\end{tabular}
\caption{HCN oscillation magnitude increases with step size, as predicted in \cite[ch.4]{CarnevaleHines}, and decays before the next AP cycle begins.}
\label{fig:hCNOscDecay}
\end{figure}

Figure \ref{fig:successfulStepSizes} illustrates stable integration step size intervals observed for each method. While the exponential Euler 
method was always stable, its waveforms were sufficiently inconsistent to prevent comparison with waveforms presented by other methods. 
So this method was only exercised for step sizes from $1-50\mu$s. \\ \\
Figure \ref{fig:fEstepSizeLimit} shows the correspondence between membrane voltage, the model's plane wave's phase angle and the Von Neumann 
step size limit for FTCS. Spectral centroid of the phase angle was computed, as described in Section \ref{subsection:SpecAnalysis}, from the FFT of 
membrane voltages of compartments from the tip of any distal dendrite to the tip of the axon. These thirteen values were then padded with nineteen 
zeros to increase low frequency resolution. \\ \\
Figure \ref{fig:RKstepsizeLimits} shows integration step size limits for RK21 and RK41 predicted two ways. The first using RK stability analysis from 
\cite{Butcher} for a true ODE, not a spatially discretized PDE, does not consider the plane wave's phase angle. The second, Von Neumann 
stability analysis applied to RK methods cast as a quasi-FD methods, is in terms of the phase angle. Figure \ref{fig:successfulStepSizes} and Figure 
\ref{fig:RKstepsizeLimits} show that step size limits from Von Neumann analysis match those we observed for RK21 and RK41, $14\mu$sec 
and $20\mu$sec respectively. \\ \\
Magnitude and location of oscillations presented by each integration method, as defined in Section \ref{sec:OscDet}, are shown in Figures 
\ref{fig:fE2OToscVsTimestep}--\ref{fig:hCNoscVsTimestep}. Note that FTCS method oscillation magnitudes grow significantly one or two 
microseconds in step size before stability is lost. Exponential Euler's oscillations approach $10^{-2}m$volts shortly after integration step size 
increases beyond $30\mu$s. \\ \\
While HCN was stable for all step sizes from $1$ to $99\mu$s, small amplitude decaying oscillations grew in duration and magnitude as step 
sizes grew larger than $50\mu$s as illustrated in Figure \ref{fig:hCNoscVsTimestep}. HCN's Von Neumann growth factor, spanning one AP 
cycle, is plotted versus integration step size in Figure \ref{fig:hCNGrwthFctrVsTimestep}. Key oscillation locations are illustrated in Figures 
\ref{fig:ADP}--\ref{fig:maxPolarization}. Figure \ref{fig:hCNOscDecay} illustrates that while oscillations are decaying, their magnitude increases 
with step size. Videos in the supplemental resources show HCN's oscillations evolving for step sizes $1-99\mu$s on three intervals, AP spikes, 
ADP pulse and maximum membrane polarization separating AP cycles.   

\begin{figure}[h!]
\centering
\setlength{\tabcolsep}{0.7mm}
\begin{tabular}{l}
\includegraphics[height=1.8in,width=4.5in]{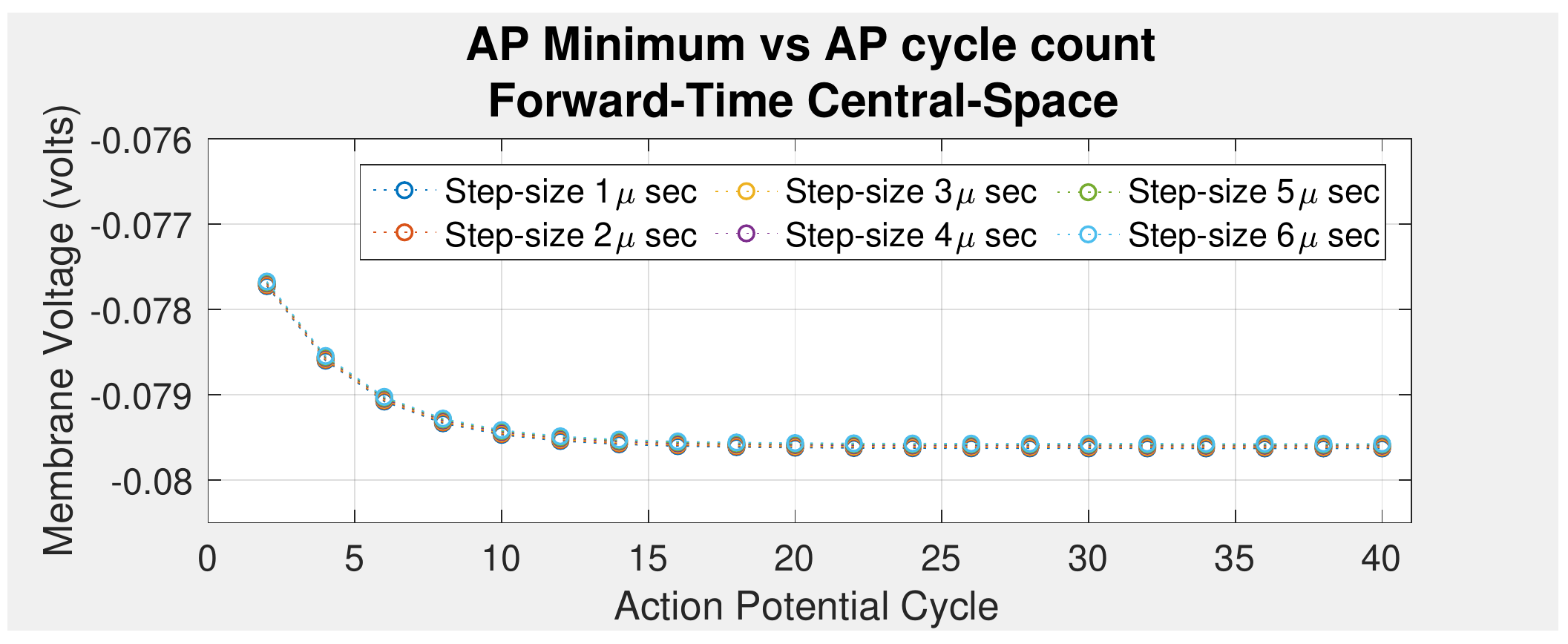} \\
\includegraphics[height=1.8in,width=4.5in]{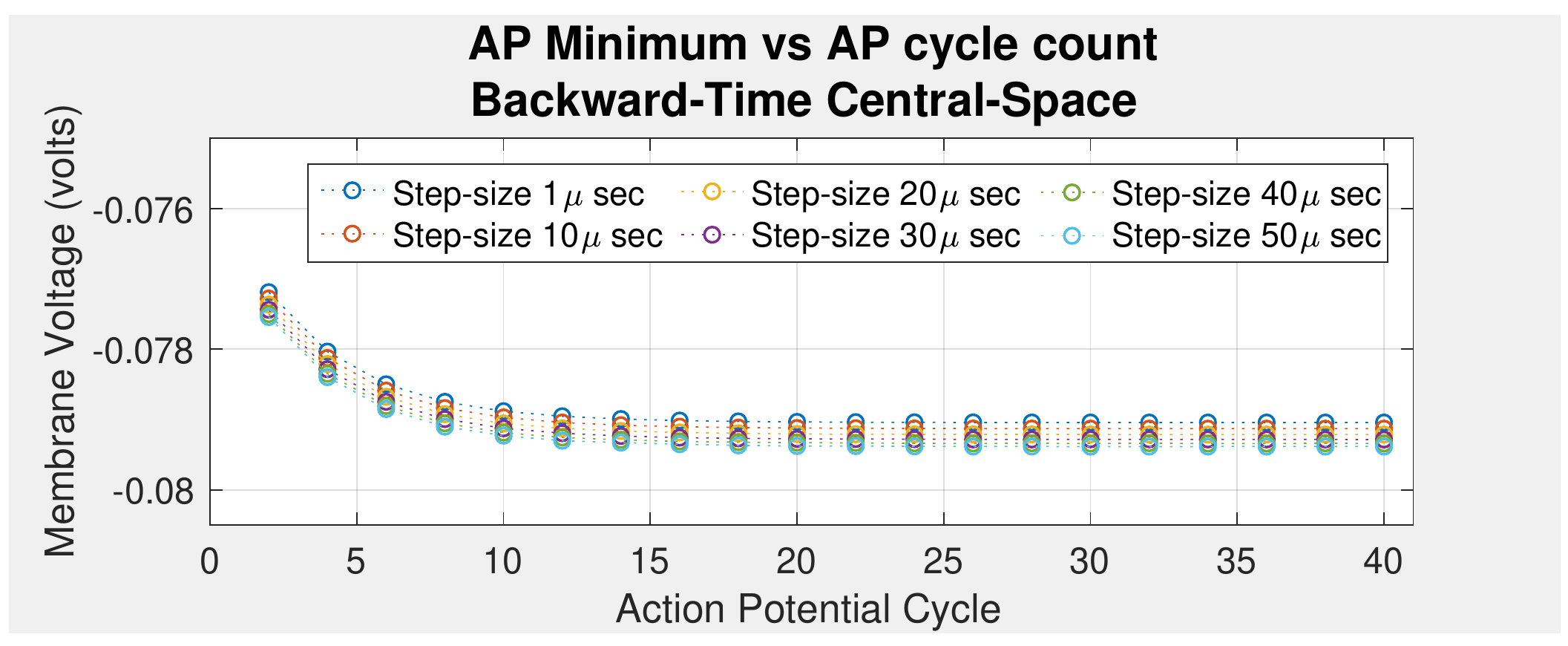} \\
\includegraphics[height=1.8in,width=4.5in]{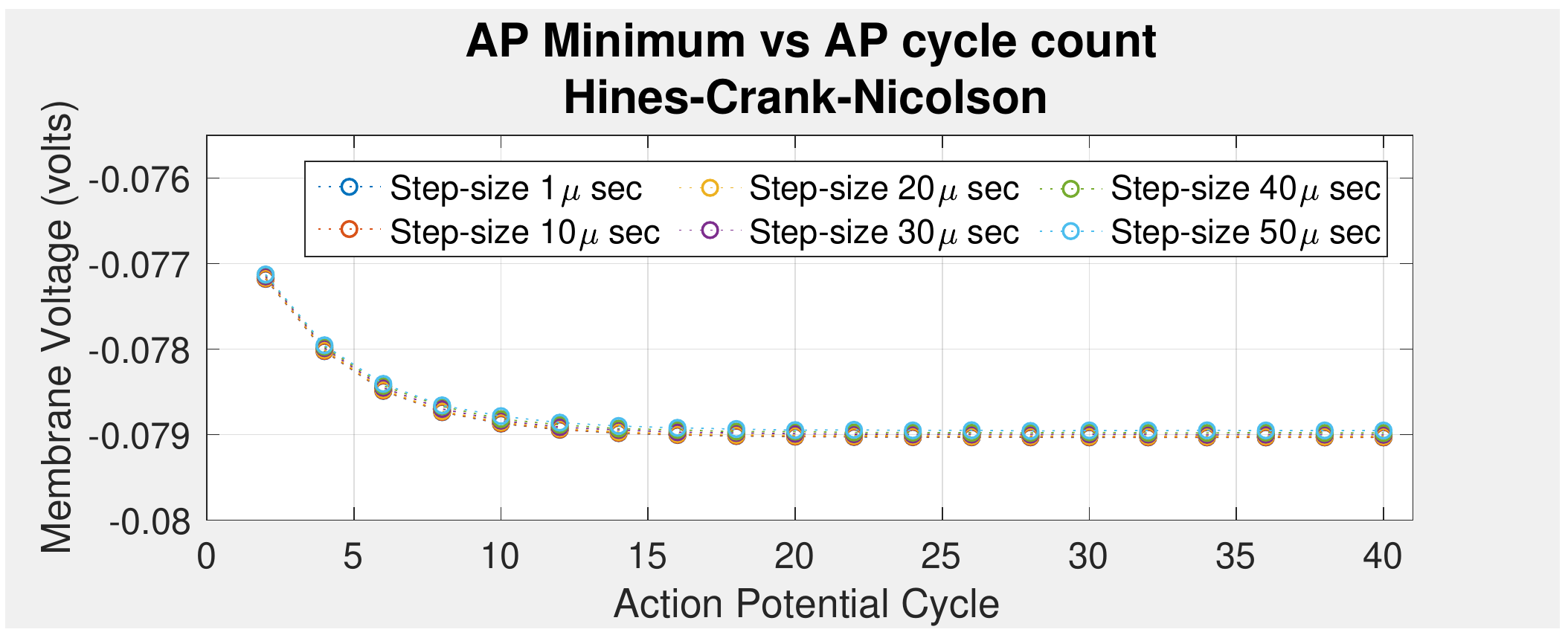} \\
\end{tabular}
\caption{Relative minimums between AP cycles stabilize by the twentieth AP cycle.}
\label{fig:minsVapidx}
\end{figure}

\begin{figure}[h!]
\centering
\setlength{\tabcolsep}{0.7mm}
\begin{tabular}{l}
\includegraphics[height=1.8in,width=4.5in]{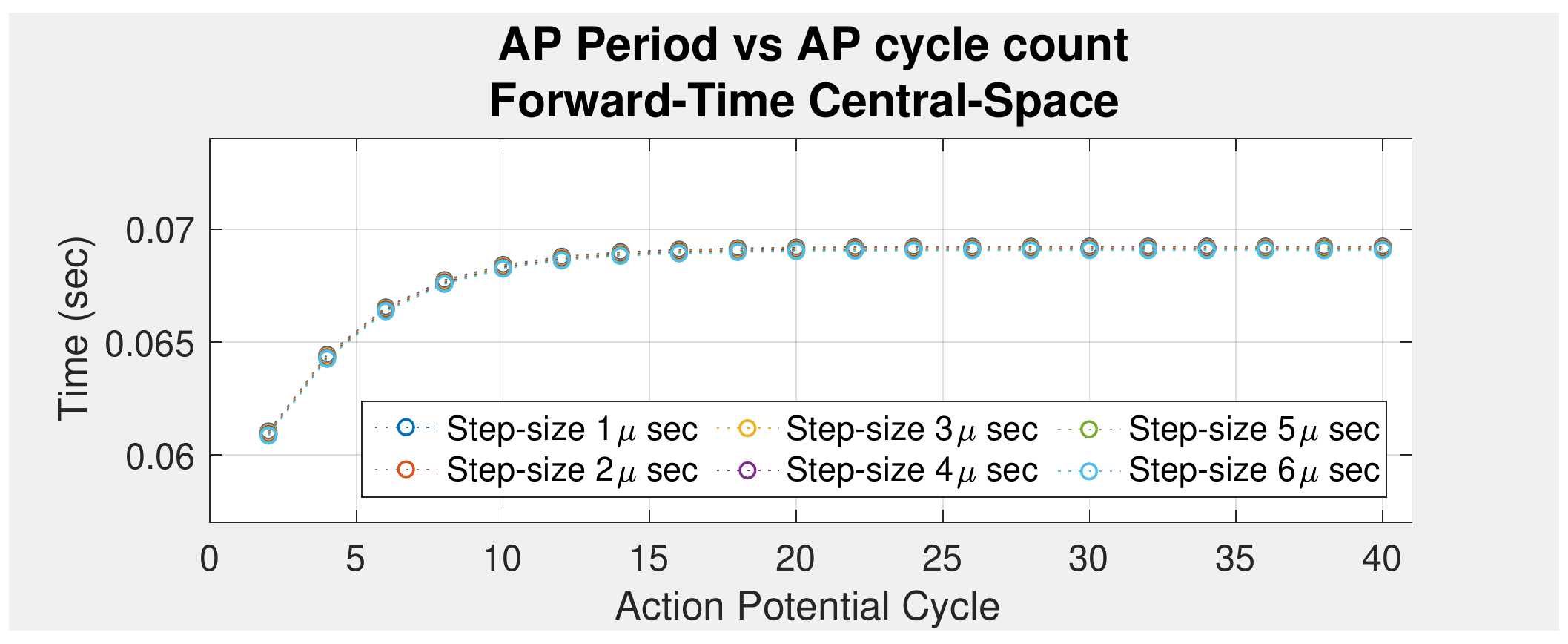} \\
\includegraphics[height=1.8in,width=4.5in]{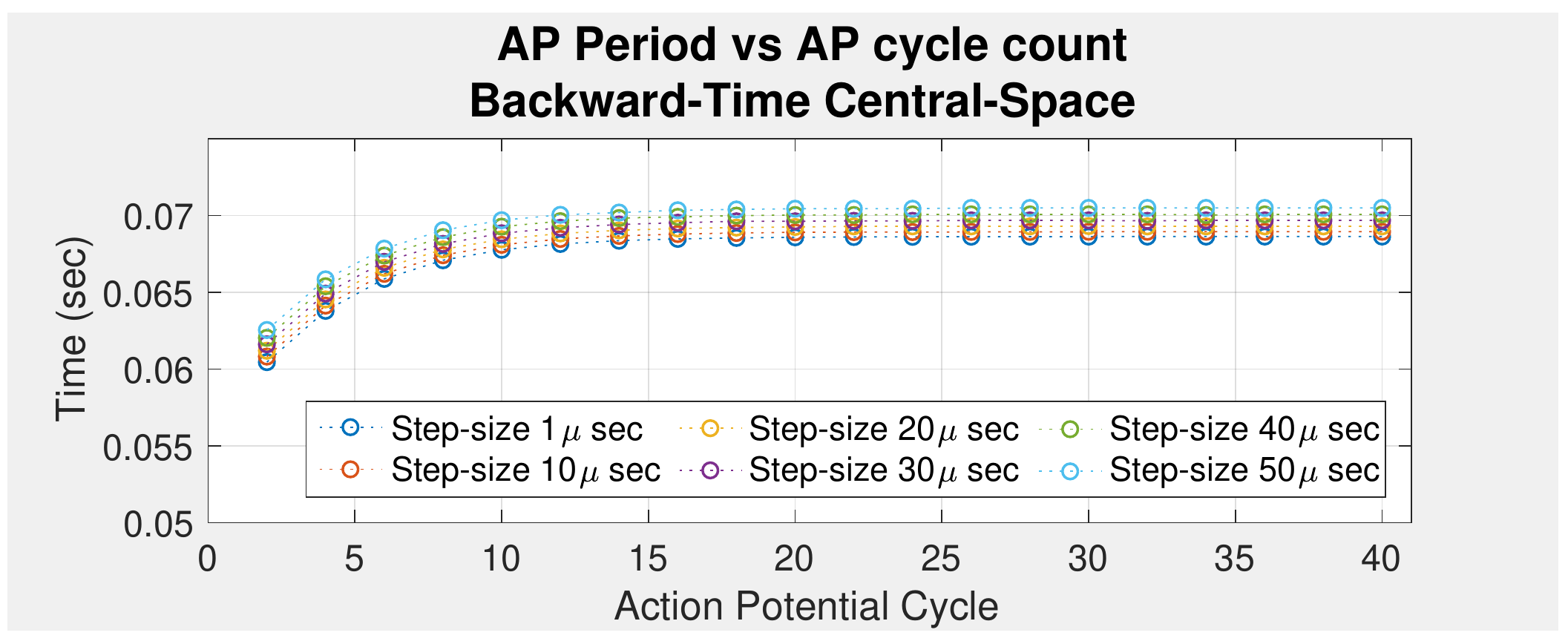} \\
\includegraphics[height=1.8in,width=4.5in]{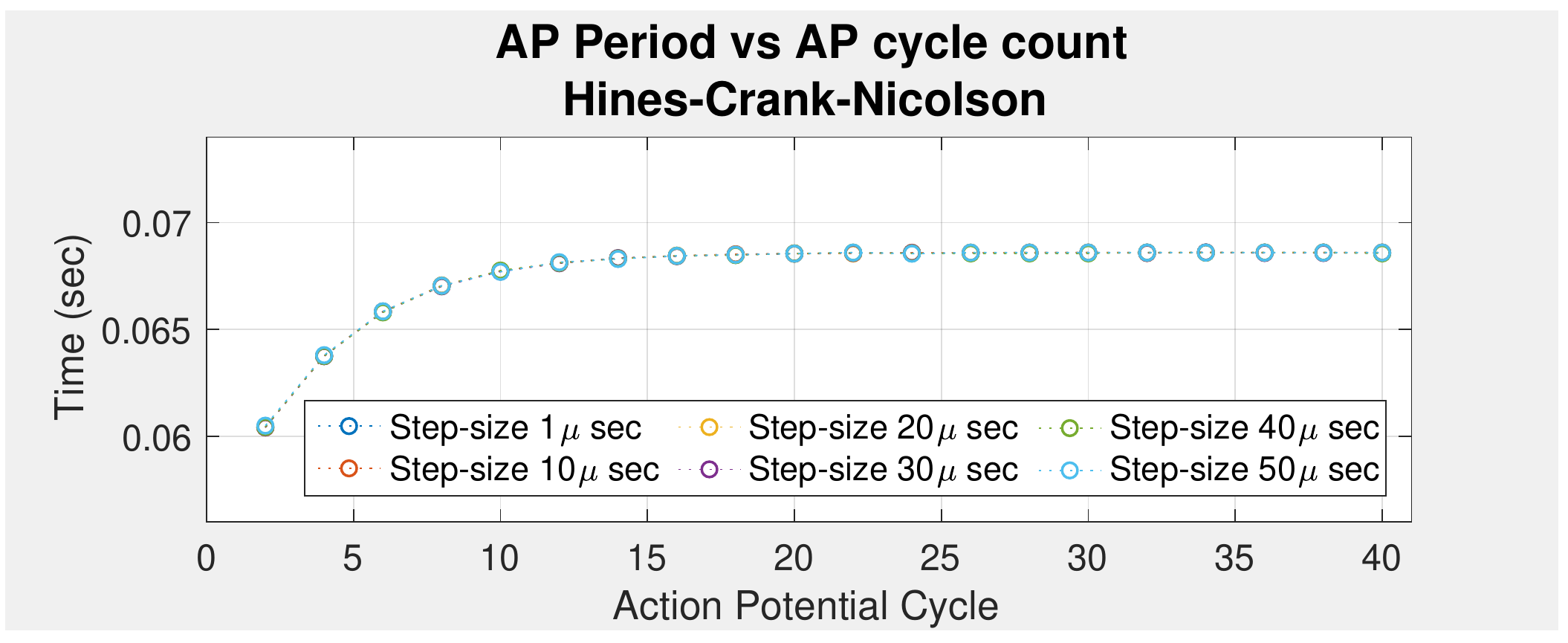} \\
\end{tabular}
\caption{Action potential (AP) periods stabilize by the twentieth AP cycle.}
\label{fig:periodsVapidx} 
\end{figure}

\clearpage


\clearpage

\clearpage

\begin{figure}[h!]
\centering
\setlength{\tabcolsep}{0.7mm}
\begin{tabular}{l}
\includegraphics[height=6.0in,width=4.5in]{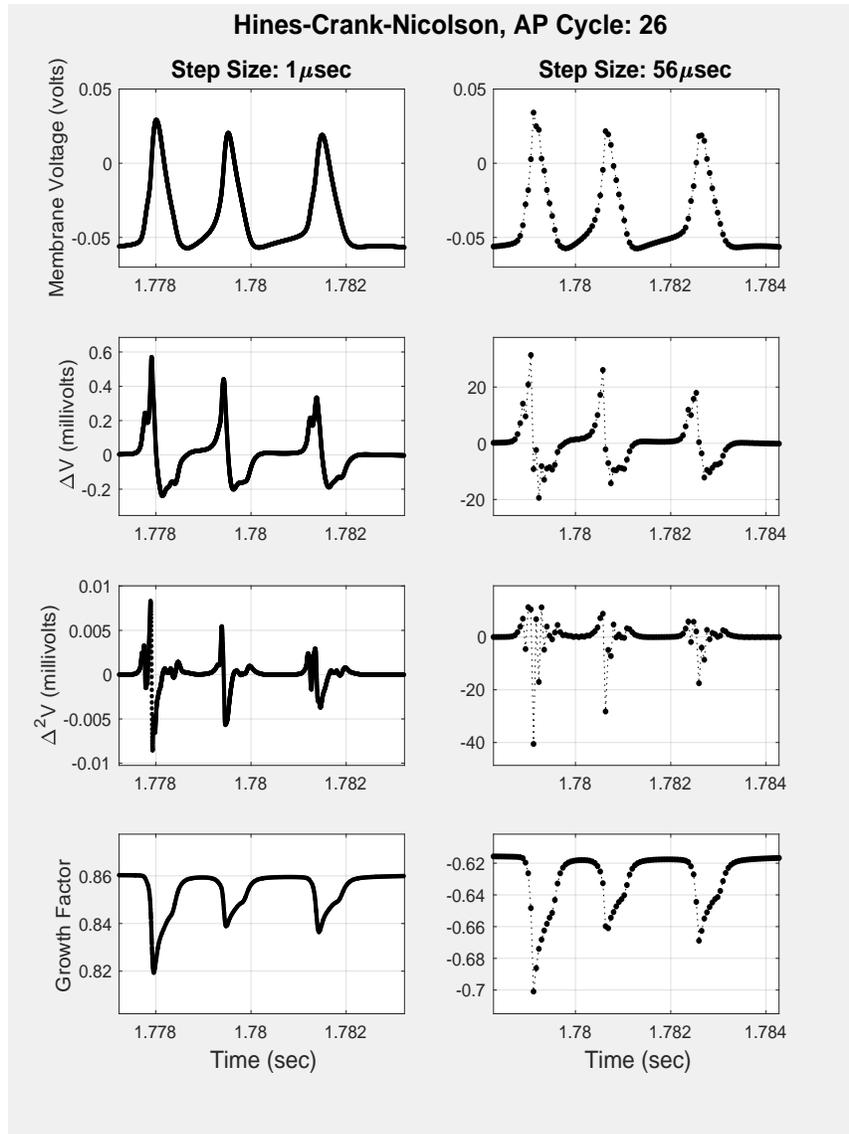}
\end{tabular}
\caption{Comparison of spiking phase membrane voltage, $V$, the first and second undivided differences of membrane voltage, and 
the Von Neumann growth factor for integration step sizes $1$ and $51\mu$sec. The onset of oscillation by the time step size has increased 
to $51\mu$sec is also demonstrated by HCN's AP maximums shown in Figure \ref{fig:somaStatsWithStd}, by RMS oscillation magnitudes 
in Figure \ref{fig:hCNoscVsTimestep} and justified by the trend of growth factor illustrated in Figure \ref{fig:hCNGrwthFctrVsTimestep}.}
\label{fig:Spiking}
\end{figure}

\clearpage
\section{Summary and Discussion}
\label{sec:SumAndDiscuss} 
Our model's computational sequence is composed of two sub-steps, the first to integrate each compartment's eleven ion channel state ODEs, 
represented in equation (\ref{eqn:gateStateODE}), and second to integrate each compartment's membrane voltage PDE, in equation
(\ref{CableEquationRoot}), as a function of the updated channel states. The solution to our model's dynamical system is periodic, a sequence of AP 
cycles, each with its own sub-sequence of three states, Spiking, Re-polarization and a long slow Depolarization, as shown in Figure 
\ref{fig:idealAPcyclesHCNsts1}. We expect numerical methods integrating our model's dynamical system to maintain periodicity and synchronization 
of all channel state cycles in order to ensure consistent AP cycle waveforms and statistics such as relative minimums and maximums, period lengths, 
spike phases and multiplicities. \\ \\
One of our first observations is that roughly twenty AP cycles were required before state statistics reached a mature equilibrium, regardless of 
integration method, as shown in Figures \ref{fig:minsVapidx}-\ref{fig:periodsVapidx}. Membrane ion gate phase portraits included in the supplemental 
resources also support this observation. These are examples of common characteristics that combine to delay the onset of a mature dynamical system 
state. Such 'warm-up' time could be reduced by saving a model's state variables at the end of a deterministic run-up to a mature state for use as initial 
values prior to varying conditions and introducing stochastic regimes. For this reason research articles would be more helpful if they indicate how far 
into a simulation their data was collected \cite{Ascoli10,Sousa15,Kuznet10}. \\ \\
The purpose of our study, a comparison of numerical integration methods, prompts the question, ``Why should we expect any method-dependent 
integration differences at all?'' While an integration method's truncation error is not the same as the inaccuracy of the solution it generates, terms in 
the former help show why the later is different for each integration method. The first and second terms of truncation errors, listed in Table 
\ref{tab:schemeTruncError} for the membrane voltage integration, are second and third derivatives of the membrane voltage which change throughout 
the AP cycle and ensure error differences between methods. The following two sections describe integration method-specific differences in accuracy 
and stability. 
\subsection{Accuracy}
\label{subsec:AccuracyDiscussion}
Perhaps the simplest generalization worth making is that every explicit ODE method integrating a multi-compartment neural model sacrifices accuracy 
by failing to account for the changing spatial derivative which is no different than assuming it remains constant for the duration of each step. Any hope 
of recovering this loss by replacing one explicit integrator with another of higher order, given our selection of integrators, is lost either because instability 
takes over at such small step sizes, as shown in Figure \ref{fig:successfulStepSizes}, or the cable equation renders RK methods first-order, a matter we 
discuss below. \\ \\  
Our study assessed integration method accuracy in the following ways. First, we computed the RMS difference of AP waveforms relative to the 
reference waveform integrated by the same method at $1\mu$s step size, as described in Section \ref{sec:AccMeasure}. These measurements 
of difference, plotted versus step size for each method in Figure \ref{fig:methodErrVsts} corroborate most truncation error orders presented in Table \ref{tab:schemeTruncError}. The abrupt shocks in the same figure presented by BTCS, RK21 and RK41 %, likely due to bifurcation incursions, 
are discussed below.  \\ \\
As summarized in Table \ref{tab:schemeTruncError} and confirmed in Figure \ref{fig:methodErrVsts}, RK21 and RK41 both present, at first glance, 
linearly increasing first-order error. Of course, RK21 and RK41 are indeed second- and fourth-order methods, respectively, when integrating each 
compartment's ion channel state ODEs. It is when integrating the discretized cable equation, as shown in Section \ref{subsubsec:RK21truncErr}, 
that truncation errors become first-order. In contrast, HCN presents second-order error growth, concave up as step size increases. \\ \\
There remains one unusual nonlinear anomaly in Figure \ref{fig:methodErrVsts} affecting BTCS and RK41. If FTCS and RK21 had remained stable 
long enough their data may also have shown a similar trend. The anomaly is that their error plots present an error order less than unity evidenced 
by the error curve's negative concavity. Perhaps BTCS was most affected as it was the only first-order method to remain stable. The second most 
significantly affected was RK41, also a, albeit situationally rendered, first-order method. We suggest the most likely reason for this anomaly is the 
same reason HCN does {\it{not}} present it. That is, to minimize error accrued by integrating highly nonlinear HH type cable equations, HCN has 
always prescribed staggering gate state integrations by one half-step relative to the membrane voltage integration. By not staggering channel state 
updates between membrane voltage updates, all other methods were more likely to accrue more integration error than did HCN. \\ \\
The second way we assessed accuracy of numerical integration methods was with AP waveform statistics. These statistics are a measure of an 
integration method's ability to consistently control ion channel constituent gate states periodicity over broad step size intervals. Without such 
consistency the reliability of an adaptive step method would be uncertain. Figure \ref{fig:somaStatsWithStd} illustrates two method-specific 
inconsistencies which present alone and together. We first thought that the meandering mean and discontinuous variation of AP maximums 
presented by HCN after $50\mu$s step sizes was one of these inconsistencies. In fact, the source of this artifact is a decaying oscillation and 
is discussed below. \\ \\
The first method-specific inconsistency seen in Figure \ref{fig:somaStatsWithStd} is a smooth nonzero rate-of-change over increasing integration 
step size, error characteristic of first order methods. With the exception of erratic AP maximums after $50\mu$s step sizes, HCN presents perfect 
consistency without any change to AP minimums, maximums and periods. While not as perfectly controlled as HCN, BTCS statistics present a very 
small rate-of-change over the same interval and also present control of AP maximums after $50\mu$s step sizes superior to HCN. Remaining 
integrators, all explicit methods, present significant and sometimes discontinuous change to AP minimums, maximums and periods as a function 
of integration step size. \\ \\
% Removed ...
%A second type of inconsistency seen in Figure \ref{fig:somaStatsWithStd} is what would be expected when the dynamical system encounters a 
%bifurcation in state space \cite{Guck93, Wang07}. RK21 and RK41 both enter into, and then exit, a different state space from step 
%sizes $9-11\mu$sec, and $12-15\mu$sec respectively. Similarly, BTCS takes the dynamical system past a likely bifurcation near step size 
%$81\mu$sec and does not recover. Exponential Euler also presents similar behavior evidenced by abrupt changes to AP minimums and periods 
%at step sizes of $3-5\mu$sec, and often whenever step sizes are larger than $30\mu$sec. Corroborating evidence for these accuracy shocks
%can be seen in Figure \ref{fig:methodErrVsts}, AP spiking multiplicity variations illustrated in Figures 
%\ref{fig:p23rsAPclassFE_2OT_RK21RK41}--\ref{fig:p23rsAPclassExpoE} and exponential Euler's waveform plots in Figure 
%\ref{fig:corruptExpoEulerAPs}. Also supporting the hypothesis of bifurcation incursion are simulation videos included with the supplemental 
%materials for each method over all step sizes. These videos show sudden significant changes to the sequence of AP waveforms at the given step sizes. \\ \\
A second type of inconsistency seen in the top and bottom panes of Figure \ref{fig:somaStatsWithStd} presents as significant jumps in the data
and affects all methods but HCN. RK21 and RK41 both present these data shocks from step sizes $9-11\mu$s, and $12-15\mu$s respectively. 
While BTCS data shocks appear limited to $81-90\mu$s step sizes, the exponential Euler method is affected throughout much of its step size 
testing interval, $1-50\mu$s. Notable for exponential Euler are perhaps the step size intervals where its waveform inconsistencies are not affected,
$7-8\mu$s and $13-19\mu$s.  These same AP spiking multiplicity variations are illustrated from different vantage points in Figure \ref{fig:methodErrVsts}, 
Figures \ref{fig:p23rsAPclassFE_2OT_RK21RK41}--\ref{fig:p23rsAPclassExpoE} and exponential Euler's waveform plots in Figure 
\ref{fig:corruptExpoEulerAPs}. Simulation videos showing the full simulation train of spiking waveforms for each method over all step sizes are 
included with the supplemental materials. These videos provide a more global perspective of sudden significant changes to AP waveforms as 
step sizes increase. \\ \\
The cause for these waveform anomalies becomes apparent by watching the evolution of membrane voltage, axial and membrane ion channel 
currents as step size increases with $1\mu$s resolution. As step size increases, integration of first order linear membrane channel ODEs with a 
first order numerical method leads to decreasing rise times and increasing widths of pulses with diminished magnitudes. The result in our case 
is for compartment membrane voltage to be pushed just above an unstable equilibrium of the dynamical system at a time that delayed membrane 
ion channels have not yet returned to resting states and results in another AP spike. While we are not hypothesizing the presence of a bifurcation, 
otherwise similar mechanics are described in \cite[ch.3]{Izh}. A video showing the onset of spiking multiplicity variation for one AP cycle integrated 
by RK21 is included in the supplemental materials. \\ \\
The third way we assess the accuracy of a numerical integration method is by observing AP spike phase changes as step size increases. Note that 
each method's AP spike phases are different from each other even when the step size is $1\mu$s.  We also see that every method except HCN 
allows its spike phases to change with increasing step size as shown in Figure \ref{fig:SpikePhaseComp}. \\ \\
Performance of a network model is often evaluated by its power spectrum. So, the fourth way we assessed the accuracy of a numerical integration 
method was by observing three method's power spectrums. Figures \ref{fig:HCNandBEPSD}-\ref{fig:HCNandExpoEPSD} show HCN's power spectral 
densities (PSD) consistent for step sizes $1-99\mu$s. BTCS PSDs were consistent for step sizes up to $80\mu$s when spiking waveform anomalies 
developed. HCN and BTCS PSDs were also consistent with each other. Exponential Euler's power spectral peaks varied between $5$ and $10$ Hz 
as step size grew from $1\mu$s to $50\mu$s.
\subsection{Stability}
\label{subsec:Stability} 
Figure \ref{fig:successfulStepSizes} illustrates step size intervals successfully integrated without oscillation by each numerical method. As 
predicted in Table \ref{tab:schemeStability}, the implicit first-order method BTCS was stable over all step sizes. Exponential Euler was 
also predicted, and found, to be stable for all step sizes, but was only exercised up to a step size of $50\mu$s 
% because its AP waveform statistics were too erratic. \\ \\
due to inconsistent waveform morphology. \\ \\
Von Neumann stability analysis \cite[pg.126]{Str} of first-order FTCS provided the step size upper limit shown in Table \ref{tab:schemeStability} 
and confirmed in Figure \ref{fig:successfulStepSizes}. When FTCS step size increased to $7\mu$s, oscillations appeared on the trailing edges 
of AP spikes as shown in the second pane of Figure \ref{fig:exOsc}, an outcome different than predicted in \cite{Borgers}. \\ \\
RK methods are intended to integrate native ODEs and not spatially discretized PDEs. RK method predictions of stable step size upper limits, 
described by \cite[pg.100]{Butcher}, shown in Figure \ref{fig:RKstepsizeLimits}, do not match the accurate predictions made by Von Neumann 
stability analysis of the quasi-FD RK methods, derived in \ref{subsubsec:RK21stability} and \ref{subsubsec:RK41stability}, listed in Table 
\ref{tab:schemeStability} and also shown in Figure \ref{fig:RKstepsizeLimits}, because their upper step size limit expressions do not consider 
the model's spatial dimension, in particular the model's plane wave's phase angle. \\ \\
Figures \ref{fig:fE2OToscVsTimestep}-\ref{fig:hCNoscVsTimestep} illustrate the RMS of oscillation magnitude, as defined in Section 
\ref{sec:OscDet}, for all integration methods. Except for the last step size prior to expected instability, magnitudes are far below $10^{-6}$mv.
That each method would have any oscillation is not unexpected given the neural model's governing DE changes after every step. \\ \\
While HCN is formally unconditionally stable \cite[pg.62]{CarnevaleHines}, Figure \ref{fig:hCNGrwthFctrVsTimestep} shows its growth factor, 
$g(k,\theta)\Rightarrow-1$ from above as step size increases. Associated oscillation amplitude increases shown in Figures 
\ref{fig:hCNoscVsTimestep},  \ref{fig:ADP}--\ref{fig:hCNOscDecay} are the source of widely varying AP cycle maximums in Figure 
\ref{fig:somaStatsWithStd} which indicate oscillation amplitudes can be as high as $10$ millivolts. Amplitude magnitudes outside the spiking 
interval only rise to a few tenths of a millivolt as shown in Figure \ref{fig:hCNOscDecay}. But given such small oscillation magnitudes it remains 
to be seen what, if any, untoward outcomes are likely. A clearer picture of relative oscillation magnitudes is presented by videos included in the 
supplementary material that show the data presented in Figure \ref{fig:SpikePhaseComp} and \ref{fig:ADP}--\ref{fig:hCNOscDecay} as step size 
changes from $1$ to $99\mu$s. \\ \\
Finally, while our model was deterministic, varying oscillation amplitudes of just one AP cycle for each step size shown in Figures 
\ref{fig:fE2OToscVsTimestep}--\ref{fig:hCNoscVsTimestep} clearly indicate the model's state space is large enough that AP cycles are unlikely 
to ever be identical.
\subsection{Conclusions}
\label{subsec:Conclusions}
The traditional advantage of explicit methods, ease of implementation, did not make up for their short span of stable step sizes 
(Figure \ref{fig:successfulStepSizes}), low order error growth \mbox{(Figure \ref{fig:methodErrVsts})}, and variation in waveform statistics 
\mbox{(Figure \ref{fig:somaStatsWithStd})}. \\ \\
We can only cite two potentially unwanted outcomes of HCN's decaying oscillations. Fluctuation of maximums would matter if, during a 
network simulation, the magnitude of APs reaching a synapse through an axon are proportional to the source neuron's maximum spiking 
magnitude and its attenuation over randomly determined distances. Furthermore, ADP oscillations would matter if their presence misleads 
software intended to detect the ADP's peak just as we experienced (see lower right pane of Figure \ref{fig:p23rsAPclassBE_HCN}). 
Otherwise, HCN's accuracy, stability and immunity to the data shocks experienced by the first order methods, are unmatched in our study. \\ \\
Our results also demonstrate that whenever neural model simulation studies name their simulation platform but not the numerical integration 
method or merely describe their integration step size as ``variable'' without further details or not mention elapsed simulation time prior to 
data collection, repeatability of a study's results is no longer a certainty \cite{MBL}.
\section{Future Work}
A similar examination of the implicit integrators' performance after hyper/depolarizing and synaptic currents have been added would be instructive.
HCN's method of staggering gate state updates should be applied to all other methods to see if their less than unity error order in Figure 
\ref{fig:methodErrVsts} will vanish. The simulations run in this study should also be repeated after substituting symmetric axial resistances for the 
asymmetrically modeled resistances used here. 
\section{Acknowledgements} 
This research did not receive any grants from funding agencies in the public, commercial, or not-for-profit sectors.

\appendix
\section{Integration Methods, Accuracy and Stability} 
\label{Appdx:MethodsAccuracyStability}
After changing the derivative notation in equation (\ref{CableEquationRoot}) we define the cable equation's differential operator
\begin{multline}
\label{CableEquationOperator}
P(V) = V_t + \Big(\alpha+\sum_i\beta_i\Big)V -\Big(\alpha E_L+\sum_i\beta_iE_i\Big) -\gamma\big(a^2V_x\big)_x  = 0,
%~ ... \hfill \\ \\[-1.5ex] %-\gamma\Big(V_{j-1}-2V_j+V_{j+1}\Big)=0.% ~ \cite[p.332,B.19]{Ster}
\end{multline}
where
\[ \alpha=\frac{1}{r_mc_m}, \quad \beta_i=\frac{g_i}{c_m} \quad \text{and} \quad \gamma=\frac{1}{2ar_Lc_m}. \] \\
Defined below as difference operators, $P_{k,h}(\cdot)$, numerical integration schemes used in this study temporally discretize equation 
(\ref{CableEquationOperator}) on a grid of points in the time-space plane. If we let $h$ and $k$ be positive rational numbers, then this grid 
will be the points $(t_n,x_j)=(nk,jh)$ for arbitrary integers $n$ and $j$. For the membrane voltage function, $V$, defined on the grid, 
$V^n_j$ represents the value of $V$ at the grid point $(t_n,x_j)$. \\ \\
After deriving truncation errors below we perform Von Neumann stability analysis as described in \cite{Str} by replacing each occurrence of the 
space-time discretized dependent variable, $V_j^n$, with the corresponding Fourier amplification expression, $g^ne^{ij\theta}$, and then solving 
for the amplification factor, $g(\theta, k, h)$. For this analysis constant terms without the dependent variable do not play a role and are discarded. 
Stability of the two Runge-Kutta methods will be classified using Butcher's methods. \cite[pg.100]{Butcher} \\ \\
In this appendix we use several definitions and identities. Parent and child specific axial conductances will be represented, respectively, as
\begin{equation}
\label{specificAxialResDefs}
c_1=\frac{1}{R_aC_m}=\frac{a^2/a}{2r_Lc_mh^2}~~\text{and}~~c_2=\frac{1}{R'_aC_m}=\frac{a'^2/a}{2r_Lc_mh^2},
\end{equation}
\\
where $a$ is the radius of a given compartment with subscript $j$ and $a'$ the radius of compartment with subscript $j+1$. Clearly the spatial 
derivative in the last term of equation (\ref{CableEquationOperator}) may be expanded as  
\begin{multline}
\label{identitiesForTEAndStab}
\big(a^2V_x\big)_x\equiv\overline{a}^2V_{xx}+\big(\overline{a}^2\big)_xV_x, \\ \\[-1.5ex]
\quad \text{where} \quad\overline{a}^2=\frac{a'^2+a^2}{2}, \quad \text{and}\quad\big(\overline{a}^2\big)_x=\displaystyle\lim_{h\rightarrow 0}\frac{a'^2-a^2}{h}. 
\end{multline}
Finally, we will find it useful to see equation (\ref{CableEquationOperator}) rearranged as 
\begin{equation}
\label{eqn:truncErrIdent}
\Big(\alpha+\sum_i\beta_i\Big)V^n_j - \Big(\alpha E_L+\sum_i\beta_iE_i\Big)-\gamma\big(\overline{a}^2V_x\big)_x  = -V_t.
\end{equation}
\subsection{Forward-Time Central-Space (FTCS)}
\label{subsec:FTCS}
\begin{multline}
\label{FtCsScheme}
P_{k,h}(V)\equiv\frac{V_j^{n+1}-V_j^n}{k} + \Big(\alpha+\sum_i\beta_i\Big)V_j^n -\Big(\alpha E_L+\sum_i\beta_iE_i\Big) ~ ... \hfill \\ \\[-1.5ex]
-\Bigg(\frac{V^n_{j-1}-V^n_j}{R_aC_m}+\frac{V^n_{j+1}-V^n_j}{R'_aC_m}\Bigg)=0.% ~ \cite[p.332,B.19]{Ster}
\end{multline}
\subsubsection{FTCS Accuracy}
\label{subsubsec:FTCSaccuracy}
Multiplying both sides of equation (\ref{FtCsScheme}) by $k$ and expanding discretized dependent variables with their Taylor equivalents we get
\begin{multline}
\label{fTcSAcc1}
\left[V_j^n+kV_t+\frac{k^2}{2}V_{tt}+\frac{k^3}{3!}V_{ttt}+O(k^4)\right]-V_j^n ~...  \hfill \\    
+\Bigg(\Big(\alpha+\sum_i\beta_i\Big)V^n_j - \Big(\alpha E_L+\sum_i\beta_iE_i\Big)~... \Bigg. \hfill \\
-\frac{1}{R_aC_m}\bigg[\bigg(V_j^n-hV_x +\frac{(-h)^2}{2}V_{xx}+\frac{(-h)^3}{3!}V_{xxx}+\frac{(-h)^4}{4!}V_{xxxx}+O(h^5)\bigg) -V_j^n \bigg]~ ... \hfill \\
\Bigg. -\frac{1}{R'_aC_m}\bigg[\bigg(V_j^n+hV_x+\frac{h^2}{2}V_{xx}+\frac{h^3}{3!}V_{xxx}+\frac{h^4}{4!}V_{xxxx}+O(h^5)\bigg)
 -V_j^n\bigg]\Bigg)k=0, \hfill 
 \end{multline}
Simplifying and dividing by $k$ then leads to 
\begin{multline}
\label{fTcSAcc2}
V_t+\frac{k}{2}V_{tt}+\frac{k^2}{6}V_{ttt}+O(k^3)+\Big(\alpha+\sum_i\beta_i\Big)V^n_j - \Big(\alpha E_L+\sum_i\beta_iE_i\Big) ~... \hfill \\
-\frac{a^2/a}{2r_Lc_mh^2}\bigg(-hV_x +\frac{h^2}{2}V_{xx}-\frac{h^3}{6}V_{xxx}+\frac{h^4}{24}V_{xxxx}+O(h^5)\bigg) ~ ... \hfill \\ \\[-2.5ex]
\hfill-\frac{a'^2/a}{2r_Lc_mh^2}\bigg(hV_x+\frac{h^2}{2}V_{xx}+\frac{h^3}{6}V_{xxx}+\frac{h^4}{24}V_{xxxx}+O(h^5)\bigg) =0. 
\end{multline}
Multiplication over addition, regrouping and eliminating third order terms gives us
\begin{multline}
\label{fTcSAcc3}
V_t+\Big(\alpha+\sum_i\beta_i\Big)V^n_j - \Big(\alpha E_L+\sum_i\beta_iE_i\Big)+\frac{V_{tt}}{2}k+\frac{V_{ttt}}{6}k^2 ~ ... \hfill \\
\hfill -\gamma\bigg[\big(a'^2+a^2\big)\bigg(\frac{V_{xx}}{2}+\frac{V_{xxxx}}{24}h^2\bigg)
+\big(a'^2-a^2\big)\bigg(\frac{V_x}{h} + \frac{V_{xxx}}{6}h\bigg)\bigg]=0.
\end{multline}
After rearranging terms we have
\begin{multline}
\label{fTcSAcc4}
V_t+\Big(\alpha+\sum_i\beta_i\Big)V^n_j - \Big(\alpha E_L+\sum_i\beta_iE_i\Big)+\frac{V_{tt}}{2}k+\frac{V_{ttt}}{6}k^2 ~ ... \hfill \\
-\gamma\bigg[\bigg(\frac{a'^2+a^2}{2}V_{xx}+\frac{a'^2-a^2}{h}V_x \bigg)
+\bigg(\frac{a'^2+a^2}{2}V_{xxxx}+2\frac{a'^2-a^2}{h}V_{xxx}\bigg)\frac{h^2}{12}\bigg]=0.
\end{multline}
Replacing terms above with identities defined in (\ref{identitiesForTEAndStab}) gives us  
\begin{multline}
\label{fTcSAcc6}
V_t+\Big(\alpha+\sum_i\beta_i\Big)V^n_j - \Big(\alpha E_L+\sum_i\beta_iE_i\Big)-\gamma\big(\overline{a}^2V_x\big)_x ~ ... \hfill \\
+\frac{V_{tt}}{2}k+\frac{V_{ttt}}{6}k^2
-\gamma\bigg(\overline{a}^2V_{xxxx}+2\big(\overline{a}^2\big)_xV_{xxx}\bigg)\frac{h^2}{12}.
\end{multline}
Therefore, the truncation error is 
\begin{multline}
\label{fTcSAcc7}
P_{k,h}(V)-P(V) =  \frac{V_{tt}}{2}k+\frac{V_{ttt}}{6}k^2-\gamma\bigg(\overline{a}^2V_{xxxx}+2\big(\overline{a}^2\big)_xV_{xxx}\bigg)\frac{h^2}{12}
\end{multline}
\subsubsection{FTCS Stability}
\label{subsubsec:FTCSstability}
After rearranging the FTCS scheme terms in \ref{FtCsScheme} and dropping the constant term we have
\begin{multline}
\label{fEstab1}
V_j^{n+1}= V_j^n + \Big(\alpha E_L+\sum_i\beta_iE_i\Big)k ~ ...  \hfill \\ 
\hfill -\Big(\alpha+\sum_i\beta_i\Big)kV_j^n+\bigg(\frac{V_{j-1}^n-V_j^n}{R_aC_m}+\frac{V_{j+1}^n-V_j^n}{R'_aC_m}\bigg)k. 
\end{multline}
\\ \\
Then after setting $V_j^n=g^ne^{ij\theta}$ we have
\begin{multline}
\label{fEstab2}
g^{n+1}e^{ij\theta} = g^ne^{ij\theta}\bigg(1-\Big(\alpha+\sum_i\beta_i\Big)k 
+\Big(c_1e^{-i\theta}-\big(c_1+c_2\big)+c_2e^{i\theta}\Big)k \bigg). \hfill \\ \\[-1.5ex]
\Rightarrow g(\cdot)=1-\Big(\alpha+\sum_i\beta_i\Big)k +\Big(c_1e^{-i\theta}-\big(c_1+c_2\big)+c_2e^{i\theta}\Big)k \hfill 
\end{multline}
\begin{multline}
\label{fEstab3}
\Rightarrow g(\cdot)=1-\Big(\alpha+\sum_i\beta_i\Big)k + \Big(c_1(\cos\theta-i\sin\theta)-(c_1+c_2)+c_2(\cos\theta+i\sin\theta)\Big)k \hfill \\ \\[-1.5ex]
\Rightarrow g(\cdot)=1-\Big(\alpha+\sum_i\beta_i\Big)k + \Big(\big(\cos\theta-1\big)(c_1+c_2)-i\sin\theta(c_1-c_2)\Big)k \hfill \\ \\[-1.5ex]
\Rightarrow g(\cdot)=1-\Big(\alpha+\sum_i\beta_i\Big)k  + \Big(-2\sin^2{\left(\frac{\theta}{2}\right)}(c_1+c_2)-i\sin\theta(c_1-c_2)\Big)k \hfill
\end{multline}
Finally, 
\begin{multline}
\label{fEstab4}
g(\cdot)~=~\Big(1-\big(K+2L\big)k\Big)-iMk,\hfill 
\end{multline}
\\ \\[-1.5ex]
where $K=\Big(\alpha+\sum_i\beta_i\Big),~ L=\big(c_1+c_2\big)\sin^2{\left(\frac{\theta}{2}\right)}$, and $M=\big(c_1-c_2\big)\sin{\theta}$. \\ \\ \\[-1.5ex]
Because the model's endpoints float, basis function spatial derivatives must be zero at the boundaries which implies the only basis functions 
are cosines and we may disregard the imaginary sine term. Then, for the scheme to be stable, $|g(\theta)|^2<1.$
\begin{multline}
\label{fEstab5}
\Rightarrow |g(\theta, k, h)|^2= \Big(1-\big(K+2L\big)k\Big)^2<1\hfill \\ \\[-1.5ex]
\Rightarrow 1-2(K+2L)k+(K+2L)^2k^2<1 \\ \\[-1.5ex]
\Rightarrow (K+2L)^2k^2<2(K+2L)k \\ \\[-1.5ex]
\Rightarrow k<\frac{2}{(K+2L)}.
\end{multline}
As shown in Figure \ref{fig:fEstepSizeLimit}, oscillations are expected when step size grows larger than $6\mu$sec.
\subsection{Exponential Euler}
\label{subsec:EE}
If we rewrite equation (\ref{FtCsScheme}) as
\begin{multline}
\label{eqn:SpatiallyDscrtzdCable}
\Big(V_t\Big)_j = \Bigg(\alpha E_L+\sum_i\beta_iE_i + \bigg(\frac{V_{j-1}}{R_aC_m}+\frac{V_{j+1}}{R'_aC_m}\bigg)\Bigg)~ ... \hfill \\ \\[-2.0ex]
-\Bigg(\Big(\alpha+\sum_i\beta_i\Big) + \bigg(\frac{1}{R_aC_m}+\frac{1}{R'_aC_m}\bigg)\Bigg)V_j,
\end{multline} 
and assume the values $V_{j-1}$ and $V_{j+1}$ are known, we now see the cable equation cast as the first order, linear, ODE,
\begin{equation}
\label{firstOlinODE}
\frac{dV}{dt}=A-BV,
\end{equation}
whose analytical solution, provided $A$ and $B$ are constants, is
\begin{equation}
\label{firstOlinODEsol}
V(t-t_0)=\frac{A}{B}+\left(V(t_0)-\frac{A}{B}\right)\exp\Big(-B(t-t_0)\Big).
\end{equation}
The exponential Euler method's approximation operator is derived by rearranging terms and temporally discretizing the membrane 
voltage in equation (\ref{firstOlinODEsol}), after which we see the difference operator defined as   
\setlength\abovedisplayskip{5pt}
\begin{multline}
\label{expoEAcc0}
P_{k,h}(V)=V_j^{n+1} - \frac{A}{B}\left(1 - \exp(-Bk)\right) - V_j^n~\exp(-Bk) = 0. \hfill 
\end{multline}
\subsubsection{Exponential Euler Accuracy}
\label{subsubsec:EEaccuracy}
After replacing A with its definition from equation (\ref{eqn:SpatiallyDscrtzdCable}) and replacing discretized as well as exponential 
terms with their Taylor expansion equivalents, we may rewrite equation (\ref{expoEAcc0}) as
\begin{multline}
\label{expoEAcc1}
\left[V_j^n+kV_t+\frac{k^2}{2}V_{tt}+\frac{k^3}{3!}V_{ttt}+O(k^4)\right] - \frac{1}{B}\Bigg[\alpha E_L + \sum_i \beta_iE_i  \Bigg. ...  \hfill \\
+ \frac{1}{R_aC_m}\bigg(V_j^n-hV_x+\frac{(-h)^2}{2}V_{xx}+\frac{(-h)^3}{3!}V_{xxx}+\frac{(-h)^4}{4!}V_{xxxx}+O(h^5)\bigg) ~ ... \hfill \\
\Bigg. + \frac{1}{R'_aC_m}\bigg(V_j^n+hV_x+\frac{h^2}{2}V_{xx}+\frac{h^3}{3!}V_{xxx}+\frac{h^4}{4!}V_{xxxx}+O(h^5)\bigg)\Bigg] ~ ... \hfill \\ \\[-1.5ex]
*\bigg(Bk-\frac{(Bk)^2}{2} + \frac{(Bk)^3}{3!}-O(Bk)^4\bigg)-V_j^n\bigg(1-Bk+\frac{(Bk)^2}{2}-\frac{(Bk)^3}{3!}+O(Bk)^4\bigg) =0. 
\end{multline}
Simplifying and dividing by $B$ leads to 
\begin{multline}
\label{expoEAcc2}
kV_t+\frac{k^2}{2}V_{tt}+\frac{k^3}{6}V_{ttt}+O(k^4) - \Bigg[\alpha E_L + \sum_i \beta E_i + ~ ...  \Bigg. \hfill \\
\frac{a^2/a}{2r_Lc_mh^2}\bigg(V_j^n-hV_x+\frac{h^2}{2}V_{xx}-\frac{h^3}{6}V_{xxx}+\frac{h^4}{24}V_{xxxx}\bigg) ~ ... \hfill \\
\Bigg.+\frac{a^2/a}{2r_Lc_mh^2}\bigg(V_j^n+hV_x+\frac{h^2}{2}V_{xx}+\frac{h^3}{6}V_{xxx}+\frac{h^4}{24}V_{xxxx}\bigg)
+O(h^3)\Bigg] ~ ... \hfill \\ \\[-1.5ex]
* \bigg(k-\frac{Bk^2}{2}+\frac{B^2k^3}{6} -kO(Bk)^3\bigg)+\bigg(Bk-\frac{(Bk)^2}{2}+\frac{(Bk)^3}{6}-O(Bk)^4\bigg)V_j^n=0.
\end{multline}
Dividing by $k$ and discarding fourth order terms gives us
\begin{multline}
\label{expoEAcc3}
V_t+\frac{k}{2}V_{tt}+\frac{k^2}{6}V_{ttt}- \bigg(1-\frac{Bk}{2}+\frac{(Bk)^2}{6} -O(Bk)^3\bigg)*\Bigg[\alpha E_L + \sum_i \beta E_i \Bigg. ~ ... \hfill \\
\Bigg. + \frac{\gamma}{h^2}\Bigg(\big(a'^2+a^2\big)\bigg(V_j^n+\frac{h^2}{2}V_{xx}+\frac{h^4}{24}V_{xxxx}\bigg)
+\big(a'^2-a^2\big)\bigg(hV_x+\frac{h^3}{6}V_{xxx}\bigg)\Bigg)\Bigg]~... \\
+\bigg(1-\frac{Bk}{2}+\frac{(Bk)^2}{6}-O(Bk)^3\bigg)BV_j^n=O(h^3). 
\end{multline}
Discarding third order terms and factoring out the common term leaves us with
\begin{multline}
\label{expoEAcc4}
V_t+\frac{k}{2}V_{tt}+\frac{k^2}{6}V_{ttt}- \bigg(1-\frac{Bk}{2}+\frac{(Bk)^2}{6} \bigg)*\left[\alpha E_L + \sum_i \beta_i E_i \right. ~ ... \hfill \\
+ \gamma\Bigg(\frac{a'^2+a^2}{2}\bigg(\frac{2V_j^n}{h^2}+V_{xx}+\frac{h^2}{12}V_{xxxx}\bigg)
+\frac{a'^2-a^2}{h}\bigg(V_x+\frac{h^2}{6}V_{xxx}\bigg)\Bigg) ~ ... \\
\Bigg. -BV_j^n\Bigg]=0.
\end{multline}
Replacing terms above with identities defined in equation (\ref{identitiesForTEAndStab}) and the last occurrence of B with its definition from 
equation (\ref{eqn:SpatiallyDscrtzdCable}) gives us 
\begin{multline}
\label{expoEAcc5}
V_t+\frac{k}{2}V_{tt}+\frac{k^2}{6}V_{ttt}- \bigg(1-\frac{Bk}{2}+\frac{(Bk)^2}{6} \bigg)*\left[ \alpha E_L + \sum_i \beta E_i \right. ~ ... \hfill \\
 + \gamma\Bigg(\Big(\overline{a}^2V_{xx}+\big(\overline{a}^2\big)_xV_x\Big)+\Big(\overline{a}^2V_{xxxx}
+2\big(\overline{a}^2\big)_xV_{xxx}\Big)\frac{h^2}{12}+\overline{a}^2\frac{2V_j^n}{h^2}\Bigg) ~ ... \\
\left.-\Bigg(\alpha+\sum_i \beta_i + \frac{2\gamma}{h^2}\overline{a}^2\Bigg)V_j^n\right]=0.
\end{multline}
Regrouping terms gives us \\ 
\begin{multline}
\label{expoEAcc6}
V_t+\frac{k}{2}V_{tt}+\frac{k^2}{6}V_{ttt} ~ ... \hfill \\
+\bigg(1-\frac{Bk}{2}+\frac{(Bk)^2}{6} \bigg)*\Bigg[\alpha(V_j^n-E_L) + \sum_i \beta_i(V_j^n-E_i)-\gamma\big(\overline{a}^2V_x\big)_x \Bigg. ~ ... \\
\Bigg. -\gamma\Big(\overline{a}^2V_{xxxx}+2\big(\overline{a}^2\big)_xV_{xxx}\Big)\frac{h^2}{12}\Bigg]=0. 
\end{multline}
Rearranging terms and using the identity in equation (\ref{eqn:truncErrIdent}) to reduce terms inside the square brackets we get 
\begin{multline}
\label{expoEAcc7}
V_t+\alpha(V_j^n-E_L) + \sum_i \beta_i(V_j^n-E_i)-\gamma\big(\overline{a}^2V_x\big)_x ~... \hfill \\ \\[-2.0ex]
+\frac{V_{tt}}{2}k+\frac{V_{ttt}}{6}k^2 +\left(\frac{Bk}{2} - \frac{(Bk)^2}{6}\right)V_t ~... \hfill \\ \\[-2.0ex]
-\gamma\Big(\overline{a}^2V_{xxxx}+2\big(\overline{a}^2\big)_xV_{xxx}\Big)\frac{h^2}{12}\bigg(1-\frac{Bk}{2}+\frac{(Bk)^2}{6}\bigg) =0. 
\end{multline}
Truncation error magnitude for the exponential Euler integration method is therefore
\begin{multline}
\label{expoEAcc8}
P_{h,k}(V)-P(V)=\frac{V_{tt}}{2}k+\frac{V_{ttt}}{6}k^2+\left(\frac{Bk}{2} - \frac{(Bk)^2}{6}\right)V_t~... \hfill \\ 
-\gamma\Big(\overline{a}^2V_{xxxx}+2\big(\overline{a}^2\big)_xV_{xxx}\Big)\frac{h^2}{12}\bigg(1-\frac{Bk}{2}+\frac{(Bk)^2}{6}\bigg).
\end{multline}
\subsubsection{Exponential Euler Stability}
\label{subsubsec:EEstability}
Determining whether the exponential Euler method is stable will be easier if we first divide the term A by B, as they were identified in 
equations (\ref{eqn:SpatiallyDscrtzdCable}) and (\ref{firstOlinODE}). This gives us 
\begin{multline}
\label{eEstab0}
\frac{A}{B}=\frac{R_aR'_aE_L+R_mR_aR'_a\sum_i G_iE_i +R_mR'_aV^n_{j-1} + R_mR_aV^n_{j+1}}
{R_aR'_a+R_mR_aR'_a\sum_i G_i + R_mR'_a + R_mR_a}.\hfill 
\end{multline}
\\ \\
After substituting equation (\ref{eEstab0}) into the discretized exponential Euler approximation operator, equation (\ref{expoEAcc0}), and 
setting $V_j^n=g^ne^{ij\theta}$, we have 
\begin{multline}
\label{eEstab1}
g^{n+1}e^{ij\theta} = g^ne^{ij\theta}e^{-Bk} + ... \\
\left(\frac{R_aR'_aE_L+R_mR_aR'_a\sum_i G_iE_i +R_mR'_ag^ne^{i(j-1)\theta} + R_mR_ag^ne^{i(j+1)\theta}}
{R_aR'_a+R_mR_aR'_a\sum_i G_i + R_mR'_a + R_mR_a}\right)~ ... \\
*\left(1-e^{-Bk}\right).
\end{multline}
\\
Disregarding the constants in the numerator and simplifying gives us 
\begin{multline}
\label{eEstab2}
g(\cdot) = e^{-Bk} + \left(\frac{R_m(R'_a+R_a)\cos\theta  +iR_m(R_a-R'_a)\sin\theta}
{R_aR'_a+R_mR_aR'_a\sum_i G_i + R_m(R'_a + R_a)}\right)(1-e^{-Bk}).
\end{multline}
\\ \\
The magnitude of the exponential Euler method's growth factor is
\begin{multline}
\label{eEstab3}
|g(\cdot)|^2 = Re(g(\cdot))^2+Im(g(\cdot))^2 = ... \hfill \\
\left[\bigg(e^{-Bk}\Big(R_aR'_a+R_mR_aR'_a\sum_i G_i + R_m(R'_a + R_a)\Big)~ ... \bigg. \right. \hfill \\
\Bigg. \bigg. +R_m(R'_a+R_a)\cos\theta(1-e^{-Bk})\bigg)^2  + \bigg(R_m(R_a-R'_a)\sin\theta(1-e^{-Bk})\bigg)^2\Bigg]  ~ ... \\
\bigg/ \bigg(R_aR'_a+R_mR_aR'_a\sum_i G_i + R_m(R'_a + R_a)\bigg)^2. 
\end{multline}
\\
After squaring the numerator's first term we have 
\begin{multline}
\label{eEstab4}
\text{numerator}(|g(\cdot)|^2) = e^{-2Bk}\Big(R_aR'_a+R_mR_aR'_a\sum_i G_i + R_m(R'_a + R_a)\Big)^2~ ...\hfill \\
+2e^{-Bk}\Big(R_aR'_a+R_mR_aR'_a\sum_i G_i + R_m(R'_a + R_a)\Big)\Big(R_m(R'_a+R_a)\cos\theta(1-e^{-Bk})\Big)~ ... \\
+\Big(R_m(R'_a+R_a)\cos\theta(1-e^{-Bk})\Big)^2+\bigg(R_m(R'_a-R_a)\sin\theta(1-e^{-Bk})\bigg)^2. \hfill
\end{multline}
Replacing the exponential terms with the first three terms of their Taylor series and rejoining numerator with denominator leads to 
\begin{multline}
\label{eEstab5}
|g(\cdot)|^2 = (1-2Bk+O(Bk)^2)~ ... \hfill \\ \\[-1.5ex]
+\frac{2(1-Bk+O(Bk)^2)(Bk+O(Bk)^2)R_m(R'_a+R_a)\cos\theta}{\bigg(R_aR'_a+R_mR_aR'_a\sum_i G_i + R_m(R'_a + R_a)\bigg)}~ ... \\ \\[-1.5ex]
\hfill +\frac{R^2_m\Big(Bk+O(Bk)^2\Big)^2\Big(R'^2_a+2R_aR'_a\cos2\theta+R^2_a \Big)}
{\bigg(R_aR'_a+R_mR_aR'_a\sum_i G_i + R_m(R'_a + R_a)\bigg)^2} 
\end{multline}
\begin{multline}
\label{eEstab6}
\Rightarrow |g(\cdot)|^2 \le (1-2Bk+O(Bk)^2)~ ... \hfill \\ \\[-1.5ex]
+\left(\frac{R_m(R'_a+R_a)\cos\theta}{R_aR'_a+R_mR_aR'_a\sum_i G_i + R_m(R'_a + R_a)}\right)\left(2Bk+O(Bk)^2\right)~ ...\\ \\[-1.5ex]
\hfill +\left(\frac{R_m(R'_a+R_a)}{R_aR'_a+R_mR_aR'_a\sum_i G_i + R_m(R'_a + R_a)}\right)^2\left((Bk)^2+O(Bk)^3\right) 
\end{multline}
\begin{multline}
\label{eEstab7}
\Rightarrow |g(\cdot)|^2< 1-2Bk+2Bk+O(Bk)^2+(Bk)^2-O(Bk)^3 \hfill 
\end{multline}
\\
Clearly, $\quad |g(\theta, h, k)|^2 \le 1+O(Bk)^2$, meaning that the growth factor's magnitude is less than or equal to unity. Therefore the 
exponential Euler method is unconditionally stable.
\subsection{Backward-Time Central-Space (BTCS)} 
\label{subsec:BTCS}
\begin{multline}
\label{bTcSscheme}
P_{k,h}(V)\equiv\frac{V_j^{n+1}-V_j^n}{k} + \Big(\alpha+\sum_i\beta_i\Big)V_j^{n+1} -\Big(\alpha E_L+\sum_i\beta_iE_i\Big) ~ ... \hfill \\
\hfill -\Bigg(\frac{V_{j-1}^{n+1}-V_j^{n+1}}{R_aC_m}+\frac{V_{j+1}^{n+1}-V_j^{n+1}}{R'_aC_m}\Bigg)=0. 
\end{multline}
\subsubsection{BTCS Accuracy}
\label{subsubsec:BTCSaccuracy}
As shown in equation (\ref{bTcSscheme}), we define the BTCS scheme's difference operator as \cite[p.~332, B.24]{Ster}
\begin{multline}
\label{bTcSConsist0}
P_{k,h}(V) = \frac{V_j^{n+1}-V_j^n}{k} + \Big(\alpha+\sum_i\beta_i\Big)V_j^{n+1} -\Big(\alpha E_L+\sum_i\beta_iE_i\Big) ~ ... \hfill \\
\hfill -\Bigg(\frac{V_{j-1}^{n+1}-V_j^{n+1}}{R_aC_m}+\frac{V_{j+1}^{n+1}-V_j^{n+1}}{R'_aC_m}\Bigg)=0.
\end{multline}
Multiplying both sides of equation (\ref{bTcSConsist0}) by $k$ and expanding discretized dependent variables with their Taylor equivalents we get %\\
\begin{multline}
\label{bTcSAcc0}
\left[V_j^n+kV_t+\frac{k^2}{2}V_{tt}+\frac{k^3}{3!}V_{ttt}+O(k^4)\right]-V_j^n~ ...  \hfill \\ 
+\left(\Big(\alpha+\sum_i\beta_i\Big)\Big(V^n_j+kV_t+\frac{k^2}{2}V_{tt}+\frac{k^3}{3!}V_{ttt}+O(k^4)\Big) - \Big(\alpha E_L+\sum_i\beta_iE_i\Big)~ ... \right.\\
-\frac{1}{R_aC_m}\Bigg[\bigg(V_j^n-hV_x+kV_t +\frac{1}{2}\Big((-h)^2V_{xx}+2(-h)kV_{xt}+k^2V_{tt}\Big)\bigg. \Bigg. ~ ... \hfill \\
\bigg. +\frac{1}{3!}\Big((-h)^3V_{xxx}+3(-h)^2kV_{xxt}+3(-h)k^2V_{xtt}+k^3V_{ttt}\Big) ~ ... \hfill \\
\bigg.+\frac{1}{4!}\Big((-h)^4V_{xxxx}+4(-h)^3kV_{xxxt}+6(-h)^2k^2V_{xxtt}+4(-h)k^3V_{xttt}+k^4V_{tttt}\Big)\bigg) ~ ... \\
\Bigg. -\bigg(V_j^n+kV_t+\frac{k^2}{2}V_{tt}+\frac{k^3}{3!}V_{ttt}+\frac{k^4}{4!}V_{tttt}\bigg)+O(h^rk^s)^{r+s=5}\Bigg]~ ... \\
-\frac{1}{R'_aC_m}\Bigg[\bigg(V_j^n+hV_x+kV_t +\frac{1}{2}\Big(h^2V_{xx}+2hkV_{xt}+k^2V_{tt}\Big) \bigg. \Bigg. ~ ... \hfill \\
\bigg. +\frac{1}{3!}\Big(h^3V_{xxx}+3h^2kV_{xxt}+3hk^2V_{xtt}+k^3V_{ttt}\Big) ~ ... \hfill \\
\bigg.+\frac{1}{4!}\Big(h^4V_{xxxx}+4h^3kV_{xxxt}+6h^2k^2V_{xxtt}+4hk^3V_{xttt}+k^4V_{tttt}\Big) \bigg)~ ... \hfill \\
\Bigg. \bigg. -\bigg(V_j^n+kV_t+\frac{k^2}{2}V_{tt}+\frac{k^3}{3!}V_{ttt}+\frac{k^4}{4!}V_{tttt}\bigg)+O(h^rk^s)^{r+s=5}\Bigg]\Bigg)k=0.  
\end{multline}
Simplifying and dividing by k leads to \\
\begin{multline}
\label{bTcSAcc1}
V_t+\frac{k}{2}V_{tt}+\frac{k^2}{6}V_{ttt}+O(k^3)~... \hfill \\
+\Big(\alpha+\sum_i\beta_i\Big)\bigg(V^n_j+kV_t+\frac{k^2}{2}V_{tt}+O(k^3)\bigg)-\Big(\alpha E_L+\sum_i\beta_iE_i\Big) ~... \\
-\frac{a^2/a}{2r_Lc_mh^2}\bigg(-hV_x+\frac{h^2}{2}V_{xx}-hkV_{xt}-\frac{h^3}{6}V_{xxx}+\frac{h^2k}{2}V_{xxt}-\frac{hk^2}{2}V_{xtt}\bigg. ~ ... \hfill \\
\hfill \bigg. +\frac{h^4}{24}V_{xxxx}-\frac{h^3k}{6}V_{xxxt}+\frac{h^2k^2}{4}V_{xxtt}-\frac{hk^3}{6}V_{xttt}+O(h^rk^s)^{r+s=5}\bigg) ~... \\
-\frac{a'^2/a}{2r_Lc_mh^2}\bigg(hV_x+\frac{h^2}{2}V_{xx}+hkV_{xt}+\frac{h^3}{6}V_{xxx}+\frac{h^2k}{2}V_{xxt}+\frac{hk^2}{2}V_{xtt} \bigg. ~ ... \hfill \\
\hfill \bigg. +\frac{h^4}{24}V_{xxxx}+\frac{h^3k}{6}V_{xxxt}+\frac{h^2k^2}{4}V_{xxtt}+\frac{hk^3}{6}V_{xttt}+O(h^rk^s)^{r+s=5}\bigg)=0. 
\end{multline}
\\
Multiplication over addition, regrouping and rearranging gives us
\begin{multline}
\label{bTcSAcc2}
V_t+\Big(\alpha+\sum_i\beta_i\Big)V^n_j-\Big(\alpha E_L+\sum_i\beta_iE_i\Big) ~ ... \hfill \\
+\Big(\alpha+\sum_i\beta_i\Big)\bigg(kV_t+\frac{k^2}{2}V_{tt}+O(k^3)\bigg)+\frac{k}{2}V_{tt}+\frac{k^2}{6}V_{ttt}+O(k^3) ~ ... \hfill \\
-\gamma\bigg[\big(a'^2+a^2\big)
\bigg(\frac{V_{xx}}{2}+\frac{k}{2}V_{xxt}+\frac{V_{xxxx}}{24}h^2+\frac{k^2}{4}V_{xxtt}+O(h^rk^s)^{r+s=3}\bigg)\bigg. ~ ... \hfill \\
\bigg. +\big(a'^2-a^2\big)\bigg(\frac{V_x}{h} +\frac{k}{h}V_{xt}+\frac{V_{xxx}}{6}h+\frac{k^2}{2h}V_{xtt}+\frac{hk}{6}V_{xxxt}+\frac{k^3}{6h}V_{xttt}+O(h^rk^s)^{r+s=3}\bigg)\bigg]=0.
\end{multline}
Rearranging further leads to
\begin{multline}
\label{bTcSAcc3}
V_t+\Big(\alpha+\sum_i\beta_i\Big)V^n_j -\Big(\alpha E_L+\sum_i\beta_iE_i\Big) ~... \hfill \\
+\left(\frac{V_{tt}}{2}+\Big(\alpha+\sum_i\beta_i\Big)V_t\right)k+\bigg(\frac{V_{ttt}}{3}+\Big(\alpha+\sum_i\beta_i\Big)V_{tt}\bigg)\frac{k^2}{2}~ ... \\ \\[-2.0ex]
-\gamma\Bigg[\Big(\overline{a}^2V_{xx}+\big(\overline{a}^2\big)_xV_x\Big)+\Big(\overline{a}^2V_{xxt}+\big(\overline{a}^2\big)_xV_{xt}\Big)k
+\Big(\overline{a}^2V_{xxtt}+\big(\overline{a}^2\big)_xV_{xtt}\Big)\frac{k^2}{2} ~ ... \\
\overline{a}^2\frac{V_{xxxx}}{12}h^2+\big(\overline{a}^2\big)_x\bigg(\frac{V_{xxx}}{6}h^2+\frac{V_{xxxt}}{6}h^2k+\frac{V_{xttt}}{6}k^3 \bigg)\Bigg]=0.
\end{multline}
Rearranging again, replacing terms above with identities defined in (\ref{identitiesForTEAndStab}) and eliminating third order terms gives us
\begin{multline}
\label{bTcSAcc4}
V_t+\Big(\alpha+\sum_i\beta_i\Big)V^n_j -\Big(\alpha E_L+\sum_i\beta_iE_i\Big)-\gamma\big(\overline{a}^2V_x\big)_x ~... \hfill \\
+\left(\frac{V_{tt}}{2}+\Big(\alpha+\sum_i\beta_i\Big)V_t-\big(\overline{a}^2V_x\big)_{xt}\right)k
+\bigg(\frac{V_{ttt}}{3}+\Big(\alpha+\sum_i\beta_i\Big)V_{tt}-\big(\overline{a}^2V_x\big)_{xt}\bigg)\frac{k^2}{2}~ ... \\ \\[-2.0ex]
-\gamma\Bigg[\overline{a}^2\frac{V_{xxxx}}{12}h^2+\big(\overline{a}^2\big)_x\bigg(\frac{V_{xxx}}{6}h^2+\frac{V_{xxxt}}{6}h^2k
+\frac{V_{xttt}}{6}k^3 \bigg)\Bigg]=0.
\end{multline}
Differentiating equation (\ref{eqn:truncErrIdent}) with respect to time gives us an identity with which we reduce the coefficient of $k$. 
Differentiating equation (\ref{eqn:truncErrIdent}) twice with respect to time yields yet another identity with which we reduce the coefficient 
of $k^2$, leaving us with \\
\begin{multline}
\label{bTcSAcc5}
V_t+\Big(\alpha+\sum_i\beta_i\Big)V^n_j -\Big(\alpha E_L+\sum_i\beta_iE_i\Big) -\gamma\big(\overline{a}^2V_x\big)_x
-\frac{V_{tt}}{2}k-\frac{V_{ttt}}{3}k^2 ~... \\
-\gamma\Bigg[\overline{a}^2\frac{V_{xxxx}}{12}h^2+\big(\overline{a}^2\big)_x\bigg(\frac{V_{xxx}}{6}h^2+\frac{V_{xxxt}}{6}h^2k
+\frac{V_{xttt}}{6}k^3 \bigg)\Bigg]=0,
\end{multline}
\\
The truncation error for the BTCS scheme is therefore,  
\begin{multline}
\label{bTcSAcc6}
P_{k,h}(V)-P(V)= ~ ... \hfill \\ \\[-1.5ex]
-\frac{V_{tt}}{2}k-\frac{V_{ttt}}{3}k^2 -\gamma\Bigg[\overline{a}^2\frac{V_{xxxx}}{12}h^2
+\big(\overline{a}^2\big)_x\bigg(\frac{V_{xxx}}{6}h^2+\frac{V_{xxxt}}{6}h^2k+\frac{V_{xttt}}{6}k^3 \bigg)\Bigg].
\end{multline}
\subsubsection{BTCS Stability}
\label{subsubsec:BTCSstability}
By rearranging the BTCS scheme terms in equation (\ref{bTcSConsist0}) we have
\begin{multline}
\label{bEstab0}
V_j^{n+1}+\Big(\alpha+\sum_i\beta_i\Big)kV_j^{n+1} -\Big(\alpha E_L+\sum_i\beta_iE_i\Big)k ~ ...  \hfill \\ 
\hfill - \bigg(\frac{V_{j-1}^{n+1}-V_j^{n+1}}{R_aC_m}+\frac{V_{j+1}^{n+1}-V_j^{n+1}}{R'_aC_m}\bigg)k=V_j^n.
\end{multline}
\\ \\
After dropping constants and replacing dependent variables with their phased growth factor equivalents we have
\begin{multline}
\label{bEstab1}
g^{n+1}e^{ij\theta}\bigg(1+\Big(\alpha+\sum_i\beta_i\Big)k\bigg)~... \hfill \\
\hfill -\left(c_1g^{n+1}e^{i(j-1)\theta}-(c_1+c_2)g^{n+1}e^{ij\theta}+c_2g^{n+1}e^{i(j+1)\theta}\right)k = g^ne^{ij\theta}
\end{multline}
\begin{multline}
\label{bEstab2}
\Rightarrow g\bigg(1+\Big(\alpha+\sum_i\beta_i\Big)k-\Big(c_1e^{-i\theta}-(c_1+c_2)+c_2e^{i\theta}\Big)k\bigg)=1 \hfill \\ \\[-1.5ex]
\Rightarrow g\bigg(1+\Big(\alpha+\sum_i\beta_i\Big)k~... \bigg. \hfill \\ 
\bigg. -\Big(c_1(\cos\theta-i\sin\theta)-(c_1+c_2)+c_2(\cos\theta+i\sin\theta)\Big)k\bigg)=1 \\ \\[-2.0ex]
\Rightarrow g\bigg(1+\Big(\alpha+\sum_i\beta_i\Big)k-\Big((\cos{\theta}-1)(c_1+c_2)+i(c_2-c_1)\sin\theta\Big)k\bigg)=1 
\end{multline}
\\
\begin{multline}
\label{bEstab3}
\Rightarrow g\bigg(1+\Big(\alpha+\sum_i\beta_i\Big)k
+\Big(2(c_1+c_2)\sin^2{\textstyle\left(\frac{\theta}{2}\right)}+i(c_1-c_2)\sin\theta\Big)k\bigg)=1 \hfill 
\end{multline}
\begin{multline}
\label{bEstab4}
\Rightarrow g(\theta, k)=\frac{1}{1+\Big(K+2L+iM\Big)k}, \qquad \parbox{1cm}{where} \hfill \\ \\[-1.5ex]
K=\Big(\alpha+\sum_i\beta_i\Big),~ L=\big(c_1+c_2\big)\sin^2{\textstyle\left(\frac{\theta}{2}\right)}, ~ \text{and } M=(c_1-c_2)\sin{\theta}. \hfill 
\end{multline}
\\
Because the model's endpoints float, basis function spatial derivatives must be zero at the boundaries which implies the only basis functions 
are cosines and we may disregard the imaginary sine term. So, we are left with
\begin{equation}
\label{bEstab5}
g(\theta, k)=\frac{1}{1+\big(K+2L\big)k}, 
\end{equation}
which is always positive.  The magnitude of the growth factor will then always be less than unity and the BTCS scheme 
is unconditionally stable for all $\theta$ and $k$.
\subsection{$2^{nd}$-Order Taylor}
\label{subsec:2OT}
The Taylor expansion of membrane voltage one time step after some reference time $t_0$ is
\begin{equation}
\label{2OTAcc0}
V(t_0+k)=V(t_0)+kV_t+\frac{k^2}{2}V_{tt}+O(k^3).
\end{equation}
Section \ref{subsubsec:EEaccuracy} described how the cable equation could be cast as the first order linear ODE, 
\begin{equation}
\label{2OTAcc1}
V_t=A-BV, \quad \text{where} 
\end{equation}
\begin{multline}
\label{2OTAcc2}
A=\alpha E_L +\sum_i \beta_i E_i +\bigg(\frac{V^n_{j-1}}{R_aC_m} +\frac{V^n_{j+1}}{R'_aC_m}\bigg), \\
\text{ and } \quad B=\alpha +\sum_i \beta_i +\bigg(\frac{1}{R_aC_m} +\frac{1}{R'_aC_m}\bigg).
\end{multline}
For every Taylor method we must choose how to approximate the higher derivatives. We are therefore in search of an explicit second 
order, or better, approximation of the second derivative. One option is a three point centered finite difference. While this may first appear 
as an implicit method, it is a simple matter to solve for $V_j^{n+1}$ explicitly. However the new approximation is the midpoint method 
and known to be unstable. A second option is a lagging three point stencil, using $V_j^n,V_j^{n-1}$ and $V_j^{n-2}$, unfortunately
this approximation is only first order. A third option would be to approximate the spatial derivative with a lagging four point stencil. But 
while the method is second order it is also inconsistent. \\ \\ 
Regarding equation (\ref{2OTAcc1}) we know that $A$ and $B$ are functions of membrane voltage which changes with respect to 
time. Therefore, the second derivative of the membrane voltage, or first derivative of equation (\ref{2OTAcc1}) is
\begin{equation}
\label{eqn:Vtt}
V_{tt}=A_t - \Big(BV_t+VB_t\Big),
\end{equation}
and may be rendered
\begin{multline}
\label{2OTAcc3}
\displaystyle\lim_{k\rightarrow 0} \bigg(\frac{A^n-A^{n-2}}{2k}\bigg) ~ ... \hfill \\ \\[-1.5ex]
\hfill -\bigg(\frac{B^n+B^{n-2}}{2}\bigg)\bigg(\frac{V^n-V^{n-2}}{2k}\bigg)-\bigg(\frac{V^n+V^{n-2}}{2}\bigg)\bigg(\frac{B^n-B^{n-2}}{2k}\bigg)= ~ ... \\ \\[-1.5ex]
\displaystyle\lim_{k\rightarrow 0}\frac{\Big(A^n-B^nV^n\Big)-\Big(A^{n-2}-B^{n-2}V^{n-2}\Big)}{2k}=\frac{(V^n)_t-(V^{n-2})_t}{2k}.
\end{multline}
\\
Taylor expanding the last expression gives us   
\begin{multline}
\label{2OTAcc4}
%\displaystyle V_{tt}=\lim_{k \rightarrow 0} \frac{(V^n)_t-(V^{n-2})_t}{2k} ~ ... \hfill \\ \\[-1.5ex]
\bigg(\Big(\big(V^{n-1}\big)_t+k\big(V^{n-1}\big)_{tt}+\frac{k^2}{2}\big(V^{n-1}\big)_{ttt}+O\big(k^3\big)\Big)~ ... \\ \\[-1.5ex]
-\Big(\big(V^{n-1}\big)_t-k\big(V^{n-1}\big)_{tt}+\frac{k^2}{2}\big(V^{n-1}\big)_{ttt}+O\big(k^3\big)\Big)\bigg) \bigg/2k  ~ ... \\ \\[-1.5ex]
= \lim_{k \rightarrow 0}\Big(\big(V^{n-1}\big)_{tt}+O\big(k^2\big)\Big).
\end{multline}
Our Taylor method is derived by substituting equations (\ref{2OTAcc1}) and (\ref{2OTAcc3}) for the temporal derivatives in equation (\ref{2OTAcc0}), 
discretizing the dependent variable and rearranging gives us the $2^{nd}$ order Taylor method's approximation operator,
\begin{multline}
\label{2OTAcc5}
P_{k,h}(V)\equiv ~ ... \\ \\[-1.0ex]
\frac{V_j^{n+1}-V_j^n}{k} - (A-BV_j^n) -\frac{k}{2}\bigg[\frac{\big(A^n-B^nV^n\big)-\big(A^{n-2}-B^{n-2}V^{n-2}\big)}{2k}\bigg]=0. \hfill
\end{multline}
\subsubsection{$2^{nd}$-Order Taylor Accuracy}
\label{subsec:2OTaccuracy}
After multiplying both sides of equation (\ref{2OTAcc5}) by $k$ we have 
\begin{multline}
\label{2OTAcc6}
V_j^{n+1}-V_j^n -k(A^n-B^nV_j^n)-\frac{k}{4}\bigg( \Big(A^n-B^nV^n\Big)-\Big(A^{n-2}-B^{n-2}V^{n-2}\Big)\bigg)=0. 
\end{multline}
Replacing $A$ and $B$ with their definitions and discretized terms with their Taylor expansions gives us
\begin{multline}
\label{2OTAcc8}
\bigg(V_j^n+kV_t+\frac{k^2}{2}V_{tt}+\frac{k^3}{3!}V_{ttt}+O(k^4)\bigg) -V_j^n 
-\frac{5k}{4}\Bigg[\alpha E_L+\sum_i\beta_iE_i-\Big(\alpha+\sum_i\beta_i\Big)V^n_j ~... \Bigg. \hfill  \\
+\frac{1}{R_aC_m}\Bigg(\bigg[V_j^n-hV_x +\frac{(-h)^2}{2}V_{xx}+\frac{(-h)^3}{3!}V_{xxx}+\frac{(-h)^4}{4!}V_{xxxx}+O(h^5)\bigg] -V_j^n \Bigg)~ ... \hfill \\
\Bigg. +\frac{1}{R'_aC_m}\Bigg(\bigg[V_j^n+hV_x+\frac{h^2}{2}V_{xx}+\frac{h^3}{3!}V_{xxx}+\frac{h^4}{4!}V_{xxxx}+O(h^5)\bigg]-V_j^n\Bigg)\Bigg]~... \hfill \\
+\frac{k}{4}\Bigg[\alpha E_L+\sum_i\beta_iE_i-\Big(\alpha+\sum_i\beta_i\Big)V^{n-2}_j ~... \Bigg. \hfill  \\
+\frac{1}{R_aC_m}\Bigg(\bigg[V_j^{n-2}-hV_x +\frac{(-h)^2}{2}V_{xx}+\frac{(-h)^3}{3!}V_{xxx}+\frac{(-h)^4}{4!}V_{xxxx}+O(h^5)\bigg] -V_j^{n-2} \Bigg)~ ... \hfill \\
\Bigg. +\frac{1}{R'_aC_m}\Bigg(\bigg[V_j^{n-2}+hV_x+\frac{h^2}{2}V_{xx}+\frac{h^3}{3!}V_{xxx}+\frac{h^4}{4!}V_{xxxx}+O(h^5)\bigg]-V_j^{n-2}\Bigg)\Bigg]=0.
\end{multline}
After simplifying and dividing both sides by $k$ we have
\begin{multline}
\label{2OTAcc9}  
V_t+\frac{k}{2}V_{tt}+\frac{k^2}{6}V_{ttt}
+\frac{5}{4}\Bigg(\Big(\alpha+\sum_i\beta_i\Big)V^n_j-\Big(\alpha E_L+\sum_i\beta_iE_i\Big) \Bigg. ~... \hfill  \\
-\frac{a^2/a}{2r_Lc_mh^2}\bigg(-hV_x +\frac{h^2}{2}V_{xx}-\frac{h^3}{6}V_{xxx}+\frac{h^4}{24}V_{xxxx} \bigg)~ ... \\
\Bigg. -\frac{a'^2/a}{2r_Lc_mh^2}\bigg(hV_x+\frac{h^2}{2}V_{xx}+\frac{h^3}{3!}V_{xxx}+\frac{h^4}{4!}V_{xxxx}\bigg) \Bigg)~...\\
-\frac{1}{4}\Bigg(\Big(\alpha+\sum_i\beta_i\Big)V^{n-2}_j-\Big(\alpha E_L+\sum_i\beta_iE_i\Big) \Bigg. ~... \hfill  \\
-\frac{a^2/a}{2r_Lc_mh^2}\bigg(-h\big(V_j^{n-2}\big)_x +\frac{h^2}{2}\big(V_j^{n-2}\big)_{xx}-\frac{h^3}{6}\big(V_j^{n-2}\big)_{xxx}
+\frac{h^4}{24}\big(V_j^{n-2}\big)_{xxxx} \bigg)~ ... \\
\Bigg. -\frac{a'^2/a}{2r_Lc_mh^2}\bigg(h\big(V_j^{n-2}\big)_x+\frac{h^2}{2}\big(V_j^{n-2}\big)_{xx}+\frac{h^3}{3!}\big(V_j^{n-2}\big)_{xxx}
+\frac{h^4}{4!}\big(V_j^{n-2}\big)_{xxxx}\bigg) \Bigg) ~ ... \\
=O(k^3,h^3).   
\end{multline}
Rearranging and eliminating third order terms gives us
\begin{multline}
\label{2OTAcc10}
V_t+\frac{k}{2}V_{tt}+\frac{V_{ttt}}{6}k^2+\frac{5}{4}\Bigg(\Big(\alpha+\sum_i\beta_i\Big)V^n_j - \Big(\alpha E_L+\sum_i\beta_iE_i\Big) \Bigg. ~... \hfill \\
\Bigg. -\gamma\bigg[\Big(a'^2+a^2\Big)\bigg(\frac{V_{xx}}{2}+\frac{V_{xxxx}}{24}h^2 \bigg) 
+\Big(a'^2-a^2\Big)\bigg(\frac{V_x}{h}+\frac{V_{xxx}}{6}h\bigg)\bigg] \Bigg)~ ... \\
-\frac{1}{4}\Bigg(\Big(\alpha+\sum_i\beta_i\Big)V_j^{n-2} - \Big(\alpha E_L+\sum_i\beta_iE_i\Big) \Bigg. ~... \hfill \\
\Bigg. -\gamma\bigg[\Big(a'^2+a^2\Big)\bigg(\frac{\big(V_j^{n-2}\big)_{xx}}{2}+\frac{\big(V_j^{n-2}\big)_{xxxx}}{24}h^2 \bigg) ~ ... \hfill \\
\hfill +\Big(a'^2-a^2\Big)\bigg(\frac{\big(V_j^{n-2}\big)_x}{h}+\frac{\big(V_j^{n-2}\big)_{xxx}}{6}h\bigg)\bigg] \Bigg)=0. 
\end{multline}
Replacing terms above with identities defined in equation (\ref{identitiesForTEAndStab}) leads to
\begin{multline}
\label{2OTAcc11}
V_t+\frac{k}{2}V_{tt}+\frac{V_{ttt}}{6}k^2+\frac{5}{4}\Bigg(\Big(\alpha+\sum_i\beta_i\Big)V^n_j - \Big(\alpha E_L+\sum_i\beta_iE_i\Big) \Bigg. ~... \hfill \\
\Bigg. -\gamma\bigg[\Big(\overline{a}^2V_{xx}+\big(\overline{a}^2\big)_xV_x\Big)+\Big(\overline{a}^2V_{xxxx} 
+2\big(\overline{a}^2\big)_xV_{xxx}\Big)\frac{h^2}{12}\bigg]\Bigg) ~... \\
-\frac{1}{4}\Bigg(\Big(\alpha+\sum_i\beta_i\Big)V_j^{n-2} - \Big(\alpha E_L+\sum_i\beta_iE_i\Big) \Bigg. 
\Bigg. -\gamma\bigg[\Big(\overline{a}^2\big(V_j^{n-2}\big)_{xx}+\big(\overline{a}^2\big)_x\big(V_j^{n-2}\big)_x\Big) ~ ... \hfill \\
+\Big(\overline{a}^2\big(V_j^{n-2}\big)_{xxxx} +2\big(\overline{a}^2\big)_x\big(V_j^{n-2}\big)_{xxx}\Big)\frac{h^2}{12}\bigg] \Bigg)=0.
\end{multline}
Moving $\frac{1}{4}$ of the term whose coefficient is $\frac{5}{4}$ to the term below whose coefficient is $-\frac{1}{4}$ and applying the identity in 
equation (\ref{eqn:truncErrIdent}) gives us
\begin{multline}
\label{2OTAcc12}
V_t+\Big(\alpha+\sum_i\beta_i\Big)V^n_j - \Big(\alpha E_L+\sum_i\beta_iE_i\Big) -\gamma\big(\overline{a}^2V_x\big)_x~... \hfill \\
\hfill +\frac{V_{tt}}{2}k +\frac{V_{ttt}}{6}k^2-\gamma\Big(\overline{a}^2V_{xxxx} +2\big(\overline{a}^2\big)_xV_{xxx}\Big)\frac{h^2}{12} ~ ... \\
-\frac{1}{4}\Bigg(\big(V_j^n\big)_t-\big(V_j^{n-2}\big)_t 
+\gamma\bigg[\overline{a}^2\big(V_j^{n-2}\big)_{xxxx} +2\big(\overline{a}^2\big)_x\big(V_j^{n-2}\big)_{xxx}\bigg]\frac{h^2}{12}\Bigg)=0.
\end{multline}
Rearranging terms leaves us with 
\begin{multline}
\label{2OTAcc13}
V_t+\Big(\alpha+\sum_i\beta_i\Big)V^n_j - \Big(\alpha E_L+\sum_i\beta_iE_i\Big) -\gamma\big(\overline{a}^2V_x\big)_x~... \hfill \\
\hfill -\frac{1}{4}\Big(\big(V_j^n\big)_t-\big(V_j^{n-2}\big)_t\Big) +\frac{V_{tt}}{2}k +\frac{V_{ttt}}{6}k^2 ~ ... \\
-\gamma \Bigg( \Big[\overline{a}^2V_{xxxx} +2\big(\overline{a}^2\big)_xV_{xxx}\Big]
+\frac{1}{4}\bigg[\overline{a}^2\big(V_j^{n-2}\big)_{xxxx} +2\big(\overline{a}^2\big)_x\big(V_j^{n-2}\big)_{xxx}\bigg]\Bigg)\frac{h^2}{12}=0.
\end{multline}
Taylor expanding $\big(V_j^{n-2}\big)_{xxxx}$ and $\big(V_j^{n-2}\big)_{xxx}$ about $\big(V_j^n\big)_{xxxx}$ and $\big(V_j^n\big)_{xxx}$, respectively, 
turns the last expression into a $O(h^2k)$ term, which we disregard. Unfortunately, the second order approximation from equation (\ref{2OTAcc3}), %\\
\mbox{$\Big(\big(V_j^n\big)_t-\big(V_j^{n-2}\big)_t\Big)= 2k\big(V_j^n\big)_{tt}$}, is only true in the limit, as step size $k\rightarrow0$.  Truncation 
error for the \mbox{$2^{nd}$-Order Taylor} method is therefore,
\begin{multline}
\label{2OTAcc14}
P_{k,h}(V)-P(V) =~... \\ \\[-1.5ex]
\frac{1}{2}\Bigg(V_{tt}-\bigg(\frac{\big(V_j^n\big)_t-\big(V_j^{n-2}\big)_t}{2}\bigg)\neq0\Bigg)k + \frac{V_{ttt}}{6}k^2 
-\gamma\Big(\overline{a}^2V_{xxxx} +2\big(\overline{a}^2\big)_xV_{xxx}\Big)\frac{h^2}{12}, \hfill
\end{multline}
\\
which is only first order in time, second in space.
\subsubsection{$2^{nd}$-Order Taylor Stability}
\label{subsubsec:2OTstability}
\begin{figure}[h!]
\centering
\begin{tabular}{l}
\includegraphics[height=2.0in,width=4.5in]{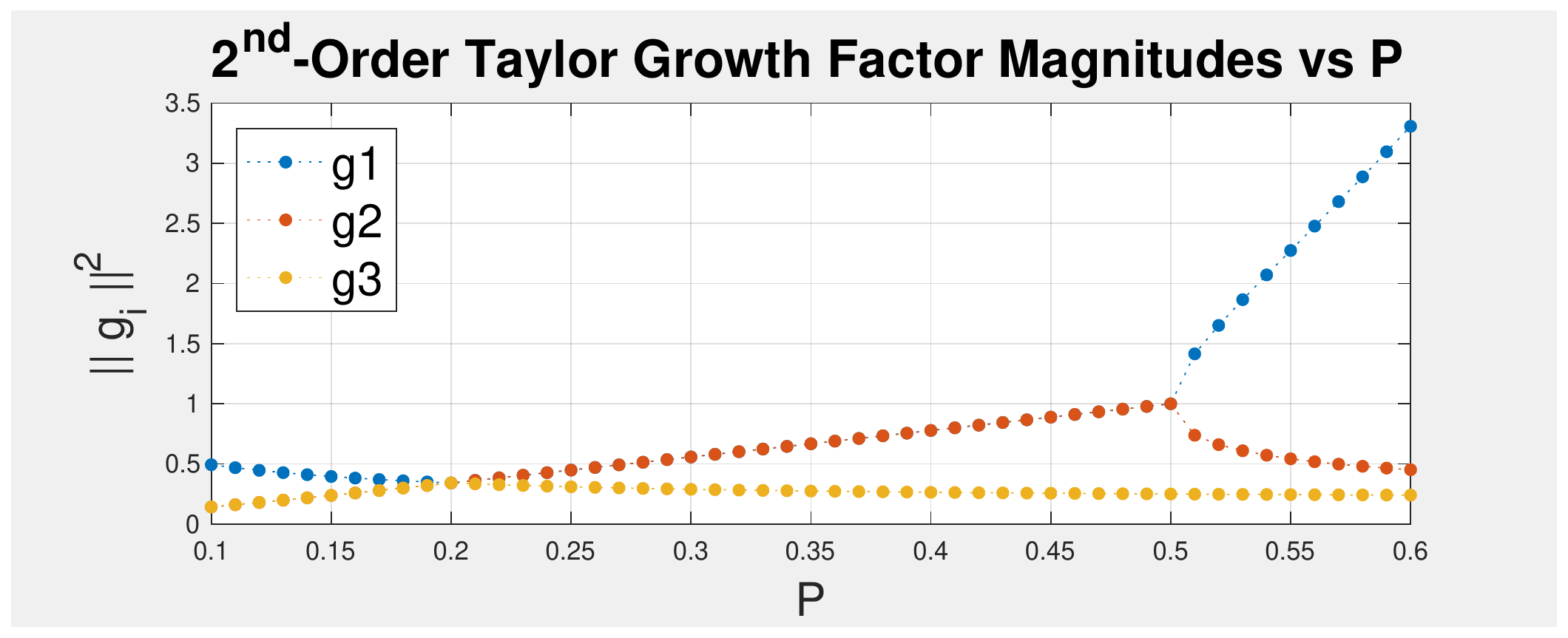}
\end{tabular}
\caption{}
\label{fig:grwthFctrMagsVsP}
\end{figure}

Starting from equation (\ref{2OTAcc6}) we have
\begin{multline}
\label{2OTStab0}
V_j^{n+1}-V_j^n -k\Bigg(\bigg[\alpha E_L+\sum_i\beta_iE_i +\bigg(\frac{V^n_{j-1}}{R_aC_m} +\frac{V^n_{j+1}}{R'_aC_m}\bigg)\bigg]\Bigg. ~...\hfill \\
\Bigg. - \bigg[\alpha +\sum_i \beta_i +\bigg(\frac{1}{R_aC_m} +\frac{1}{R'_aC_m}\bigg) \bigg] V_j^n \Bigg)~ ... \\
-\frac{k}{4}\Bigg(\Big(A^n-B^nV^n\Big)-\Big(A^{n-2}-B^{n-2}V^{n-2}\Big)\Bigg)=0. 
\end{multline}
After dropping constants, replacing dependent variables with their phased growth factor equivalents and combining like terms we have
\begin{multline}
\label{2OTstab1}
g^{n+1}e^{ij\theta}= g^ne^{ij\theta}~...\hfill \\
+\frac{5k}{4}\Bigg(-\Big(\alpha+\sum_i\beta_i\Big)g^ne^{ij\theta}+\bigg[\frac{g^ne^{i(j-1)\theta}-g^ne^{ij\theta}}{R_aC_m}
+\frac{g^ne^{i(j+1)\theta}-g^ne^{ij\theta}}{R'_aC_m}\bigg]\Bigg)~... \\
-\frac{k}{4}\Bigg(-\Big(\alpha+\sum_i\beta_i\Big)g^{n-2}e^{ij\theta} \Bigg. ~ ... \hfill \\
\Bigg. +\bigg[\frac{g^{n-2}e^{i(j-1)\theta}-g^{n-2}e^{ij\theta}}{R_aC_m}+\frac{g^{n-2}e^{i(j+1)\theta}-g^{n-2}e^{ij\theta}}{R'_aC_m}\bigg]\Bigg)
\end{multline}
\begin{multline}
\label{2OTstab2}
g=1+\frac{5k}{4}\Bigg(-\Big(\alpha+\sum_i\beta_i\Big)+\bigg(\frac{e^{-i\theta}-1}{R_aC_m}+\frac{e^{i\theta}-1}{R'_aC_m}\bigg)\Bigg) ~ ... \hfill \\
-\frac{k}{4}\Bigg(-\Big(\alpha+\sum_i\beta_i\Big)+\bigg(\frac{e^{-i\theta}-1}{R_aC_m}+\frac{e^{i\theta}-1}{R'_aC_m}\bigg)\Bigg)g^{-2} 
\end{multline}
\begin{multline}
\label{2OTstab3}
g^3-\Bigg[1+\frac{5k}{4}\Bigg(-\Big(\alpha+\sum_i\beta_i\Big)+ \bigg(\frac{\cos{\theta}-i\sin{\theta}-1}{R_aC_m} 
+\frac{\cos{\theta}+i\sin{\theta}-1}{R'_aC_m}\bigg)\Bigg)\Bigg]g^2 ~ ... \\
+\frac{k}{4}\Bigg(-\Big(\alpha+\sum_i\beta_i\Big)+ \bigg(\frac{\cos{\theta}-i\sin{\theta}-1}{R_aC_m} 
+\frac{\cos{\theta}+i\sin{\theta}-1}{R'_aC_m}\bigg)\Bigg)=0. %\\ \\
\end{multline}
\begin{multline}
\label{2OTstab4}
g^3-\Bigg(1+\frac{5k}{4}\bigg(-\Big(\alpha+\sum_i\beta_i\Big)+ \Big(\big(c_1+c_2\big)\big(\cos{\theta}-1\big)+i\big(c_2-c_1\big)\sin{\theta}\Big)\bigg)\Bigg)g^2 ~ ... \\
+\frac{k}{4}\bigg(-\Big(\alpha+\sum_i\beta_i\Big)+ \Big(\big(c_1+c_2\big)\big(\cos{\theta}-1\big)+i\big(c_2-c_1\big)\sin{\theta}\Big)\bigg)=0. \\ \\
g^3-\bigg(1-\frac{5k}{4}\Big(K+2L-iM \Big)\bigg)g^2 - \frac{k}{4}\Big(K+2L-iM\Big) = 0, \qquad \parbox{1cm}{where} \hfill \\ \\[-1.5ex]
K=\Big(\alpha+\sum_i\beta_i\Big),~ L=\big(c_1+c_2\big)\sin^2{\textstyle\left(\frac{\theta}{2}\right)}, ~ \text{and } M=(c_2-c_1)\sin{\theta}. \hfill 
\end{multline}
Because the model's endpoints float, basis function spatial derivatives must be zero at the boundaries which implies the only basis functions 
are cosines and we may disregard the imaginary sine term. We now seek values of $g$ such that
\begin{equation}
\label{2OTstab5}
g^3+(5P-1)g^2-P=0, \quad \text{where} ~P=\frac{k}{4}(K+2L).
\end{equation}
To locate the roots of this cubic polynomial we use the MATLAB function $roots([a_3~a_2~a_1~a_0])$. Figure \ref{fig:grwthFctrMagsVsP}
illustrates that for the magnitude of the growth factor to be less than or equal to one, P can be no larger than $0.5$. Putting this value back into 
equation (\ref{2OTstab5}) we get the same step size limitation relationship as for FTCS in equation (\ref{fEstab5}).
\subsection{Hines-Crank-Nicolson (HCN)}
\label{subsec:HCN}
The Hines adaptation of the Crank-Nicolson scheme achieves second order accuracy in time and space by leading with a half-step implicit 
backward-Euler-like scheme and following with a half-step explicit forward-Euler-like scheme as shown in \cite{Hines}. The first implicit half-step 
difference operator for the cable equation was defined in \cite{Hines} as
\begin{multline}
\label{hCNAcc0}
P_{k,h}(V) = \frac{V_j^{n+\frac{1}{2}}-V_j^n}{k/2} + \Big(\alpha+\sum_i\beta_i\Big)V_j^{n+\frac{1}{2}} -\Big(\alpha E_L+\sum_i\beta_iE_i\Big) ~ ... \hfill \\
\hfill -\Bigg(\frac{V_{j-1}^{n+\frac{1}{2}}-V_j^{n+\frac{1}{2}}}{R_aC_m}+\frac{V_{j+1}^{n+\frac{1}{2}}-V_j^{n+\frac{1}{2}}}{R'_aC_m}\Bigg)=0,  
\end{multline}
after multiplying by $\frac{k}{2}$ and rearranging terms we have
\begin{multline}
\label{hCNAcc1}
V_j^{n+\frac{1}{2}} - V_j^n = -\Bigg(\Big(\alpha+\sum_i\beta_i\Big)V_j^{n+\frac{1}{2}}-\Big(\alpha E_L+\sum_i\beta_iE_i\Big) \Bigg.~... \hfill \\
\hfill \Bigg. -\bigg(\frac{V_{j-1}^{n+\frac{1}{2}}-V_j^{n+\frac{1}{2}}}{R_aC_m}
+\frac{V_{j+1}^{n+\frac{1}{2}}-V_j^{n+\frac{1}{2}}}{R'_aC_m}\bigg)\Bigg)\frac{k}{2}. 
\end{multline}
\\ \\
The subsequent HCN explicit half-step was also defined in \cite{Hines} as 
\begin{equation}
\label{eqn:hinesExplHalfStep}
V^{n+1}_j=2\left(V^{n+\frac{1}{2}}_j-V^n_j\right)+V^n_j. 
\end{equation}
\subsubsection{HCN Accuracy}
\label{subsubsec:HCNaccuracy}
Combining equations (\ref{hCNAcc0}) and (\ref{eqn:hinesExplHalfStep}) leads to 
\begin{multline}
\label{hCNAcc2}
V^{n+1}_j=-2\Bigg(\Big(\alpha+\sum_i\beta_i\Big)V_j^{n+\frac{1}{2}}-\Big(\alpha E_L+\sum_i\beta_iE_i\Big)\Bigg. ~... \hfill \\
\hfill \Bigg. -\bigg(\frac{V_{j-1}^{n+\frac{1}{2}}-V_j^{n+\frac{1}{2}}}{R_aC_m}
+\frac{V_{j+1}^{n+\frac{1}{2}}-V_j^{n+\frac{1}{2}}}{R'_aC_m}\bigg)\Bigg)\frac{k}{2} +V^n_j.  
\end{multline}
\\ \\
Replacing the first two discretized dependent variables with their Taylor expansion gives us
\begin{multline}
\label{hCNAcc3}
\Big[V_j^n+kV_t+\frac{k^2}{2}V_{tt}+\frac{k^3}{3!}V_{ttt}+O(k^4)\Big] - V_j^n ~ ... \hfill \\
+\Bigg( \Big(\alpha+\sum_i\beta_i\Big)\left(V_j^n+\frac{k}{2}V_t+\frac{(k/2)^2}{2}V_{tt}+\frac{(k/2)^3}{3!}V_{ttt}+O(k^4)\right) ... \Bigg. \\ 
\hfill \Bigg. -\Big(\alpha E_L+\sum_i\beta_iE_i\Big) - \bigg(\frac{V_{j-1}^{n+\frac{1}{2}}-V_j^{n+\frac{1}{2}}}{R_aC_m}
+\frac{V_{j+1}^{n+\frac{1}{2}}-V_j^{n+\frac{1}{2}}}{R_aC_m}\bigg) \Bigg)k =0. 
\end{multline}
Simplifying, dividing by $k$ and expanding the remaining terms leads to 
\begin{multline}
\label{hCNAcc4}
V_t+\frac{k}{2}V_{tt}+\frac{k^2}{6}V_{ttt}+O(k^3)~ ...  \hfill \\ 
+\Big(\alpha+\sum_i\beta_i\Big)\bigg(V_j^n+\frac{k}{2}V_t+\frac{k^2}{8}V_{tt}+\frac{k^3}{48}V_{ttt}+O(k^4)\bigg)- \Big(\alpha E_L+\sum_i\beta_iE_i\Big) ~... \hfill \\
-\frac{1}{R_aC_m}\bigg[\bigg(V_j^n-hV_x+\frac{k}{2}V_t +\frac{1}{2}\Big((-h)^2V_{xx}+2(-h)\frac{k}{2}V_{xt}+(k/2)^2V_{tt}\Big) \bigg.\bigg. ~... \hfill\\
\hfill +\frac{1}{3!}\Big((-h)^3V_{xxx}+3(-h)^2(k/2)V_{xxt}+3(-h)(k/2)^2V_{xtt}+(k/2)^3V_{ttt}\Big) ~ ... \\
\bigg. +\frac{1}{4!}\Big((-h)^4V_{xxxx}+4(-h)^3\Big(\frac{k}{2}\Big)V_{xxxt}+6(-h)^2\Big(\frac{k}{2}\Big)^2V_{xxtt}
+4(-h)\Big(\frac{k}{2}\Big)^3V_{xttt}+\Big(\frac{k}{2}\Big)^4V_{tttt}\Big)\bigg) ~ ... \\
\bigg. - \bigg(V_j^n+\frac{k}{2}V_t+\frac{(k/2)^2}{2}V_{tt}+\frac{(k/2)^3}{3!}V_{ttt}+\frac{(k/2)^4}{4!}V_{tttt}\bigg)+O(h^rk^s)^{r+s=5}\bigg]~ ... \\
- \frac{1}{R_aC_m}\bigg[\bigg(V_j^n+hV_x+\frac{k}{2}V_t +\frac{1}{2}\Big(h^2V_{xx}+2h\frac{k}{2}V_{xt}+(k/2)^2V_{tt}\Big)~... \bigg.\bigg. \hfill \\
\qquad +\frac{1}{3!}\Big(h^3V_{xxx}+3h^2(k/2)V_{xxt}+3h(k/2)^2V_{xtt}+(k/2)^3V_{ttt}\Big) ~ ... \hfill \\
\bigg.+\frac{1}{4!}\Big(h^4V_{xxxx}+4h^3\Big(\frac{k}{2}\Big)V_{xxxt}+6h^2\Big(\frac{k}{2}\Big)^2V_{xxtt}
+4h\Big(\frac{k}{2}\Big)^3V_{xttt}+\Big(\frac{k}{2}\Big)^4V_{tttt}\Big)\bigg) ~ ... \hfill \\
- \bigg(V_j^n+\frac{k}{2}V_t+\frac{(k/2)^2}{2}V_{tt}+\frac{(k/2)^3}{3!}V_{ttt}+\frac{(k/2)^4}{4!}V_{tttt}\bigg)+O(h^rk^s)^{r+s=5}\bigg]=0. \hfill
\end{multline}
\\ \\
Simplifying leads to 
\begin{multline}
\label{hCNAcc5}
V_t+\frac{k}{2}V_{tt}+\frac{k^2}{6}V_{ttt} 
+\Big(\alpha+\sum_i\beta_i\Big)\bigg(V_j^n+\frac{k}{2}V_t+\frac{k^2}{8}V_{tt}\bigg) - \Big(\alpha E_L+\sum_i\beta_iE_i\Big) ~ ... \\
-\frac{a^2/a}{2r_Lc_mh^2}\bigg(-hV_x +\frac{h^2}{2}V_{xx}-h\frac{k}{2}V_{xt}-\frac{h^3}{6}V_{xxx}+\frac{h^2k}{4}V_{xxt}-\frac{hk^2}{8}V_{xtt}\bigg.~ ... \hfill \\
\hfill \bigg.+\frac{h^4}{24}V_{xxxx}-\frac{h^3k}{12}V_{xxxt}+\frac{h^2k^2}{16}V_{xxtt}-\frac{hk^3}{48}V_{xttt}+O(h^rk^s)^{r+s=5}\bigg)~ ... \\
-\frac{a'^2/a}{2r_Lc_mh^2}\bigg(hV_x +\frac{h^2}{2}V_{xx}+h\frac{k}{2}V_{xt}+\frac{h^3}{6}V_{xxx}+\frac{h^2k}{4}V_{xxt}+\frac{hk^2}{8}V_{xtt} ~ ... \hfill \\
\hfill \bigg.+\frac{h^4}{24}V_{xxxx}+\frac{h^3k}{12}V_{xxxt}+\frac{h^2k^2}{16}V_{xxtt}+\frac{hk^3}{48}V_{xttt}+O(h^rk^s)^{r+s=5}\bigg)=0. 
\end{multline}
\\ \\
Dividing by $h^2$ and combining terms related to $a'$ and $a$ gives us 
\begin{multline}
\label{hCNAcc6}
V_t+\frac{k}{2}V_{tt}+\frac{k^2}{6}V_{ttt} +\Big(\alpha+\sum_i\beta_i\Big)\bigg(V_j^n+\frac{k}{2}V_t+\frac{k^2}{8}V_{tt}\bigg) 
- \Big(\alpha E_L+\sum_i\beta_iE_i\Big) ... \hfill \\
- \gamma\Bigg(\Big(a'^2-a^2\Big)\bigg(\frac{V_x}{h}+\frac{k}{2h}V_{xt} + \frac{h}{6}V_{xxx}+\frac{k^2}{8h}V_{xtt}
+\frac{hk}{12}V_{xxxt}+\frac{k^3}{48h}V_{xttt}\bigg)~ ...\hfill \\
+\Big(a'^2+a^2\Big)\bigg(\frac{V_{xx}}{2}+\frac{k}{4}V_{xxt}+\frac{h^2}{24}V_{xxxx}+\frac{k^2}{16}V_{xxtt} \bigg) ~ ... \\
\hfill \Bigg.+O(h^rk^s)^{r+s=5}\Bigg)=0. 
\end{multline}
\\ \\ 
Factoring an $h$ from the sum of terms into the denominator of $\big(a'^2-a^2\big)$, likewise a $2$ from the sum of terms into the denominator 
of $\big(a'^2+a^2\big)$, applying identities defined in equation (\ref{identitiesForTEAndStab}), eliminating third order terms and simplifying leads to
\begin{multline}
\label{hCNAcc7}
V_t+\Big(\alpha+\sum_i\beta_i\Big)V^n_j -\Big(\alpha E_L+\sum_i\beta_iE_i\Big) ~... \hfill \\
+\Bigg(V_{tt}+\Big(\alpha+\sum_i\beta_i\Big)V_t\Bigg)\frac{k}{2}
+\Bigg(\frac{4V_{ttt}}{3}+\Big(\alpha+\sum_i\beta_i\Big)V_{tt}\Bigg)\frac{k^2}{8}~ ... \\
- \gamma\Bigg[\Big(\overline{a}^2V_{xx}+\big(\overline{a}^2\big)_xV_x\Big)
+\Big(\overline{a}^2V_{xxt}+\big(\overline{a}^2\big)_xV_{xt}\Big)\frac{k}{2}
+\Big(\overline{a}^2V_{xxtt}+\big(\overline{a}^2\big)_xV_{xtt}\Big)\frac{k^2}{8} ~ ... \\
+\overline{a}^2\frac{V_{xxxx}}{12}h^2+\big(\overline{a}^2\big)_x\bigg(\frac{V_{xxx}}{6}h^2+\frac{V_{xxxt}}{6}h^2k\bigg)\Bigg]=0. 
\end{multline}
\\
Rearranging again gives us
\begin{multline}
\label{hCNAcc8}
V_t+\Big(\alpha+\sum_i\beta_i\Big)V^n_j -\Big(\alpha E_L+\sum_i\beta_iE_i\Big)-\gamma\big(\overline{a}^2V_x\big)_x ~... \hfill \\
+\Bigg(V_{tt}+\Big(\alpha+\sum_i\beta_i\Big)V_t-\gamma\big(\overline{a}^2V_x\big)_{xt}\Bigg)\frac{k}{2}
+\Bigg(\frac{4V_{ttt}}{3}+\Big(\alpha+\sum_i\beta_i\Big)V_{tt}-\gamma\big(\overline{a}^2V_x\big)_{xtt}\Bigg)\frac{k^2}{8}~ ... \\
- \gamma\bigg(\overline{a}^2V_{xxxx}+2\big(\overline{a}^2\big)_xV_{xxx}\bigg)\frac{h^2}{12}=0. 
\end{multline}
Differentiating equation (\ref{eqn:truncErrIdent}) with respect to time yields an identity with which we zero-out the coefficient of $\frac{k}{2}$. 
Differentiating equation (\ref{eqn:truncErrIdent}) twice with respect to time yields a similar identity with which we reduce the coefficient of 
$\frac{k^2}{8}$, after which we have \\
\begin{multline}
\label{hCNAcc9}
V_t+\Big(\alpha+\sum_i\beta_i\Big)V^n_j -\Big(\alpha E_L+\sum_i\beta_iE_i\Big) -\gamma\big(\overline{a}^2V_x\big)_x ~ ... \\
+\frac{V_{ttt}}{24}k^2- \gamma\bigg(\overline{a}^2V_{xxxx}+2\big(\overline{a}^2\big)_xV_{xxx}\bigg)\frac{h^2}{12}=0. 
\end{multline}
Hence, truncation error for the combined implicit and explicit half-steps in the HCN numerical integration method is 
\begin{multline}
\label{hCNAcc10}
P_{k,h}(V)-P(V) = \frac{V_{ttt}}{24}k^2- \gamma\bigg(\overline{a}^2V_{xxxx}+2\big(\overline{a}^2\big)_xV_{xxx}\bigg)\frac{h^2}{12}. \hfill
\end{multline}
\subsubsection{HCN Stability}
\label{subsubsec:HCNstability}
Each step of the HCN numerical integration method for the cable equation was described in \cite{Hines} as a half-step of the implicit 
backward Euler scheme followed by a generic explicit forward Euler half-step. The implicit half-step being \\ \\[-1.5ex]
\begin{multline}
\label{hCNimplicitStep}
\frac{V_j^{n+\frac{1}{2}}-V_j^n}{k/2} + \Big(\alpha+\sum_i\beta_i\Big)V_j^{n+\frac{1}{2}} -\Big(\alpha E_L+\sum_i\beta_iE_i\Big) ~ ... \hfill \\
\hfill -\Bigg(\frac{V_{j-1}^{n+\frac{1}{2}}-V_j^{n+\frac{1}{2}}}{R_aC_m}+\frac{V_{j+1}^{n+\frac{1}{2}}-V_j^{n+\frac{1}{2}}}{R'_aC_m}\Bigg)=0, 
\end{multline}
\\
\cite{Hines}, which when rearranged becomes
\begin{multline}
\label{hCNstab0}
V_j^{n+\frac{1}{2}} - V_j^n = -\Bigg(\Big(\alpha+\sum_i\beta_i\Big)V_j^{n+\frac{1}{2}}-\Big(\alpha E_L+\sum_i\beta_iE_i\Big) \Bigg.~... \hfill \\
\hfill \Bigg. -\bigg(c_1\Big(V_{j-1}^{n+\frac{1}{2}}-V_j^{n+\frac{1}{2}}\Big)+c_2\Big(V_{j+1}^{n+\frac{1}{2}}-V_j^{n+\frac{1}{2}}\Big)\bigg)\Bigg)\frac{k}{2}. 
\end{multline}
\\
The subsequent HCN explicit half-step was defined as
\begin{multline}
\label{hCNstab1}
V^{n+1}_j=2\left(V^{n+\frac{1}{2}}_j-V^n_j\right)+V^n_j=2V^{n+\frac{1}{2}}_j-V^n_j. \hfill 
\end{multline}
\\
As suggested in \cite{Str} we treat this two-step scheme by defining the following growth factor identities
\begin{multline}
\label{hCNstab2}
V_j^{n+\frac{1}{2}}=\tilde{g}g^n e^{ij\theta} \quad \parbox{1cm} {and} \quad V_j^n=g^n e^{ij\theta} \hfill
\end{multline}
\\
Solving for the HCN scheme's growth factor requires two steps. First we will solve for $\tilde{g}$, the growth factor associated with the 
first, implicit, half-step. To proceed we substitute instances of $V_j^{n+\frac{1}{2}}$ and $V_j^n$ in equation (\ref{hCNstab1}) with their definitions 
in equation (\ref{hCNstab3}). This leads to \\ \\
\begin{multline}
\label{hCNstab3}
\tilde{g}g^ne^{ij\theta}-g^ne^{ij\theta}=-\Bigg(\Big(\alpha+\sum_i\beta_i\Big)\tilde{g}g^ne^{ij\theta} ~ ... \hfill \\ \\[-1.5ex]
-\bigg(c_1\Big(\tilde{g}g^ne^{i(j-1)\theta}-\tilde{g}g^ne^{ij\theta}\Big)
+c_2\Big(\tilde{g}g^ne^{i(j+1)\theta}-\tilde{g}g^ne^{ij\theta}\Big)\bigg)\Bigg)\frac{k}{2}. %\\
\end{multline}
Dividing both sides by $g^ne^{ij\theta}$
\begin{multline}
\label{hCNstab4}
\tilde{g}=1-\tilde{g}\Bigg(\Big(\alpha+\sum_i\beta_i\Big)-\bigg(c_1(e^{-i\theta}-1)+c_2(e^{i\theta}-1)\bigg)\Bigg)\frac{k}{2} \hfill \\ \\[-1.5ex]
\Rightarrow ~\tilde{g}\Bigg(1+\bigg(\Big(\alpha+\sum_i\beta_i\Big) 
-\Big(c_1(\cos{\theta}-i\sin{\theta}-1)+c_2(\cos{\theta}+i\sin{\theta}-1)\Big)\bigg)\frac{k}{2}\Bigg) = 1. \hfill \\
\Rightarrow ~\tilde{g}=2 \Bigg/\bigg(2+\Big(\alpha+\sum_i\beta_i\Big)k+\Big((c_1+c_2)(1-\cos{\theta})-i\left(c_2-c_1\right)\sin{\theta}\Big)k\bigg) \hfill \\
\Rightarrow~\tilde{g}=2 \Bigg/\bigg(2+\Big(\alpha+\sum_i\beta_i\Big)k
+\Big(2(c_1+c_2)\sin^2{\textstyle\left(\frac{\theta}{2}\right)}-i\left(c_2-c_1\right)\sin{\theta}\Big)k\bigg). \hfill 
\end{multline}
Next, to account for growth in the second, explicit, half-step, in terms of the first implicit half-step, we substitute instances of $V_j^{n+\frac{1}{2}}$ 
and $V_j^n$ in equation (\ref{hCNstab2}) with their definitions in equation (\ref{hCNstab3}). This gives us 
\begin{multline}
\label{hCNstab5}
g^{n+1}e^{ij\theta}=2\tilde{g}g^n e^{ij\theta}-g^n e^{ij\theta} \quad \Rightarrow \quad g = 2\tilde{g}-1. \hfill 
\end{multline}
After substituting $\tilde{g}$, representing growth in the first half-step, in equation (\ref{hCNstab5}) with equation (\ref{hCNstab4}), we find 
the scheme's composite growth factor to be
\begin{multline}
\label{hCNstab8}
g(k,\theta)=\frac{2-\Big(K+2L-iM\Big)k}{2+\Big(K+2L-iM\Big)k}, \qquad \text{ where } \hfill \\ \\[-1.5ex]
K=\Big(\alpha+\sum_i\beta_i\Big), ~~ L=\left(c_1+c_2\right)\sin^2{\left(\frac{\theta}{2}\right)}, \text{ and } M=\left(c_2-c_1\right)\sin{\theta}.
\end{multline}
Because the model's endpoints float, basis function spatial derivatives must be zero at the boundaries which implies the only basis functions 
are cosines and we may disregard the imaginary sine term. So we are left with  
\begin{equation}
\label{hCNstab9}
g(k,\theta)=\frac{1-0.5\big(K+2L\big)k}{1+0.5\big(K+2L\big)k}.
\end{equation}
Since $|g(k,\theta)|<1$ for all step sizes, $k$, HCN is unconditionally stable. However, 
\begin{equation}
\label{}
k>\frac{2}{K+2L}~\Rightarrow~ g(k,\theta)<0,
\end{equation}
and as $k$ continues to grow, $g(k,\theta)~\Rightarrow~-1$ from above, which implies the solution will oscillate. Technically the amplitude of 
this oscillation will decrease over time, but only if the model's governing equation did not change. The fact that the governing equation does 
change at every time step leads to oscillations whose magnitudes grow as the solution's rate of change decreases. Oscillation magnitudes 
reach a maximum where the solution magnitude's rate of change approaches zero and changes polarity. This is illustrated in the middle pane 
of Figures \ref{fig:somaStatsWithStd}, \ref{fig:Spiking}-\ref{fig:maxPolarization} and in videos included with supplemental materials.
\subsection{RK21}
\subsubsection{RK21 Accuracy}
\label{subsubsec:RK21accuracy}
Expressed in the context of a generic Hodgkin-Huxley model, we define the format of an RK method as 
\begin{multline}
\label{rk21TruncErr0}
V_j^{n+1}=V_j^n+kF(t_n,V_j^n,k;f) ~~\text{where}~~F(t_n,V_j^n,k;f)=\sum_{i=1}^s b_iK_i, \\ \\[-2.0ex]
K_i=f\big(t_n+c_ik,V_j^n+k\sum_{j=1}^sa_{ij}K_j\big), ~~ f(t_n,V_j^n) = \big(V_t\big)_j^n = A-BV_j^n. 
\end{multline}
Coefficients $a_{ij}, b_i$ and $c_i$ are listed in \mbox{Table \ref{tab:RK21}}, and terms $A,B$ are defined in equations 
(\ref{eqn:SpatiallyDscrtzdCable}) and (\ref{firstOlinODE}). 
\begin{table}[h!]\small
\caption{RK21 Tableau}
\ \\
\centering
\begin{tabular}{l l}
\toprule \\[-1.0ex]
\begin{tabular}{l | l l}
$\displaystyle c_1$ & $a_{11}$  & \\
\hline \\[-2.0ex]
 & $b_1$ & $b_2$
\end{tabular}
&
$\quad \Rightarrow \quad$
\begin{tabular}{l | l l}
$1$ & $1$ & \\
\hline \\[-2.0ex]
 & $\frac{1}{2}$ &  $\frac{1}{2}$
\end{tabular}
\\ \\[-1.0ex]
\bottomrule
\end{tabular}
\label{tab:RK21}
\end{table}
\\
RK21 rendered as a quasi-FD method becomes
\begin{multline}
\label{rk21TruncErr1}
V_j^{n+1}=V_j^n+\frac{k}{2}\Big( K_1+K_2 \Big), ~~ \text{where} \hfill \\ \\[-2.0ex]
K_1=f\big(t_n,V_j^n\big)= A-BV_j^n, ~~ \text{and} \hfill \\ \\[-2.0ex]
K_2= f\big(t_n,V_j^n+kK_1\big)  = A-B\Big(V_j^n+k\big(A-BV_j^n\big)\Big). \hfill
\end{multline}
\\
The RK21 difference operator is 
\begin{multline}
\label{rk21TruncErr2}
P_{h,k}(V) = ~ ... \hfill \\ \\[-2.0ex]
V_j^{n+1}-V_j^n-\frac{k}{2}\bigg( \Big[A-BV_j^n \Big] + \Big[ A-B\Big(V_j^n+k\big(A-BV_j^n\big)\Big) \Big] \bigg)=0, \hfill \\ \\[-2.0ex]
\Rightarrow  V_j^{n+1}-V_j^n-k\Big(A-BV_j^n \Big) + \frac{k^2}{2}\Big(B\big(A-BV_j^n\big)\Big)=0. \hfill 
\end{multline}
\\
Taylor expanding the difference operator gives us
\begin{multline}
\label{rk21TruncErr3}
\bigg(V_j^n+kV_t+\frac{k^2}{2}V_{tt}+\frac{k^3}{3!}V_{ttt}+O(k^4)\bigg) -V_j^n ~...\hfill  \\ 
-k\Big(A-BV_j^n \Big) + \frac{k^2}{2}\Big(B\big(A-BV_j^n\big)\Big)=0.
\end{multline}
Simplifying, dividing by $k$ and substituting $A$ and $B$ with their definitions in equations (\ref{eqn:SpatiallyDscrtzdCable}) and 
(\ref{firstOlinODE}) leads to
\begin{multline}
\label{rk21TruncErr4}
V_t+\frac{k}{2}V_{tt}+\frac{k^2}{3!}V_{ttt}+O(k^3) ~ ... \\ \\[-2.0ex]
-\Bigg(\bigg[\alpha E_L+\sum_i\beta_iE_i +\bigg(\frac{V^n_{j-1}}{R_aC_m} +\frac{V^n_{j+1}}{R'_aC_m}\bigg)\bigg]\Bigg. ~...\\
\Bigg. - \bigg[\alpha +\sum_i \beta_i +\bigg(\frac{1}{R_aC_m} +\frac{1}{R'_aC_m}\bigg) \bigg] V_j^n \Bigg) +\frac{Bk}{2}\Big( A-BV_j^n \Big)=0. 
\end{multline}
\\
Taylor expanding the discretized spatial derivative gives us
\begin{multline}
\label{rk21TruncErr5}
V_t+\frac{k}{2}V_{tt}+\frac{k^2}{3!}V_{ttt}+O(k^3) -\Bigg[\alpha E_L+\sum_i\beta_iE_i-\Big(\alpha+\sum_i\beta_i\Big)V^n_j ~... \Bigg. \hfill  \\ \\[-2.0ex]
+\frac{1}{R_aC_m}\Bigg(\bigg[V_j^n-hV_x +\frac{(-h)^2}{2}V_{xx}+\frac{(-h)^3}{3!}V_{xxx}+\frac{(-h)^4}{4!}V_{xxxx}+O(h^5)\bigg] -V_j^n \Bigg)~ ... \hfill \\ \\[-2.0ex]
\Bigg. +\frac{1}{R'_aC_m}\Bigg(\bigg[V_j^n+hV_x+\frac{h^2}{2}V_{xx}+\frac{h^3}{3!}V_{xxx}+\frac{h^4}{4!}V_{xxxx}+O(h^5)\bigg]-V_j^n\Bigg)\Bigg]~... \hfill \\ \\[-2.0ex]
+\frac{Bk}{2}\Big( A-BV^n_j \Big)=0.
\end{multline}
Simplifying, applying the identity in equation (\ref{specificAxialResDefs}), leaves us with
\begin{multline}
\label{rk21TruncErr6}
V_t+\frac{k}{2}V_{tt}+\frac{k^2}{6}V_{ttt}
-\Bigg[\Big(\alpha E_L+\sum_i\beta_iE_i\Big)-\Big(\alpha+\sum_i\beta_i\Big)V^n_j\Bigg. ~... \hfill  \\
+\frac{a^2/a}{2r_Lc_mh^2}\bigg(-hV_x +\frac{h^2}{2}V_{xx}-\frac{h^3}{6}V_{xxx}+\frac{h^4}{24}V_{xxxx} \bigg)~ ... \\
\Bigg. +\frac{a'^2/a}{2r_Lc_mh^2}\bigg(hV_x+\frac{h^2}{2}V_{xx}+\frac{h^3}{3!}V_{xxx}+\frac{h^4}{4!}V_{xxxx}\bigg) \Bigg]~...\\
+\frac{Bk}{2}\Big( A-BV_j^n \Big)=O(k^3,h^3).   
\end{multline}
Rearranging and eliminating third order terms gives us
\begin{multline}
\label{rk21TruncErr7}
V_t+\frac{V_{tt}}{2}k+\frac{V_{ttt}}{6}k^2 -\Bigg(\Big(\alpha E_L+\sum_i\beta_iE_i\Big)-\Big(\alpha+\sum_i\beta_i\Big)V^n_j  \Bigg.~... \hfill \\
\Bigg. +\gamma\bigg[\Big(a'^2+a^2\Big)\bigg(\frac{V_{xx}}{2}+\frac{V_{xxxx}}{24}h^2 \bigg) 
+\Big(a'^2-a^2\Big)\bigg(\frac{V_x}{h}+\frac{V_{xxx}}{6}h\bigg)\bigg]\Bigg)~ ... \\
+\frac{Bk}{2}\Big( A-BV_j^n \Big)=0. 
\end{multline}
Replacing terms with identities defined in equation (\ref{identitiesForTEAndStab}) leads to
\begin{multline}
V_t+\frac{V_{tt}}{2}k+\frac{V_{ttt}}{6}k^2 -\Bigg(\Big(\alpha E_L+\sum_i\beta_iE_i\Big)- \Big(\alpha+\sum_i\beta_i\Big)V^n_j  \Bigg.~... \hfill \\
\Bigg. +\gamma\bigg[\Big(\overline{a}^2V_{xx}+\big(\overline{a}^2\big)_xV_x\Big)+\Big(\overline{a}^2V_{xxxx} 
+2\big(\overline{a}^2\big)_xV_{xxx}\Big)\frac{h^2}{12} \bigg] \Bigg) %~ ... \\
+\frac{Bk}{2}\Big( A-BV_j^n \Big) =0. \\ \\
\Rightarrow V_t+\frac{V_{tt}}{2}k+\frac{V_{ttt}}{6}k^2 ~ ... \hfill \\
-\Bigg(\Big(\alpha E_L+\sum_i\beta_iE_i\Big) - \Big(\alpha+\sum_i\beta_i\Big)V^n_j \Bigg.~... \hfill \\
\Bigg. +\gamma\bigg[\big(\overline{a}^2V_x\big)_x+\Big(\overline{a}^2V_{xxxx} +2\big(\overline{a}^2\big)_xV_{xxx}\Big)\frac{h^2}{12} \bigg] \Bigg) ~ ... \\
+\frac{Bk}{2}\Bigg( \Big(\alpha+\sum_i\beta_i\Big)V^n_j - \Big(\alpha E_L+\sum_i\beta_iE_i\Big)  \Bigg.~... \hfill \\
\Bigg. +\gamma\bigg[\big(\overline{a}^2V_x\big)_x+\Big(\overline{a}^2V_{xxxx} +2\big(\overline{a}^2\big)_xV_{xxx}\Big)\frac{h^2}{12} \bigg] \Bigg) =0.
\end{multline}
\begin{multline}
\label{rk21TruncErr8}
\Rightarrow V_t+\Big(\alpha+\sum_i\beta_i\Big)V^n_j - \Big(\alpha E_L+\sum_i\beta_iE_i\Big)-\gamma\big(\overline{a}^2V_x\big)_x~... \hfill \\
+\frac{V_{tt}}{2}k+\frac{V_{ttt}}{6}k^2-\gamma\Big(\overline{a}^2V_{xxxx} +2\big(\overline{a}^2\big)_xV_{xxx}\Big)\frac{h^2}{12} ~ ... \\
+\frac{Bk}{2}\bigg(V_t+\gamma\Big(\overline{a}^2V_{xxxx} +2\big(\overline{a}^2\big)_xV_{xxx}\Big)\frac{h^2}{12}\bigg) =0.
\end{multline}
\\
Discarding the third order $kh^2$ term we are left with
\begin{multline}
\label{rk21TruncErr9}
V_t+\Big(\alpha+\sum_i\beta_i\Big)V^n_j - \Big(\alpha E_L+\sum_i\beta_iE_i\Big)-\gamma\big(\overline{a}^2V_x\big)_x~... \hfill \\
+\Big(BV_t+V_{tt}\Big)\frac{k}{2}+\frac{V_{ttt}}{6}k^2-\gamma\Big(\overline{a}^2V_{xxxx} +2\big(\overline{a}^2\big)_xV_{xxx}\Big)\frac{h^2}{12}=0.
\end{multline}
\\
From equation (\ref{eqn:Vtt}) we know that $V_{tt}$ is not simply $-BV_t$, and therefore the truncation error for the RK21 method is
\begin{multline}
\label{rk21TruncErr10}
P_{k,h}(V)-P(V) = ~ ... \hfill \\
\Big(BV_t+V_{tt}\Big)\frac{k}{2}+\frac{V_{ttt}}{6}k^2-\gamma\Big(\overline{a}^2V_{xxxx} +2\big(\overline{a}^2\big)_xV_{xxx}\Big)\frac{h^2}{12},
\end{multline}
which is only first order in time, second in space.
\subsubsection{RK21 Stability}
\label{subsubsec:RK21stability}
Starting from equation (\ref{rk21TruncErr2}) we have
\begin{multline}
\label{rk21Stab0}
P_{h,k}(V) = ~ ... \hfill \\ \\[-2.0ex]
V_j^{n+1}-V_j^n-\frac{k}{2}\bigg( \Big[A-BV_j^n \Big] + \Big[ A-B\Big(V_j^n+k\big(A-BV_j^n\big)\Big) \Big] \bigg)=0, \hfill \\ \\[-2.0ex]
\Rightarrow  V_j^{n+1}-V_j^n-k\Big(A-BV_j^n \Big) + \frac{k^2}{2}\Big(B\big(A-BV_j^n\big)\Big)=0. \hfill 
\end{multline}
After dropping constants, replacing dependent variables with their phased growth factor equivalents and combining like terms we have
\begin{multline}
\label{rk21Stab1}
g^{n+1}e^{ij\theta}= g^ne^{ij\theta}~...\hfill \\
+k\Bigg(-\Big(\alpha+\sum_i\beta_i\Big)g^ne^{ij\theta}+\bigg[\frac{g^ne^{i(j-1)\theta}-g^ne^{ij\theta}}{R_aC_m}
+\frac{g^ne^{i(j+1)\theta}-g^ne^{ij\theta}}{R'_aC_m}\bigg]\Bigg)~... \\
-\frac{k^2}{2}B\Bigg(-\Big(\alpha+\sum_i\beta_i\Big)g^ne^{ij\theta}+\bigg[\frac{g^ne^{i(j-1)\theta}-g^ne^{ij\theta}}{R_aC_m}
+\frac{g^ne^{i(j+1)\theta}-g^ne^{ij\theta}}{R'_aC_m}\bigg]\Bigg)
\end{multline}
\begin{multline}
\label{rk21Stab2}
g= 1+k\Bigg(-\Big(\alpha+\sum_i\beta_i\Big)+\bigg[\frac{e^{-i\theta}-1}{R_aC_m}+\frac{e^{i\theta}-1}{R'_aC_m}\bigg]\Bigg)~... \\
-\frac{k^2}{2}B\Bigg(-\Big(\alpha+\sum_i\beta_i\Big)+\bigg[\frac{e^{-i\theta}-1}{R_aC_m}+\frac{e^{i\theta}-1}{R'_aC_m}\bigg]\Bigg)
\end{multline}
\begin{multline}
\label{rk21Stab3}
g= 1+k\Bigg(-\Big(\alpha+\sum_i\beta_i\Big)+\bigg[\frac{\cos{\theta}-i\sin{\theta}-1}{R_aC_m}+\frac{\cos{\theta}+i\sin{\theta}-1}{R'_aC_m}\bigg]\Bigg)~... \\
-\frac{k^2}{2}B\Bigg(-\Big(\alpha+\sum_i\beta_i\Big)+\bigg[\frac{\cos{\theta}-i\sin{\theta}-1}{R_aC_m}+\frac{\cos{\theta}+i\sin{\theta}-1}{R'_aC_m}\bigg]\Bigg)
\end{multline}
\begin{multline}
\label{rk21Stab4}
g= 1+k\Bigg(-\Big(\alpha+\sum_i\beta_i\Big)+\Big((c_1+c_2)(\cos{\theta}-1)+i(c_2-c_1)\sin{\theta}\Big)\Bigg)~... \\
-\frac{k^2}{2}B\Bigg(-\Big(\alpha+\sum_i\beta_i\Big)+\Big((c_1+c_2)(\cos{\theta}-1)+i(c_2-c_1)\sin{\theta}\Big)\Bigg)
\end{multline}
\begin{multline}
\label{rk21Stab5}
g=1-k\Big(K+2L-iM \Big)+\frac{k^2}{2}B\Big(K+2L-iM\Big), \qquad \parbox{1cm}{where} \hfill \\ \\[-1.5ex]
K=\Big(\alpha+\sum_i\beta_i\Big),~ L=\big(c_1+c_2\big)\sin^2{\textstyle\left(\frac{\theta}{2}\right)}, ~ \text{and } M=(c_2-c_1)\sin{\theta}. \hfill 
\end{multline}
Because the model's endpoints float, basis function spatial derivatives must be zero at the boundaries which implies the only basis functions 
are cosines and we may disregard the imaginary sine term. We now seek values of $g$ such that \\ \\
\begin{multline}
\label{rk21Stab6}
|g(k,\theta)|^2 \le 1 \Rightarrow \bigg(1-k\big(K+2L\big)+\frac{k^2}{2}B\big(K+2L\big) \bigg)^2\le 1 \hfill \\
\Rightarrow k\Bigg[\bigg(\frac{BP}{2}\bigg)^2k^3-BP^2k^2+(BP+P^2)k-2P\Bigg]\le 0  \\
\Rightarrow \bigg(\frac{BP}{2}\bigg)^2k\Bigg[k^3-\frac{4}{B}k^2+\frac{4(B+P)}{B^2P}k-\frac{8}{B^2P}\Bigg]\le 0 \\
\Rightarrow Q(k) = k^3-\frac{4}{B}k^2+\frac{4(B+P)}{B^2P}k-\frac{8}{B^2P}\le 0, 
\end{multline}
where $P=\big(K+2L\big)$. The only variable terms in $B$ and $P$ are $\sum_i \beta_i$ and $\theta$. Although a bit involved, it would 
be possible to plot $Q(k)$ over the domain $[k=0,~...~,~50\mu\text{sec}]$X$[\text{span of }\sum_i \beta_i]$ for multiple values of $\theta$.  \\ \\
Instead, we consider the Butcher perspective of stability for the model ODE problem, $y'=\lambda y$, which is \cite[pg.100]{Butcher}
\begin{equation}
\label{rk21Stab7}
R(z)=1+z+\frac{1}{2}z^2, \qquad \text{where} \quad z=\lambda\cdot k.
\end{equation}
To be stable,
\begin{equation}
\label{rk21Stab8}
\big|R(z)\big|\le 1 \quad \Rightarrow ~ -2<Re(z)<0 \quad \Rightarrow ~ 0 < k < \frac{2}{|\lambda |}.
\end{equation}
A simple change of variables turns our ODE, $V'=A-BV$, into $W'=-BW$. Clearly, $\lambda$ in equation (\ref{rk21Stab8}) is $-B$ and so for the 
numerical integration to remain stable, $k\le\frac{2}{B}$ where $B$ was defined in Table \ref{tab:schemeStability}. Unfortunately, Traditional RK 
stability analysis treats the numerical method as though it were integrating an ODE, not a spatially discretized PDE. Without $\theta$ part of the 
step size limit's expression it is significantly underestimated as shown in Figure \ref{fig:RKstepsizeLimits}.
\subsection{RK41}
\subsubsection{RK41 Accuracy}
\label{subsubsec:RK41accuracy}
Building on the development for RK21 in \ref{subsubsec:RK21accuracy}, 
\begin{table}[h!]\small
\caption{RK41 Tableau}
\ \\
\centering
\begin{tabular}{l l}
\toprule \\[-1.0ex]
\begin{tabular}{l | l l l l}
$\displaystyle c_1$ &                 &                &                  & \\
$\displaystyle c_2$ & $a_{21}$  &                &                  & \\
$\displaystyle c_3$ & $a_{31}$  & $a_{32}$ &                  & \\
$\displaystyle c_4$ & $a_{41}$  & $a_{42}$  & $a_{43}$  &  \\
\hline \\[-2.0ex]
 & $b_1$ & $b_2$ & $b_3$ & $b_4$
\end{tabular}
&
$\quad \Rightarrow \quad$
\begin{tabular}{l | l l l l}
$0$               &                      &        &           & \\
$\frac{1}{2}$ & $\frac{1}{2}$ &        &            & \\
$\frac{1}{2}$ & $0$               & $\frac{1}{2}$ &        & \\
$1$               & $0$               & $0$               & $1$ &  \\
\hline \\[-2.0ex]
 & $\frac{1}{6}$ &  $\frac{1}{3}$ & $\frac{1}{3}$ &  $\frac{1}{6}$
\end{tabular}
\\ \\[-1.0ex]
\bottomrule
\end{tabular}
\label{tab:RK41}
\end{table}
we render RK41, defined by the tableau in Table \ref{tab:RK41}, as 
\begin{multline}
\label{rk41TruncErr1}
V_j^{n+1}=V_j^n+\frac{k}{6}\Big(K_1+2K_2+2K_3+K_4 \Big), ~~ \text{where} \hfill \\ \\[-2.0ex]
K_1=f\big(t_n,V_j^n\big)= A-BV_j^n, ~~ \text{and} \hfill \\ \\[-2.0ex]
K_2= f\big(t_n,V_j^n+\frac{k}{2}K_1\big)  = A-B\Big(V_j^n+\frac{k}{2}\big(A-BV_j^n\big)\Big), \hfill \\ \\[-2.0ex]
K_3= f\big(t_n,V_j^n+\frac{k}{2}K_2\big)  = A-B\bigg(V_j^n+\frac{k}{2}\Big(A-B\big(V_j^n+\frac{k}{2}(A-BV_j^n)\big)\Big)\bigg), \hfill \\ \\[-2.0ex]
K_4= f\big(t_n,V_j^n+kK_3\big) = ~ ... \hfill \\
A-B\Bigg(V_j^n+k\bigg(A-B\Big(V_j^n+\frac{k}{2}\big(A-B(V_j^n+\frac{k}{2}[A-BV_j^n])\big)\Big)\bigg)\Bigg).
\end{multline}
\\
The RK41 quasi-FD method difference operator is 
\begin{multline}
\label{rk41TruncErr2}
P_{h,k}(V) = V_j^{n+1}-V_j^n ~ ... \hfill \\
-\frac{k}{6}\Bigg( \bigg[A-BV_j^n \bigg] + 2\bigg[ A-B\Big(V_j^n+\frac{k}{2}\big[A-BV_j^n\big]\Big) \bigg] ~ ... \hfill \\
+2\Bigg[ A-B\bigg(V_j^n+\frac{k}{2}\Big[A-B\Big(V_j^n+\frac{k}{2}[A-BV_j^n]\Big)\Big]\bigg) \Bigg] ~ ... \\
+\Bigg[A-B\bigg(V_j^n+k\bigg[A-B\Big(V_j^n+\frac{k}{2}\Big[A-B\big(V_j^n+\frac{k}{2}[A-BV_j^n]\big)\Big]\Big)\bigg]\bigg) \Bigg] \Bigg)=0, \\ \\
\Rightarrow  V_j^{n+1}-V_j^n-k\Big(A-BV_j^n \Big) + \frac{k^2}{2}\Big(B\big(A-BV_j^n\big)\Big) ~ ... \hfill \\
-\frac{k^3}{3!}\Big(B^2\big(A-BV_j^n\big)\Big)+\frac{k^4}{4!}\Big(B^3\big(A-BV_j^n\big)\Big)=0.  
\end{multline}
\\
Taylor expanding the difference operator gives us
\begin{multline}
\label{rk41TruncErr3}
\bigg(V_j^n+kV_t+\frac{k^2}{2}V_{tt}+\frac{k^3}{3!}V_{ttt}+O(k^4)\bigg) -V_j^n -k\Big(A-BV_j^n \Big) ~... \\ 
+ \frac{k^2}{2}\Big(B\big(A-BV_j^n\big)\Big) -\frac{k^3}{3!}\Big(B^2\big(A-BV_j^n\big)\Big)+\frac{k^4}{4!}\Big(B^3\big(A-BV_j^n\big)\Big)=0.
\end{multline}
Simplifying, dividing by $k$ and substituting just the first occurance of $A$ and $B$ with their definitions from equations 
(\ref{eqn:SpatiallyDscrtzdCable}) and (\ref{firstOlinODE}) leads to
\begin{multline}
\label{rk41-4}
V_t+\frac{k}{2}V_{tt}+\frac{k^2}{3!}V_{ttt}+O(k^3) ~ ... \hfill \\ \\[-2.0ex]
-\Bigg(\bigg[\alpha E_L+\sum_i\beta_iE_i +\bigg(\frac{V^n_{j-1}}{R_aC_m} +\frac{V^n_{j+1}}{R'_aC_m}\bigg)\bigg]\Bigg. ~...\hfill \\
\Bigg. - \bigg[\alpha +\sum_i \beta_i +\bigg(\frac{1}{R_aC_m} +\frac{1}{R'_aC_m}\bigg) \bigg] V_j^n \Bigg) ~ ... \\
+\frac{Bk}{2}\Big(A - BV_j^n\Big) -\frac{(Bk)^2}{6}\Big(A - BV_j^n\Big) + \frac{(Bk)^3}{24}\Big(A - BV_j^n\Big)=0. 
\end{multline}
\\
Taylor expanding the spatial derivative gives us
\begin{multline}
\label{rk41TruncErr5}
V_t+\frac{k}{2}V_{tt}+\frac{k^2}{3!}V_{ttt}+O(k^3) 
-\Bigg[\alpha E_L+\sum_i\beta_iE_i-\Big(\alpha+\sum_i\beta_i\Big)V^n_j ~... \Bigg. \hfill  \\ \\[-2.0ex]
+\frac{1}{R_aC_m}\Bigg(\bigg[V_j^n-hV_x +\frac{(-h)^2}{2}V_{xx}+\frac{(-h)^3}{3!}V_{xxx}+\frac{(-h)^4}{4!}V_{xxxx}+O(h^5)\bigg] -V_j^n \Bigg)~ ... \hfill \\ \\[-2.0ex]
\Bigg. +\frac{1}{R'_aC_m}\Bigg(\bigg[V_j^n+hV_x+\frac{h^2}{2}V_{xx}+\frac{h^3}{3!}V_{xxx}+\frac{h^4}{4!}V_{xxxx}+O(h^5)\bigg]-V_j^n\Bigg)\Bigg]~... \hfill \\ \\[-2.0ex]
+\frac{Bk}{2}\Big(A - BV_j^n\Big) -\frac{(Bk)^2}{6}\Big(A - BV_j^n\Big) + \frac{(Bk)^3}{24}\Big(A - BV_j^n\Big)=0. 
\end{multline}
\\
Simplifying, applying the identity in equation (\ref{specificAxialResDefs}), leaves us with
\begin{multline}
\label{rk41TruncErr6}
V_t+\frac{k}{2}V_{tt}+\frac{k^2}{6}V_{ttt}
-\Bigg[\Big(\alpha E_L+\sum_i\beta_iE_i\Big)-\Big(\alpha+\sum_i\beta_i\Big)V^n_j \Bigg. ~... \hfill  \\
+\frac{a^2/a}{2r_Lc_mh^2}\bigg(-hV_x +\frac{h^2}{2}V_{xx}-\frac{h^3}{6}V_{xxx}+\frac{h^4}{24}V_{xxxx} \bigg)~ ... \\
\Bigg. +\frac{a'^2/a}{2r_Lc_mh^2}\bigg(hV_x+\frac{h^2}{2}V_{xx}+\frac{h^3}{3!}V_{xxx}+\frac{h^4}{4!}V_{xxxx}\bigg) \Bigg]~...\\
+\frac{Bk}{2}\Big(A - BV_j^n\Big) -\frac{(Bk)^2}{6}\Big(A - BV_j^n\Big) + \frac{(Bk)^3}{24}\Big(A - BV_j^n\Big)=O(k^3,h^3).   
\end{multline}
Rearranging and eliminating third order terms gives us
\begin{multline}
\label{rk41TruncErr7}
V_t+\frac{V_{tt}}{2}k+\frac{V_{ttt}}{6}k^2 -\Bigg(\Big(\alpha E_L+\sum_i\beta_iE_i\Big)-\Big(\alpha+\sum_i\beta_i\Big)V^n_j  \Bigg.~... \hfill \\
\Bigg. +\gamma\bigg[\Big(a'^2+a^2\Big)\bigg(\frac{V_{xx}}{2}+\frac{V_{xxxx}}{24}h^2 \bigg) 
+\Big(a'^2-a^2\Big)\bigg(\frac{V_x}{h}+\frac{V_{xxx}}{6}h\bigg)\bigg]\Bigg)~ ... \\
+\frac{Bk}{2}\Big( A-BV_j^n \Big)-\frac{(Bk)^2}{6}\Big(A - BV_j^n\Big) + \frac{(Bk)^3}{24}\Big(A - BV_j^n\Big)=0. 
\end{multline}
After replacing terms above with identities defined in equation (\ref{identitiesForTEAndStab}) we have 
\begin{multline}
\label{rk41TruncErr8}
V_t+\frac{V_{tt}}{2}k+\frac{V_{ttt}}{6}k^2 -\Bigg(\Big(\alpha E_L+\sum_i\beta_iE_i\Big)- \Big(\alpha+\sum_i\beta_i\Big)V^n_j  \Bigg.~... \hfill \\
\Bigg. +\gamma\bigg[\Big(\overline{a}^2V_{xx}+\big(\overline{a}^2\big)_xV_x\Big)+\Big(\overline{a}^2V_{xxxx} 
+2\big(\overline{a}^2\big)_xV_{xxx}\Big)\frac{h^2}{12} \bigg] \Bigg) ~ ... \\
+\frac{Bk}{2}\Big( A-BV_j^n \Big)-\frac{(Bk)^2}{6}\Big(A - BV_j^n\Big) + \frac{(Bk)^3}{24}\Big(A - BV_j^n\Big) =0.
\end{multline}
\begin{multline}
\label{rk41TruncErr9}
\Rightarrow V_t+\frac{V_{tt}}{2}k+\frac{V_{ttt}}{6}k^2 ~ ... \hfill \\
-\Bigg(\Big(\alpha E_L+\sum_i\beta_iE_i\Big) - \Big(\alpha+\sum_i\beta_i\Big)V^n_j \Bigg.~... \hfill \\
\Bigg. +\gamma\bigg[\big(\overline{a}^2V_x\big)_x+\Big(\overline{a}^2V_{xxxx} +2\big(\overline{a}^2\big)_xV_{xxx}\Big)\frac{h^2}{12} \bigg] \Bigg) ~ ... \\
+\frac{Bk}{2}\Bigg(\Big(\alpha E_L+\sum_i\beta_iE_i\Big) - \Big(\alpha+\sum_i\beta_i\Big)V^n_j \Bigg.~... \hfill \\
\Bigg. +\gamma\bigg[\big(\overline{a}^2V_x\big)_x+\Big(\overline{a}^2V_{xxxx} +2\big(\overline{a}^2\big)_xV_{xxx}\Big)\frac{h^2}{12} \bigg] \Bigg)~ ... \\
-\frac{(Bk)^2}{6}\Bigg(\Big(\alpha E_L+\sum_i\beta_iE_i\Big) - \Big(\alpha+\sum_i\beta_i\Big)V^n_j \Bigg.~... \hfill \\
\Bigg. +\gamma\bigg[\big(\overline{a}^2V_x\big)_x+\Big(\overline{a}^2V_{xxxx} +2\big(\overline{a}^2\big)_xV_{xxx}\Big)\frac{h^2}{12} \bigg] \Bigg)~ ... \\
+\frac{(Bk)^3}{24}\Bigg(\Big(\alpha E_L+\sum_i\beta_iE_i\Big) - \Big(\alpha+\sum_i\beta_i\Big)V^n_j \Bigg.~... \hfill \\
\Bigg. +\gamma\bigg[\big(\overline{a}^2V_x\big)_x+\Big(\overline{a}^2V_{xxxx} +2\big(\overline{a}^2\big)_xV_{xxx}\Big)\frac{h^2}{12} \bigg] \Bigg) =0.
\end{multline}
Replacing terms above with identities defined in equation (\ref{eqn:truncErrIdent}) gives us
\begin{multline}
\label{rk41TruncErr10}
V_t+\Big(\alpha+\sum_i\beta_i\Big)V^n_j - \Big(\alpha E_L+\sum_i\beta_iE_i\Big)-\gamma\big(\overline{a}^2V_x\big)_x~... \hfill \\
+\frac{V_{tt}}{2}k+\frac{V_{ttt}}{6}k^2-\gamma\Big(\overline{a}^2V_{xxxx} +2\big(\overline{a}^2\big)_xV_{xxx}\Big)\frac{h^2}{12} ~ ... \\
+\frac{Bk}{2}\bigg(V_t+\gamma\Big(\overline{a}^2V_{xxxx} +2\big(\overline{a}^2\big)_xV_{xxx}\Big)\frac{h^2}{12}\bigg) ~ ... \hfill \\
-\frac{(Bk)^2}{6}\bigg(V_t+\gamma\Big(\overline{a}^2V_{xxxx} +2\big(\overline{a}^2\big)_xV_{xxx}\Big)\frac{h^2}{12}\bigg) ~ ... \\
+\frac{(Bk)^3}{24}\bigg(V_t+\gamma\Big(\overline{a}^2V_{xxxx} +2\big(\overline{a}^2\big)_xV_{xxx}\Big)\frac{h^2}{12}\bigg) =0.
\end{multline}
\\
Discarding third order terms leaves us with
\begin{multline}
\label{rk41TruncErr11}
V_t+\Big(\alpha+\sum_i\beta_i\Big)V^n_j - \Big(\alpha E_L+\sum_i\beta_iE_i\Big)-\gamma\big(\overline{a}^2V_x\big)_x~... \hfill \\
+\Big(BV_t+V_{tt}\Big)\frac{k}{2}+\frac{V_{ttt}}{6}k^2-\gamma\Big(\overline{a}^2V_{xxxx} +2\big(\overline{a}^2\big)_xV_{xxx}\Big)\frac{h^2}{12}=0.
\end{multline}
\\
Truncation error for the RK41 method is therefore,
\begin{multline}
\label{rk41TruncErr12}
P_{k,h}(V)-P(V) = ~ ... \hfill \\
\Big(BV_t+V_{tt}\Big)\frac{k}{2}+\frac{V_{ttt}}{6}k^2-\gamma\Big(\overline{a}^2V_{xxxx} +2\big(\overline{a}^2\big)_xV_{xxx}\Big)\frac{h^2}{12},
\end{multline}
\\
first order in time, second in space.
\subsubsection{RK41 Stability}
\label{subsubsec:RK41stability}
Starting from equation (\ref{rk41TruncErr2}) we have
\begin{multline}
\label{rk41Stab0}
V_j^{n+1}-V_j^n-k\Big(A-BV_j^n \Big) + \frac{k^2}{2}\Big(B\big(A-BV_j^n\big)\Big) ~ ... \hfill \\
-\frac{k^3}{3!}\Big(B^2\big(A-BV_j^n\big)\Big)+\frac{k^4}{4!}\Big(B^3\big(A-BV_j^n\big)\Big)=0.  
\end{multline}
After dropping constants, replacing dependent variables with their phased growth factor equivalents and combining like terms we have
\begin{multline}
\label{rk41Stab1}
g^{n+1}e^{ij\theta}= g^ne^{ij\theta}~...\hfill \\
+k\Bigg(-\Big(\alpha+\sum_i\beta_i\Big)g^ne^{ij\theta}+\bigg[\frac{g^ne^{i(j-1)\theta}-g^ne^{ij\theta}}{R_aC_m}
+\frac{g^ne^{i(j+1)\theta}-g^ne^{ij\theta}}{R'_aC_m}\bigg]\Bigg)~... \\
-\frac{k^2}{2}B\Big(A-BV_j^n\Big) +\frac{k^3}{6}B^2\Big(A-BV_j^n\Big)-\frac{k^4}{24}B^3\Big(A-BV_j^n\Big)
\end{multline}
\begin{multline}
\label{rk41Stab2}
g= 1+k\Bigg(-\Big(\alpha+\sum_i\beta_i\Big)+\bigg[\frac{e^{-i\theta}-1}{R_aC_m}+\frac{e^{i\theta}-1}{R'_aC_m}\bigg]\Bigg)~... \\
-\frac{k^2}{2}B\Big(A-BV_j^n\Big)+\frac{k^3}{6}B^2\Big(A-BV_j^n\Big)-\frac{k^4}{24}B^3\Big(A-BV_j^n\Big)
\end{multline}
\begin{multline}
\label{rk41Stab3}
g= 1+k\Bigg(-\Big(\alpha+\sum_i\beta_i\Big)+\bigg[\frac{\cos{\theta}-i\sin{\theta}-1}{R_aC_m}+\frac{\cos{\theta}+i\sin{\theta}-1}{R'_aC_m}\bigg]\Bigg)~... \\
-\frac{k^2}{2}B\Big(A-BV_j^n\Big)+\frac{k^3}{6}B^2\Big(A-BV_j^n\Big)-\frac{k^4}{24}B^3\Big(A-BV_j^n\Big)
\end{multline}
\begin{multline}
\label{rk41Stab4}
g= 1+k\Bigg(-\Big(\alpha+\sum_i\beta_i\Big)+\Big((c_1+c_2)(\cos{\theta}-1)+i(c_2-c_1)\sin{\theta}\Big)\Bigg)~... \\
-\frac{k^2}{2}B\Bigg(-\Big(\alpha+\sum_i\beta_i\Big)+\Big((c_1+c_2)(\cos{\theta}-1)+i(c_2-c_1)\sin{\theta}\Big)\Bigg) \\
+\frac{k^3}{6}B^2\Bigg(-\Big(\alpha+\sum_i\beta_i\Big)+\Big((c_1+c_2)(\cos{\theta}-1)+i(c_2-c_1)\sin{\theta}\Big)\Bigg) \\
-\frac{k^4}{24}B^3\Bigg(-\Big(\alpha+\sum_i\beta_i\Big)+\Big((c_1+c_2)(\cos{\theta}-1)+i(c_2-c_1)\sin{\theta}\Big)\Bigg)
\end{multline}
\begin{multline}
\label{rk41Stab5}
g=1-k\Big(K+2L-iM \Big) ~ ... \hfill \\
+\frac{k^2}{2}B\Big(K+2L-iM\Big)-\frac{k^3}{6}B^2\Big(K+2L-iM\Big)+\frac{k^4}{24}B^3\Big(K+2L-iM\Big) \\ \\
\text{where} \hfill \\ \\[-1.5ex]
K=\Big(\alpha+\sum_i\beta_i\Big),~ L=\big(c_1+c_2\big)\sin^2{\textstyle\left(\frac{\theta}{2}\right)}, ~ \text{and } M=(c_2-c_1)\sin{\theta}. \hfill 
\end{multline}
Because the model's endpoints float, basis function spatial derivatives must be zero at the boundaries which implies the only basis functions 
are cosines and we may disregard the imaginary sine term. We now seek values of $g$ such that \\ \\
\begin{multline}
\label{rk41Stab6}
|g(k,\theta)|^2\le1 ~ ... \hfill \\
\Rightarrow \bigg(1-k\big(K+2L\big)+\frac{k^2}{2}B\big(K+2L\big)-\frac{k^3}{6}B^2\big(K+2L\big)+\frac{k^4}{24}B^3\big(K+2L\big)\bigg)^2\le1 \\ \\[-1.5ex]
\Rightarrow \bigg(\frac{B^6P^2}{576}k^8+\frac{B^5P^2}{72}k^7+\frac{5B^4P^2}{72}k^6+\frac{B^3P^2}{12}k^5+\frac{B^3P-B^2P^2}{12}k^4\bigg. ~... \hfill \\ \\[-1.5ex]
\bigg. +\Big(\frac{B^2P}{3}-BP^2\Big)k^3+\big(BP+P^2\big)k^2-2Pk+1\bigg)\le1 \\ \\[-1.5ex]
\Rightarrow \frac{B^6P^2}{576}k\Bigg(k^7+\frac{8}{B}k^6+\frac{40}{B^2}k^5+\frac{48}{B^3}k^4+\frac{48(B-P)}{B^4P}k^3 \Bigg.~ ... \hfill \\
\Bigg. +\frac{192(B-3P)}{B^5P}k^2+\frac{576(B+P)}{B^6P}k-\frac{1152}{B^6P}\Bigg)\le0,
\end{multline}
where $P=\big(K+2L\big)$. The Butcher perspective of stability, which proved satisfactory in equations (\ref{rk21Stab7}) and (\ref{rk21Stab8}), 
is far more tractable, simply 
\begin{equation}
\label{eqn:RK4R(z)}
R(z)=1+z+\frac{1}{2}z^2+\frac{1}{6}z^3+\frac{1}{24}z^4, \qquad \text{where} \quad z=\lambda\cdot k.
\end{equation}
To be stable,
\begin{equation}
\label{eqn:RK41Stab}
R(z)\le 1 \qquad \Rightarrow ~ \sim -2.7853<Re(z)<0, \quad \Rightarrow ~ 0 < k < \frac{2.7853}{|\lambda |}.
\end{equation}
A simple change of variables turns our ODE, $V'=A-BV$, into $W'=-BW$. Clearly, $\lambda$ in equation (\ref{rk21Stab8}) is $-B$ and so for the 
numerical integration to remain stable, $k\le\frac{2.7853}{B}$. Unfortunately, by not considering the problem as a PDE, this step size limit is significantly 
underestimated as shown in \mbox{Figure \ref{fig:RKstepsizeLimits}}.
%\clearpage
\section{Neural Model Diagrams}
\label{Appdx:neuralModelDiags}
%\ \\
\begin{figure}[h!]
\begin{center}
\includegraphics[height=3.5in,width=4.5in]{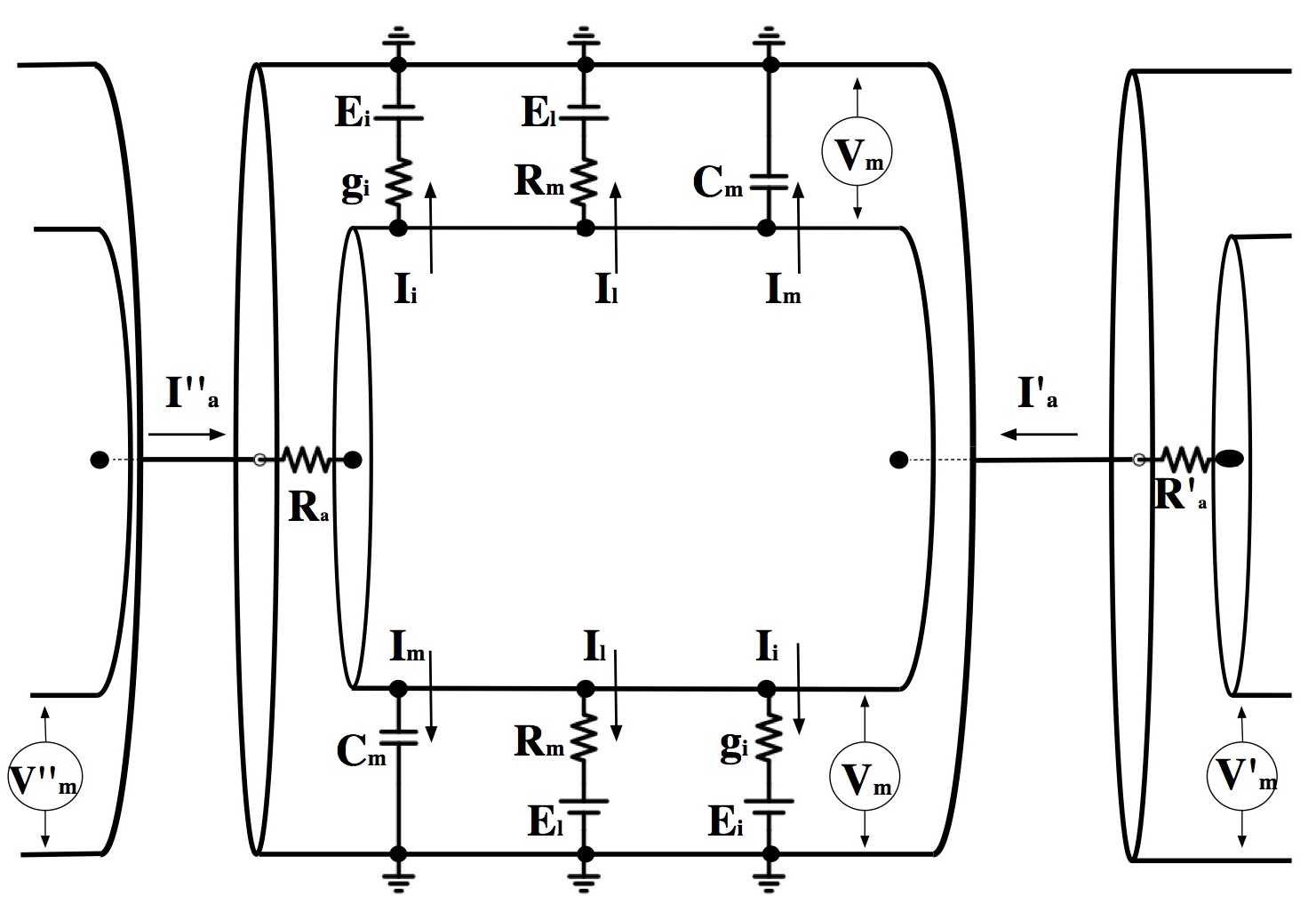}
\end{center}
\caption{Circuit diagram for an asymmetric neural model compartment, using only one sided axial resistance, whose parent compartment 
is on the left and child compartment on the right.}
\label{fig:modelCmpts}
\end{figure} 
\begin{figure}[h!]
\begin{center}
\includegraphics[height=6.5in,width=4.0in]{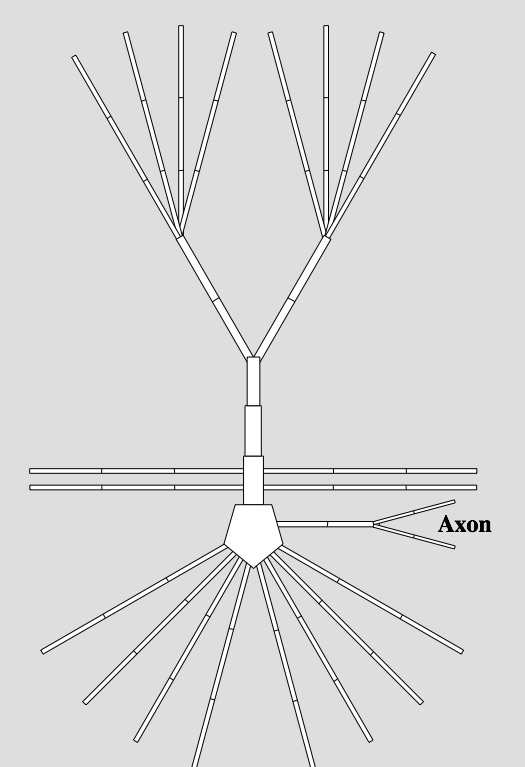}
\end{center}
\caption{Planar view of P23RS compartment structure.}
\label{fig:p23rsDiagram}
\end{figure} 
\ \\
\setlength{\tabcolsep}{4.0mm}  
\begin{table}[h!]
\vspace{-0.6cm}
\caption{Level $2$-$3$ RS Pyramidal Neuron Membrane Ion Channel Definitions} 
\centering
\begin{tabular}{l l l c c r}
\toprule \\[-2.5ex]
\multicolumn{1}{l}{Ion} & \multicolumn{2}{l}{Channel}  & \multicolumn{1}{l}{State Function}  &  \multicolumn{1}{l}{Gate}  &  \multicolumn{1}{r}{Type} \\
\multicolumn{1}{l}{Current} & \multicolumn{1}{l}{i} & \multicolumn{1}{l}{name} & \multicolumn{1}{l}{$p_i=f_i(m_i,h_i)$}  &  \multicolumn{1}{l}{Name} & \\
\hline \\[-1.5ex]
\multirow{3}{*}{Na+} & $1$ &  \multicolumn{1}{l}{ Fast Na+}  & \multicolumn{1}{l}{$\quad m_1^3h_1$}  & \multicolumn{1}{l}{$~m_1$}  &  \multicolumn{1}{r}{Activation} \\
         &        &                  &                                                      &  \multicolumn{1}{l}{$~h_1$}  &  \multicolumn{1}{r}{Inactivation} \\
         & $2$ &  Slow Na+ & \multicolumn{1}{l}{$\quad m_2$}            & \multicolumn{1}{l}{$~m_2$}  &  \multicolumn{1}{r}{Activation} \\ [1.0ex]
 \hline \\[-1.5ex]
\multirow{9}{*}{K+} & $3$ & DR\footnotemark[1] & \multicolumn{1}{l}{$\quad m_3^4$} & \multicolumn{1}{l}{$~m_3$}  &  \multicolumn{1}{r}{Activation} \\
        & $4$ &  A-trans & \multicolumn{1}{l}{$\quad m_4^4h_4$}  & \multicolumn{1}{l}{$~m_4$}  &  \multicolumn{1}{r}{Activation} \\
        &        &                  &                                                                  &  \multicolumn{1}{l}{$~h_4$}  &  \multicolumn{1}{r}{Inactivation} \\
        & $5$ &  $2$          & \multicolumn{1}{l}{$\quad m_5h_5$}        & \multicolumn{1}{l}{$~m_5$}  &  \multicolumn{1}{r}{Activation} \\
        &        &                  &                                                                  &  \multicolumn{1}{l}{$~h_5$}  &  \multicolumn{1}{r}{Inactivation} \\
 & $6$ & M\footnotemark[2] & \multicolumn{1}{l}{$\quad m_6$} & \multicolumn{1}{l}{$~m_6$}  &  \multicolumn{1}{r}{Activation} \\
 & $7$ & AHP\footnotemark[3] & \multicolumn{1}{l}{$\quad m_7$}    & \multicolumn{1}{l}{$~m_7$}  &  \multicolumn{1}{r}{Ca$^{2+}$ dep} \\
        & $8$ &  C           & \multicolumn{1}{l}{$\quad m_8\Gamma()$}& \multicolumn{1}{l}{$~m_8$}  &  \multicolumn{1}{r}{Activation} \\
        &        &                &                                                               & \multicolumn{1}{l}{$~\Gamma$}  &  \multicolumn{1}{r}{Ca$^{2+}$ dep} \\    
        & $9$ & AR\footnotemark[4]        & \multicolumn{1}{l}{$\quad m_9$}          & \multicolumn{1}{l}{$~m_9$}  &  \multicolumn{1}{r}{Activation} \\ [1.0ex]
 \hline  \\[-1.5ex]
\multirow{3}{*}{Ca$^{2+}$} & $10$ & Ca(T)\footnotemark[5] & \multicolumn{1}{l}{$\quad m_{10}^2h_{10}$} & 
\multicolumn{1}{l}{$~m_{10}$} & \multicolumn{1}{r}{Activation} \\
  &        &                  &                                                      &  \multicolumn{1}{l}{$~h_{10}$} & \multicolumn{1}{r}{Inactivation} \\% \ \\
  & $11$ & Ca(L)\footnotemark[6] & \multicolumn{1}{l}{$\quad m_{11}^2$} & \multicolumn{1}{l}{$~m_{11}$} & \multicolumn{1}{r}{Activation} \\
\bottomrule
\end{tabular}
\label{tab:p23rsChannels}
\end{table} 
\ \\
\footnotetext[1]{Delayed Rectifier}  \footnotetext[2]{Muscarinic}   \footnotetext[3]{After-Hyperpolarization} \footnotetext[4]{Anomalous Rectifier}
\footnotetext[5]{Both (T) for ``Transient'' and (L) for Low-Threshold have been used to describe this channel} 
\footnotetext[6]{Both (L) for ``Long'' and (H) for High-Threshold have been used to describe this channel}

%% For citations use: 
%%       \citet{<label>} ==> Jones et al. [21]
%%       \citep{<label>} ==> [21]
%%

%% If you have bibdatabase file and want bibtex to generate the
%% bibitems, please use
%%
%%  \bibliographystyle{elsarticle-num-names} 
%%  \bibliography{<your bibdatabase>}

%% else use the following coding to input the bibitems directly in the
%% TeX file.

\end{document}